%% file: main.tex
\documentclass[draft]{article} 
\usepackage[english]{babel}

\usepackage[letterpaper,top=2cm,bottom=2cm,left=2cm,right=2cm,marginparwidth=1.75cm]{geometry}

\usepackage{calrsfs}
\usepackage{amssymb}
\usepackage{amsmath}
\usepackage{amsthm}
\usepackage{caption}
\usepackage{subcaption}
\usepackage{multirow}
\usepackage[final]{graphicx}
\usepackage{amsfonts}
\usepackage{mathtools}
\usepackage{mathrsfs}
\usepackage{placeins}
\usepackage{algorithm}
\usepackage{algpseudocode}

\usepackage[mathscr]{euscript}
\usepackage[dvipsnames]{xcolor}
\usepackage[colorlinks=true, allcolors=OliveGreen]{hyperref}
\usepackage{filecontents}
\usepackage{tikz}
\usepackage{pgfplots}
\usetikzlibrary{external}
\usepackage{authblk}
\usepackage{todonotes}
\usepackage{array}
\newcolumntype{C}{>{\centering\arraybackslash}p{2cm}}

\usepackage{cleveref}
\usepackage{fixme}
\fxsetup{mode=multiuser}
\fxsetup{innerlayout=inline}
\fxusetheme{color}
\FXRegisterAuthor{fb}{envfb}{FB}

\newcommand{\R}{\mathbb{R}}
\newcommand{\N}{\mathbb{N}}
\renewcommand{\H}{\mathcal{H}}
\newcommand{\qoi}{\mathcal{Q}}
\newcommand{\SG}{\mathcal S}
\newcommand{\U}{\mathcal U}

\DeclareMathAlphabet{\mathcalligra}{T1}{calligra}{m}{n}

\title{Uncertainty Quantification for Fisher-Kolmogorov Equation on Graphs with Application to Patient-Specific Alzheimer Disease \footnote{\textbf{Funding}: PFA has been partially funded by the research grants PRIN2017 n. 201744KLJL funded by MUR. PFA and AMQ has been partially funded by the research grants PRIN2020 n. 20204LN5N5 funded by MUR. PFA has been partially supported by ICSC—Centro Nazionale di Ricerca in High Performance Computing, Big Data, and Quantum Computing funded by European Union—NextGenerationEU. FB is partially funded by “INdAM - GNCS Project”, codice CUP E53C22001930001. MC, FB, PFA and AMQ are members of INdAM-GNCS.}}

\author[1]{Mattia Corti}

\affil[1]{MOX-Dipartimento di Matematica, Politecnico di Milano, Piazza Leonardo da Vinci 32, Milan, 20133, Italy}

\author[1]{Francesca Bonizzoni}
\author[1]{Paola F. Antonietti}
\author[1,2]{Alfio M. Quarteroni}

\affil[2]{Institute of Mathematics, \'{E}cole Polytechnique F\'{e}d\'{e}rale de Lausanne, Station 8, Av. Piccard, Lausanne, CH-1015, Switzerland (Professor Emeritus)}

\begin{document}
\maketitle

\begin{abstract}
The Fisher-Kolmogorov equation is a diffusion-reaction PDE that is used to model the accumulation of prionic proteins, which are responsible for many different neurological disorders. Likely, the most important and studied misfolded protein in literature is the Amyloid-$\beta$, responsible for the onset of Alzheimer disease. Starting from medical images we construct a reduced-order model based on a graph brain connectome. The reaction coefficient of the proteins is modelled as a stochastic random field, taking into account all the many different underlying physical processes, which can hardly be measured. Its probability distribution is inferred by means of the Monte Carlo Markov Chain method applied to clinical data. The resulting model is patient-specific and can be employed for predicting the disease's future development. Forward uncertainty quantification techniques (Monte Carlo and sparse grid stochastic collocation) are applied with the aim of quantifying the impact of the variability of the reaction coefficient on the progression of protein accumulation within the next 20 years.
\end{abstract}

\section{Introduction}

Neurodegenerative disorders are neurological diseases that yield progressive damage to the neuronal tissue, compromising the communication between cells and then the entire cell structure. Clinical evidence suggests that a mechanism of spreading toxic misfolded prionic proteins is at the basis of neurodegenerative disorders \cite{walker_neurodegenerative_2015}. One of the pathologies in this group is the most common form of dementia nowadays, namely Alzheimer disease (AD) \cite{scheltens_Alzheimers_2021}. In AD, the pathological onset is connected to the accumulation of two proteins, namely the Amyloid$-\beta$ and the tau protein. Their accumulation causes damage with destruction of the synapses that mediate memory and cognition, due to the inability of living neurons to maintain functional axons and dendrites or by their death \cite{bloom2014amyloid}. 
\par
From a clinical point of view, the control of the disease progression is made both by performing cognition exams and by controlling the value of some critical biomarkers. Amyloid-$\beta$ is the earliest hallmark of AD, for this reason, positron emission tomography (PET) plays an important role in the clinical practice, to evaluate the concentration of the misfolded proteins. However, PET scanners are not widely available and are considered an invasive and relatively expensive clinical exam \cite{vanoostveenImagingTechniquesAlzheimer2021}, for this reason, in many cases PET is performed only for the diagnosis and not to control the progression of the disease \cite{hampelCoreCandidateNeurochemical2008}. 
\par
In this context, the numerical simulation of the misfolded protein concentrations can describe the chemical nature of this type of disease. Several mathematical models for prion dynamics have been proposed in the last few years based on partial differential equations (PDEs). Oligomer coagulation and fragmentation, which are the chemical processes to which Amyloid-$\beta$ undergoes, are typically described by means of the Smoluchowski equation \cite{franchi-lorenzani} and its suitable simplifications, such as the Fisher-Kolmogorov (FK) equation (also known as Fisher-KPP) \cite{weickenmeierPhysicsbasedModelExplains2019, fornariPrionlikeSpreadingAlzheimer2019}.
\par
Simulations of neurodegeneration on a complete brain geometry are mostly based on the FK model and on finite element discretizations \cite{weickenmeierPhysicsbasedModelExplains2019,schaferInterplayBiochemicalBiomechanical2019}. Recently, in \cite{corti:FK,bonizzoni_structure-preserving_2020} the authors proposed numerical methods based on discontinuous Galerkin methods for the resolution of FK equation. More complex models, such as the Smoluchowsky one, are at the moment unfeasible on three-dimensional meshes, due to the prohibitively computational complexity. For this reason, suitable simplified models have been proposed, based on graphs that represent the principal connections between different brain regions, namely the connectograms. These networks can be derived starting from Diffusion Weighted Images (DWI) and then used as a geometrically reduced order model to compute a numerical solution of the PDEs \cite{fornariPrionlikeSpreadingAlzheimer2019,gorielyNeuronalOscillationsEvolving2020}, upon including some patient-specific information inside the graph structure itself. The same approach is applied in this work: starting from DWI images we derive a brain connectome describing our patient's brain that we use as a discrete space domain for the FK equation. 
\par
Concerning the simulation of the FK model on graphs, diffusion directions and magnitude are derived from medical imaging. However, no information about the parameter describing the protein reaction can be obtained from clinical exams. Indeed, due to the simplified form of the equation, the reaction parameter represents many different physical processes at the same time: production, misfolding, aggregation and clearance \cite{fornariPrionlikeSpreadingAlzheimer2019}. Due to this complexity, some works describe the parameter as a random variable \cite{schaferBayesianPhysicsBasedModeling2021}. In this work, our approach for the calibration of the parameter makes use of the Amyloid-$\beta$ concentrations derived from two subsequent PET medical images to estimate the probability distribution by means of an inverse uncertainty quantification (UQ) algorithm. 
\par
Bayesian inverse UQ algorithms in the context of neurodegenerative disorders has been used in \cite{schaferBayesianPhysicsBasedModeling2021}, to describe the diffusion and reaction of tau proteins, by means of the FK equation. In \cite{schaferCorrelatingTauPathology2022}, the application of the Bayesian framework is extended to a coupling between the tau protein accumulation and the induced atrophy. In this work, we employ the Markov Chain Monte Carlo (MCMC) method \cite{robert:MCMC,hastingsMonteCarloSampling1970} to calibrate the reaction parameter of the FK equation modelling the spreading of Amyloid-$\beta$ protein. As a result, we obtain a patient-specific model describing the disease development.
\par
With the aim of predicting the disease development and assessing the impact of the reaction coefficient variability on the random protein concentration, we apply forward UQ methods. In particular, we introduce the quantity of interest (QoI) being the spatial average of the protein concentration at a set of time instances ($t=5,10,15,20$ years), and we estimate its expectation and variance by means of the Monte Carlo (MC) \cite{haltonRetrospectiveProspectiveSurvey1970} as well as the sparse grid stochastic collocation (SC) method \cite{babuskaStochasticCollocationMethod2007,nobileSparseGridStochastic2008}. To the best of our knowledge, this is the first time that the SC method is employed to infer the evolution of neurodegenerative diseases.    
\par
The paper is organized as follows. In Section \ref{sec:model}, we introduce the FK mathematical model with stochastic reaction parameter and discuss its application to neurodegeneration. Moreover, we introduce the space discretization of the problem on the reconstructed graph. In Section~\ref{eq:UQ_for_FK} we recall numerical methods for inverse and forward UQ. Numerical results are presented in Section \ref{sec:numerical_results}: we introduce the numerical discretization and the PET imaging projection (Section \ref{sec:graph_reconstruction}), we show the results of the inverse UQ (Section \ref{sec:estimation_parameter_distributions}), and we discuss the results of the forward UQ simulations over 20 years (Section \ref{sec:future_prediction}). Finally, in Section \ref{sec:conclusions} we draw some conclusions and discuss further developments. 

\section{Mathematical model}
\label{sec:model}

In this section, we present the FK equation to describe the reaction and diffusion of Amyloid-$\beta$ proteins. The problem is dependent on time $t\in(0,T]$, with final time $T>0$, and space $\boldsymbol{x}\in\Omega\subset\R^d$ ($d=2,3$). A typical assumption in the literature is that there exists a  constant baseline concentration of healthy state proteins. The problem at hand can be described by only one relative concentration of the misfolded protein $c$. Indeed, the variable $c$ is rescaled in the interval $[0,1]$, where $0$ means the absence of misfolded proteins and $1$ is the high prevalence of them. For a detailed model derivation, we refer to \cite{weickenmeierPhysicsbasedModelExplains2019}.
\par
We assume that the reaction parameter $\alpha=\alpha(\boldsymbol{x},\boldsymbol{p})$, representing the local conversion rate of the proteins from the healthy to the misfolded state, depends on the parameter vector $\boldsymbol{p} = ({p}_1,...,{p}_N)$ randomly varying in the hyperrectangle $\Gamma = \Gamma_1 \times ... \times \Gamma_N \subseteq \R^N$. In particular, we take $\alpha(\boldsymbol{x},\boldsymbol{p})$ as follows:
\begin{equation}
    \label{eq:alpha}
    \alpha(\boldsymbol{x},\boldsymbol{p})
    = \sum_{j=1}^N p_j \chi_{\Omega_j}(\boldsymbol{x}),
\end{equation}
where $\left(\Omega_j\right)_{j=1}^N$ is a non-overlapping partition of the domain $\Omega$ and $p_j$ denotes the reaction parameter in the region $\Omega_j$, for all $j=1,\ldots,N$. As a consequence, the concentration $c$ depends on the space and time variables as well as on the parameter vector $\boldsymbol{p}$, namely, $c=c(\boldsymbol{x},t,\boldsymbol{p})$.
\par
The differential model reads as follows: find $c\colon \Omega\times[0,T]\times \Gamma\rightarrow \R$ such that
\begin{equation}
\begin{dcases}
     \dfrac{\partial c}{\partial t}(\boldsymbol{x},t,\boldsymbol{p}) =\nabla \cdot(\mathbf{D}(\boldsymbol{x}) \nabla\, c(\boldsymbol{x},t,\boldsymbol{p}) + \alpha(\boldsymbol{x},\boldsymbol{p})\,c(\boldsymbol{x},t,\boldsymbol{p})(1-c(\boldsymbol{x},t,\boldsymbol{p})) & (\boldsymbol{x},t,\boldsymbol{p})\in\Omega\times(0,T]\times\Gamma, \\[8pt]
   \left(\mathbf{D}(\boldsymbol{x}) \nabla\, c(\boldsymbol{x},t,\boldsymbol{p})\right)\cdot\boldsymbol{n}(\boldsymbol{x}) = 0 & (\boldsymbol{x},t,\boldsymbol{p})\in\partial\Omega\times(0,T]\times\Gamma, \\[8pt]
     c(\boldsymbol{x},0,\boldsymbol{p}) = c_0(\boldsymbol{x})& (\boldsymbol{x},\boldsymbol{p})\in\Omega\times\Gamma.
\end{dcases}
\label{eq:fk_strong}
\end{equation}
The reaction parameter is assumed to be stochastic since its value is highly patient-specific and depends on the production, misfolding and clearance processes \cite{fornariPrionlikeSpreadingAlzheimer2019,schaferBayesianPhysicsBasedModeling2021}. The diffusion tensor $\mathbf{D}=\mathbf{D}(\boldsymbol{x})$ denotes the spreading of the misfolded protein inside the domain, which represents the volume occupied by the whole brain parenchymal tissue. From a biological point of view, most of the diffusion happens in the axonal direction; for this reason, we derive a graph, which connects brain regions from DWI \cite{weickenmeierPhysicsbasedModelExplains2019}. This graph is typically known as brain connectome \cite{DSIStudioTractographySoftware}. We denote with $ c_0=c_0(\boldsymbol{x})$ the deterministic initial condition. Finally, at $\partial\Omega$, we assume homogeneous Neumann boundary condition, $\boldsymbol{n}=\boldsymbol{n}(\boldsymbol{x})$ being the unit normal to $\partial\Omega$.

\subsection{Discrete formulation}

\begin{figure}[t]
    \centering
    \includegraphics[width=0.8\textwidth]{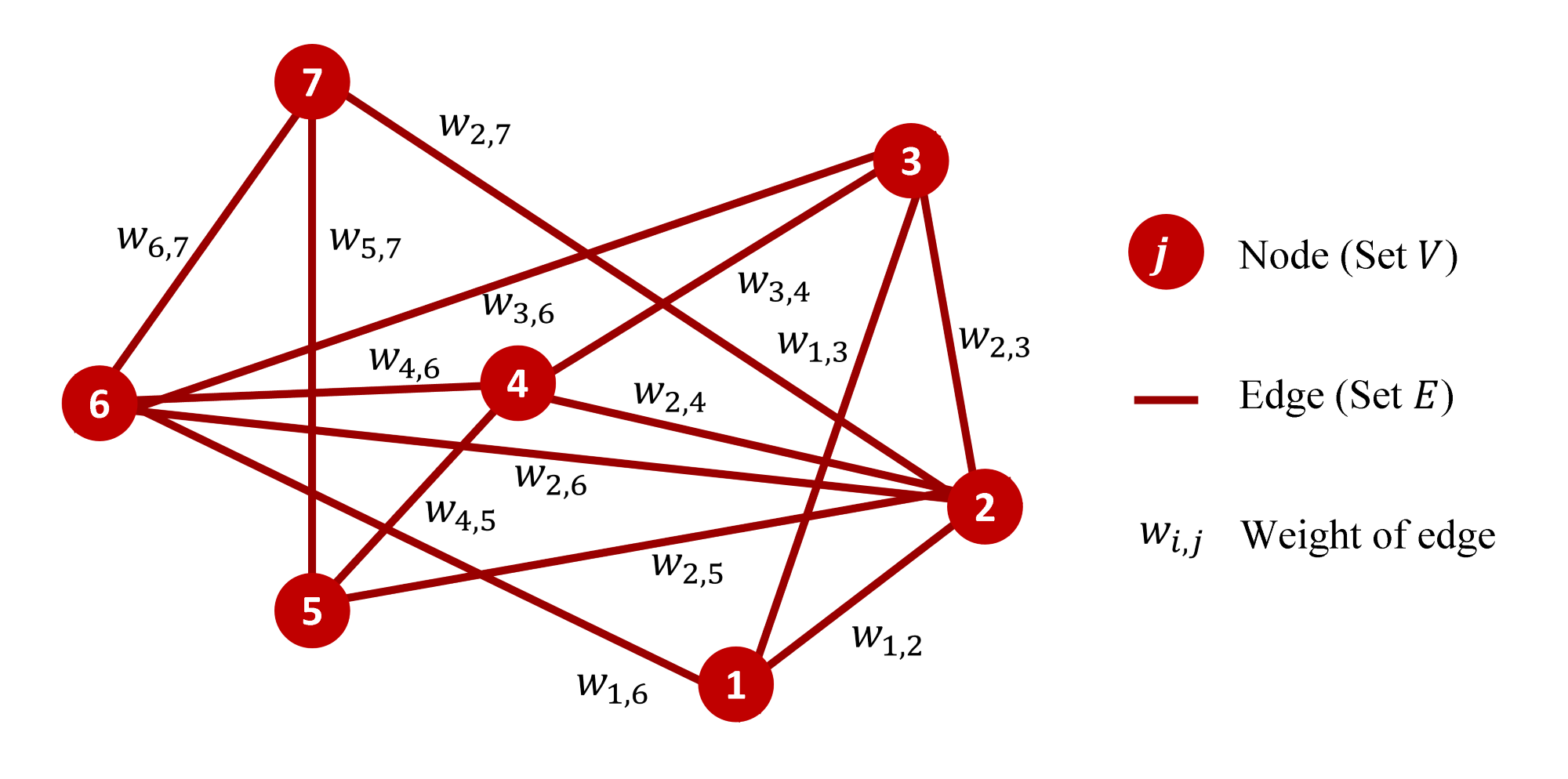}
    \caption{Schematic representation of a weighted undirected graph.}
    \label{fig:GraphImage}
\end{figure}

Due to the nature of the prionic diffusion, which is mostly directed as the axonal connections, we discretize Equation \eqref{eq:fk_strong} on the brain connectome graph, following \cite{betzelGenerativeModelsNetwork2017,fornariPrionlikeSpreadingAlzheimer2019}. More in details, we introduce the graph $\mathcal{G}=(V,E)$ representing a discretization of the domain $\Omega$, where $V = (\boldsymbol{x}_j)_{j=1}^M$ denotes the set of vertices (nodes) of the graph
and $E$ denotes the set of edges, as it can be observed in Figure \ref{fig:GraphImage}.
Applying the theory of discretization of PDEs on graphs \cite{smolaKernelsRegularizationGraphs2003} to the model problem \eqref{eq:fk_strong}, we find the following semi-discretized formulation 
\begin{equation}
\begin{dcases}
     \dfrac{\partial c_j}{\partial t}(t,\boldsymbol{p}) = -\sum_{i=1}^{M} \mathbf{L}_{ji}c_i(t,\boldsymbol{p}) + \alpha(\boldsymbol{x}_j,\boldsymbol{p})\,c_j(t,\boldsymbol{p})(1-c_j(t,\boldsymbol{p})) & (t,\boldsymbol{p})\in (0,T]\times\Gamma\quad j=1,...,M, \\[8pt]
     c_j(0,\boldsymbol{p}) = c_{0j} & \boldsymbol{p}\in\Gamma\quad j=1,...,M.
\end{dcases}
\label{eq:fk_semi}
\end{equation}
Equation \eqref{eq:fk_semi} is a system of networked ordinary differential equations (ODEs), where $\mathbf{L} = (\mathbf{L}_{ij})_{i,j=1}^M$ denotes the Laplacian matrix that discretizes the continuous Laplacian operator (see \cite{smolaKernelsRegularizationGraphs2003}), and $c_j(t,\boldsymbol{p}) \simeq c(\boldsymbol{x}_j,t,\boldsymbol{p})\colon [0,T]\times\Gamma\rightarrow\R$ is the approximated solution at the $j$-th node $\boldsymbol{x}_j$ of the graph $\mathcal{G}$. The solution evaluated at the set of nodes $\boldsymbol{x}_j\in V$ is collected into the vector $\boldsymbol{c} = (c_j)_{j=1}^M\colon [0,T]\times\Gamma\rightarrow\R^M$. In the same way, the values of the initial condition $c_0(\boldsymbol{x})$ at $V$ are collected into the vector $\boldsymbol{c}_0=(c_{0j})_{j=1}^M\in\R^M$, with $c_{0j}= c_0(\boldsymbol{x}_j)$ for all $j=1,\ldots,M$.

To discretize in time, we apply the Crank-Nicolson scheme, with a second-order extrapolation for a semi-implicit treatment of the nonlinear reactive term. In particular, let $\{t_\ell\}_{\ell=0}^{N_t}$ be the uniform partition of the time interval $[0,T]$ into $N_t$ intervals with length $dt=\frac{T}{N_t}$, namely, $0=t_0<t_1<...<t_{N_t}=T$ and $t_\ell=\frac{\ell T}{N_t}$ for $\ell=0,...,N_t$.
The fully discrete formulation of problem \eqref{eq:fk_strong} reads: give the initial conditions $\boldsymbol{c}_0$ and $\boldsymbol{c}_{-1}$, find $\boldsymbol{c}^{k+1} = \boldsymbol{c}^{k+1}(\boldsymbol{p})$, such that
\begin{equation}
    \begin{dcases}
    \dfrac{\boldsymbol{c}^{k+1} - \boldsymbol{c}^{k}}{\Delta t} = - \frac{1}{2} \mathbf{L}\left(\boldsymbol{c}^{k+1} +\boldsymbol{c}^{k}\right) + \boldsymbol{\alpha}\odot\left(\dfrac{1}{2}\boldsymbol{c}^{k+1}+\dfrac{1}{2}\boldsymbol{c}^{k}\right)\odot\left(\boldsymbol{1}-\left(\frac{3}{2}\boldsymbol{c}^{k}-\frac{1}{2}\boldsymbol{c}^{k-1}\right)\right) & k=0,...,N_t-1, \\[8pt]
    \boldsymbol{c}^0 = \boldsymbol{c}_{0}, \qquad \boldsymbol{c}^{-1} = \boldsymbol{c}_{-1},
    \end{dcases}
    \label{eq:fk_fully}
\end{equation}
where $\boldsymbol{\alpha}$ is the vector containing for each node of the graph $\mathcal{G}$, the corresponding value of the reaction parameter $\alpha$, $\boldsymbol{1}$ is a vector of ones and the symbol $\odot$ denotes the (component-by-component) Hadamard product.

\section{UQ methods for FK with stochastic reaction on graphs}
\label{eq:UQ_for_FK}

In this section, we describe inverse and forward uncertainty quantification algorithms, which we apply in Section~\ref{sec:numerical_results}. In particular, in Section~\ref{sec:inverse_UQ}, we introduce the MCMC method, and in Section~\ref{sec:forward_UQ}, we introduce the MC and the SC methods.

\subsection{Inverse Uncertainty Quantification}
\label{sec:inverse_UQ}

As detailed in Section \ref{sec:model}, the FK equation \eqref{eq:fk_strong} depends on the reaction parameter $\boldsymbol{\alpha}(\boldsymbol{x},\boldsymbol{p})$ defined in equation~\eqref{eq:alpha}. We model it as a stochastic quantity depending on the random vector $\boldsymbol{p}$, because its value encodes many different physical processes that cannot be directly estimated. However, due to the patient-specificity of the neurodegeneration, it is difficult to determine a distribution of the components $(p_1,\ldots,p_N)$ that can be used for every patient \cite{schaferBayesianPhysicsBasedModeling2021}. For this reason, in this work, we propose to estimate specific distributions of the parameters, starting from the PET medical images and using the MCMC method.
\par
Let us denote with $\qoi\in\R$ the QoI or output of the fully-discrete FK formulation \eqref{eq:fk_fully}. It can be seen as a function of the $N$ uncertain parameters, namely $\qoi\colon\Gamma\rightarrow\R$, even though generalizations to the vector-valued case are straightforward. Possible examples are the spatial average of the solution at the final time 
\begin{gather}
  \label{eq:qoi_space_average}
    \qoi(\boldsymbol{p})=\langle c\rangle(T,\boldsymbol{p}) = \frac{1}{M} \sum_{i=1}^M c_i(T,\boldsymbol{p}),
\end{gather} 
(see Section~\ref{sec:numerical_results}) or the solution at one specific node $x_{j^*}$ at the final time ($\qoi(\boldsymbol{p})=c_{j^*}(T,\boldsymbol{p}$)).
\par
In order to determine the distribution of the parameter vector $\boldsymbol{p}=(p_1,\ldots,p_N)\in \Gamma$, we want to construct a Markov chain, by extracting at every iteration an independent step from a Gaussian distribution $\boldsymbol{\delta}\sim\mathcal{N}(\boldsymbol{0},\widehat{\Sigma})$, where we assume that $\widehat{\Sigma}=\widehat{\sigma}^2 \mathbf{I}$, with $\widehat{\sigma}$ a scalar value that we assume constant for each parameter of the vector $\boldsymbol{p}$. The idea of the Metropolis-Hastings algorithm is to construct a posterior distribution of the parameter  $\pi_\mathrm{post}(\boldsymbol{p})$ starting from some assumptions on the prior distributions $\pi_\mathrm{pr}(\boldsymbol{p})$ of the parameters. In particular, in this work, we assume a uniform prior distribution of the parameters, with bounds $a_\ell$ and $b_\ell$ for the $\ell$-th component of the parameters vector $\boldsymbol{p}$.
\par
At each iteration, we choose to accept or not the new value of the chain, by using a Bernoulli distribution of parameter $\rho$, computed using the application of Bayes' theorem and defined by:
\begin{equation}
    \rho =\dfrac{\pi_\mathrm{post}(\boldsymbol{p}^*|\qoi^{\mathrm{PET}})}{\pi_\mathrm{post}(\boldsymbol{p}^{(i-1)}|\qoi^{\mathrm{PET}})} = \dfrac{\pi(\qoi^\mathrm{PET}|\boldsymbol{p}^*)}{\pi(\qoi^\mathrm{PET}|\boldsymbol{p}^{(i)})}\cdot\dfrac{\pi_\mathrm{pr}(\boldsymbol{p}^*)}{\pi_\mathrm{pr}(\boldsymbol{p}^{(i)})},
\end{equation}
where $\pi(\qoi|\boldsymbol{p})$ is the conditioned distribution of the quantity of interest, with respect to the parameters' realization. The computation of this quantity requires that we assume a model for the uncertainty in the measuring errors. In particular, we can assume that the errors are distributed as a vectorial Gaussian distribution $\mathcal{N}(\boldsymbol{0},\sigma^2\mathbf{I})$. In this way, we can obtain a formula for the likelihood distribution:
\begin{equation}
    \qoi^\mathrm{PET}|\boldsymbol{p} = \mathcal{N}(\qoi(\boldsymbol{p}),\sigma^2\mathbf{I}).
\end{equation}
The MCMC algorithm is reported in Algorithm \ref{alg:MCMC}. At the end of the iterations, we compute the vector of the mean $\mu_\ell$ and variance $\sigma^2_\ell$ of the realizations of the $p_\ell$ parameter. However, we neglect the first $\widetilde{M}$ realizations, to eliminate the dependency on the initial point of the chain. More information about MCMC methods can be found in \cite{robert:MCMC}.
\begin{algorithm}[t]
\caption{MCMC estimation of the distribution of $\boldsymbol{p}$}
\label{alg:MCMC}
\begin{algorithmic}
		\Require Initial value of input parameter $\boldsymbol{p}^{(0)}$
		\Ensure Samples of the distribution $\boldsymbol{p}^{(1:M)}$
		\For{$i=1:M$}
            \State Extract the white noise value $\boldsymbol{\delta}^*\sim\mathcal{N}(\boldsymbol{0},\Sigma_p)$
			\State Compute $\boldsymbol{p}^*=\boldsymbol{p}^{(i-1)}+\boldsymbol{\delta}^*$
            \State Compute $c^*=c^*(\boldsymbol{x},t,\boldsymbol{p}^*)$ by solving equation~\eqref{eq:fk_fully}
            \State Evaluate the QoI $\qoi^*$ at $c^*$
			\State Compute the acceptance rate $ \rho = \dfrac{\pi_\mathrm{post}(\boldsymbol{p}^*|\qoi^{\mathrm{PET}})}{\pi_\mathrm{post}(\boldsymbol{p}^{(i-1)}|\qoi^{\mathrm{PET}})}$
			\State  Extract from the Bernoulli the choice of accepting or not $y \sim Be(\rho)$
			\If{$ y=1 $}
				\State $\boldsymbol{p}^{(i)}=\boldsymbol{p}^*$
			\Else
				\State $\boldsymbol{p}^{(i)}=\boldsymbol{p}^{(i-1)}$
            \EndIf
		\EndFor
        \State Discard the first $\widetilde{M}$ realization of the Markov chain
        \State Compute the mean of the parameters' realizations: \[\mu_\ell = \frac{1}{M-\widetilde{M}}\displaystyle\sum_{j=\widetilde{M}+1}^M p^{(j)}_\ell\]
        \State Compute the variance of the parameters' realizations: \[\sigma_\ell^2 = \frac{1}{M-\tilde{M}-1}\displaystyle\sum_{j=\widetilde{M}+1}^M \left(\left(p^{(j)}_\ell\right)^2-\mu^2_\ell\right)\]
	\end{algorithmic}
 \end{algorithm}

\subsection{Forward Uncertainty Quantification}
\label{sec:forward_UQ}

As already discussed in the previous section, FK equation \eqref{eq:fk_strong} depends on the vector of parameters $\boldsymbol{p}=(p_1,\ldots,p_N)\in \Gamma$. The random distribution of each parameter $p_\ell$, $\ell=1,\ldots, N$, can be determined following the procedure presented in Section \ref{sec:inverse_UQ}. In particular, as an output of the inverse UQ problem, we have the set of pairs $\{(\mu_\ell,\sigma_\ell)\}_{\ell=1}^{N}$, where each $(\mu_\ell,\sigma_\ell)$ denotes the mean and standard deviation of the normal-distributed parameter $p_\ell\sim\mathcal N(\mu_\ell,\sigma_\ell)$. Note that each $p_\ell$ takes values in the entire real axes, namely, $\Gamma_\ell=\R$, so that $\Gamma=\R^N$. Moreover, we denote with $\rho\colon\Gamma\rightarrow\R_+$ the joint probability density function (pdf) of $\boldsymbol{p}$. If the entries of $\boldsymbol{p}$ are independent, then there holds $\rho=\prod_{\ell=1}^{N} \rho_\ell$, where $\rho_\ell\colon\Gamma_\ell\rightarrow\R_+$ denotes the pdf of a Gaussian random variable $\rho_\ell(y)=\frac{1}{\sqrt{2\pi}\sigma_\ell}\textrm{exp}\left(\frac{y-\mu_\ell}{2\sigma^2_\ell}\right)$.
\par
In this framework, we are interested in quantifying how the uncertainty in the input parameter vector $\boldsymbol{p}$ affects the QoI, by computing its statistical moments, like the expectation and variance, respectively given by
\begin{equation}
    \begin{aligned}
    \mathbb E[\qoi]& = \int_\Gamma \qoi(\boldsymbol{p}) \rho(\boldsymbol{p}) \, d\boldsymbol{p},\\
    \mathbb V[\qoi]& = \int_\Gamma \big(\qoi(\boldsymbol{p})- \mathbb E[\qoi](\boldsymbol{p})\big)^2\rho(\boldsymbol{p})\, d\boldsymbol{p}.
    \end{aligned}
    \label{eq:QoI}
\end{equation}
In particular, in relation to our medical application, we aim at quantifying the impact of varying piecewise constant reaction coefficients of prionic proteins on the onset and progression of Alzheimer disease.

\subsubsection{MC method}

The easiest and most diffused approach to compute approximations to the statistical moments defined in Equation \eqref{eq:QoI} is the MC method. The MC algorithm (in its plain formulation) computes the estimators for the expectation and variance as follows:
\begin{gather}
        \begin{aligned}
            \mathbb E[\qoi]&\simeq \mu_Q = \frac{1}{Q} \sum_{q=1}^Q\qoi(c^{(q)}),\\
            \mathbb V[\qoi]&\simeq \sigma^2_Q = \frac{1}{Q-1} \sum_{q=1}^Q
                \left(\qoi(c^{(q)})-\mu_Q\right)^2,
        \end{aligned}
\end{gather}
see also Algorithm~\ref{alg:MC}. Standard analysis results show that the MC error is $O(Q^{-1/2})$, namely it decays as the number of samples $Q$ increases, with rate $-\frac{1}{2}$.

\begin{algorithm}
\caption{MC estimators for expectation and variance}
\label{alg:MC}
\begin{algorithmic}
\Require $Q$
\Ensure MC estimators $\mu_Q$, $\sigma_Q^2$
\State Set $\mu_Q=0$, $\nu_Q=0$
\State Generate $Q$ i.i.d. samples of the input parameter $\boldsymbol{p}$, $\{\boldsymbol{p^{(q)}}\}_{q=1}^Q$
\For{$q=1:Q$}
\State Compute $c^{(q)}=c^{(q)}(\boldsymbol{x},t,\boldsymbol{p}^{(q)})$ by solving equation~\eqref{eq:fk_fully}
\State Evaluate the QoI at $c^{(q)}$
\State $\mu_Q = \mu_Q + \qoi(c^{(q)})$
\State $\nu_Q = \nu_Q + \qoi(c^{(q)})^2$
\EndFor
\State $\sigma_Q^2 = \frac{1}{Q-1} \left(\nu_Q-\mu_Q^2\right)$
\State $\mu_Q = \frac{1}{Q}\mu_Q$
\end{algorithmic}
\end{algorithm}

\subsubsection{Sparse grid SC method}

A successful alternative approach to performing forward UQ is the SC method based on sparse grids. It can be informally described as an approximation technique, obtained as a linear combination of several tensor product interpolation operators $\Gamma$, each formed by a limited number of points. The underlying idea is the so-called \emph{sparsification principle}: while none of these tensor product interpolation operators will be very accurate, an overall good approximation can be recovered by carefully combining many of them. This comes at a much lower cost than that required by considering a full tensor product approach over the entire parameter space $\Gamma$. Indeed, such an approach would involve a number of grid points exponential in $N$ hence it would be affected by the so-called curse of dimensionality. Consequently, the full tensor product approach is unfeasible, even for moderately small $N$.

\subsubsection*{One-dimensional Stochastic Collocation}

The starting point for the construction of the sparse grid approximation to the QoI relies on the introduction of one-dimensional Lagrange interpolation operators $\U^{m(i_n)}\colon C^0(\Gamma_n)\rightarrow P_{m(i_n)-1}(\Gamma_n)$ for all $n=1,\ldots, N$, given by
\begin{gather}
    \label{eq:1D_interp}
    \U^{m(i_n)}[u](p_n)= 
    \sum_{j_n=1}^{m(i_n)} u(p_n^{(j_n)})\ell_j(p_n^{(j_n)})
    \quad \forall p_n\in\Gamma_n,
\end{gather}
where $\H_{n,m(i_n)}=\{p_n^{(j_n)}\}_{j_n}^{m(i_n)}$ and $\{\ell_{j_n}\}_{j_n}^{m(i_n)}$ denote the set of Lagrange abscissas (or collocation points) and Lagrange polynomials of degree $m(i_n)$, respectively. Here $i_n\geq 1$ is an integer denoting the level of approximation and $m(i_n)$ is the number of collocation points used to build the interpolation at level $i_n$. 
The choice of the collocation points depends on the distribution of the random parameter $p_n$. In the particular case where $p_n$ follows a Gaussian distribution (as in the present paper), possible choices are Gauss-Hermite points, Genz–Keister points or Leja points. 
As $i_n$ varies, we get a sequence of Lagrange interpolation operators $\{\U^{m(i_n)}\}_{i_n\geq 1}$. As a function of $i_n$, the (level-to-knot) function $m(i_n)$ has to fulfill: $m(0)=0$ (and $\U^0=0)$, $m(1)=1$ and $m(i_n)<m(i_n+1)$ for $i_n\geq 1$.

\subsubsection*{Smolyak sparse grids}

Let $\boldsymbol{i}=(i_1,\ldots,i_N)\in \N_+^N$ be the multi-index with entries representing the approximation level along each direction $n$, for $n=1,\ldots,N$. Given the $N$-dimensional set of collocation points 
\begin{gather*}
    \H_{m(\boldsymbol{i})}= \prod_{n=1}^N \H_{n,m(i_n)}
\end{gather*}
being the Cartesian product of the one-dimensional sets of collocation points $\H_{1,m(i_1)},\ldots, \H_{N,m(i_N)}$, we introduce the tensor product interpolation operator $\U^{m(\boldsymbol{i})}= \U^{m(i_1)}\otimes\cdots\otimes\U^{m(i_N)}$ given by
\begin{gather}
    \label{eq:ND:interp}
    \U^{m(\boldsymbol{i})}[u](\boldsymbol{p})
    = \sum_{j_1=1}^{m(i_1)}\cdots\sum_{j_N=1}^{m(i_N)}
    u(p_1^{(j_1)},\ldots,p_N^{(j_N)})
    (\ell_{j_1}\otimes\cdots\otimes\ell_{j_N})(\boldsymbol{p})
\end{gather}
for all $u\in C^0(\Gamma)$. Note that the set $\H_{m(\boldsymbol{i})}$ contains $\#\H_{m(\boldsymbol{i})}=\prod_{n=1}^N m(i_n)$ collocation points. Hence, the tensor product interpolation operator $\U^{m(\boldsymbol{i})}$ requires $\#\H_{m(\boldsymbol{i})}$ evaluations of $u$.
For example, when taking $u=\qoi$ as in equation~\eqref{eq:qoi_space_average} being the spatial average of the FK solution at final time $T$, we get 
\begin{gather*}
    \qoi(\boldsymbol{p})\sim 
    \U^{m(\boldsymbol{i})}[\qoi](\boldsymbol{p})
\end{gather*}

The Smolyak sparse grid approximation to the QoI $\qoi$ is expressed by means of the combination technique as follows:
\begin{gather}
    \label{eq:sparse_grid}
    \SG_{I}[\qoi](\boldsymbol{p})
    = \sum_{\boldsymbol{i}\in I} \gamma_{\boldsymbol{i}}\U^{m(\boldsymbol{i})}[\qoi](\boldsymbol{p}),
\end{gather}
where 
\begin{itemize}
    \item 
    $\{\gamma_{\boldsymbol{i}}\}_{\boldsymbol{i}\in I}$ are the combination technique coefficients, given by $\gamma_{\boldsymbol{i}}=\sum_{\substack{\boldsymbol{j} \in \{0,1\}^N \\ \boldsymbol{i}+\boldsymbol{j} \in I}} (-1)^{|\boldsymbol{j}|}$. Note that some $c_{\boldsymbol{i}}$ might be null, in which case $\U^{m(\boldsymbol{i})}$ does not
    enter in the final approximation;

    \item
    $I\subset\N_+^N$ is a set of multi-indices that specifies which tensor interpolation operators enter in the sparse grid construction. It should be chosen according to the sparsification principle, namely, whenever one entry of $\boldsymbol{i}$ is large, the others should be kept as small as possible. Moreover, it is required to be downward-closed, namely, if $\boldsymbol{i}\in I$, then $\boldsymbol{i}-\boldsymbol{e}_n\in I$, for all $n=1,\ldots,N$, where $\boldsymbol{e}_n$ denotes the vector with all zeros expect the $k$-th component, that is equal to 1.
\end{itemize}
Note that equation \eqref{eq:sparse_grid} requires to solve equation~\eqref{eq:fk_fully} at all the collocation points $\boldsymbol{p}$ in the sparse grid
\begin{equation*}
    \H_I = \bigcup_{\boldsymbol{i}\in I}\H_{m(\boldsymbol{i})}.  
\end{equation*} 
Moreover, it becomes operative the moment we specify the three basic ``ingredients'' of the sparse grid construction, namely, the set $I$, the level-to-knot function $m(\cdot)$ and the one-dimensional knots used along each direction $n=1,\ldots,N$. In Section \ref{sec:numerical_results} we will detail our choices.

\subsubsection*{Sparse grids for forward UQ}


Associated with the Smolyak sparse grid approximant~\eqref{eq:sparse_grid} we introduce the sparse grid quadrature formula
\begin{gather}
    \label{eq:sparse_quadr}
     Q_I[\qoi]= \sum_{\boldsymbol{i}\in I} c_{\boldsymbol{i}} Q_{\boldsymbol{i}}[\qoi],
\end{gather}
where $Q_{\boldsymbol{i}}$ denotes the tensor product quadrature rule
\begin{gather*}
    \int_{\Gamma}\qoi(\boldsymbol{p})\rho(\boldsymbol p)\, d\boldsymbol p
    \approx Q_{\boldsymbol{i}}[\qoi]\coloneqq
    \sum_{j_1=1}^{m(i_1)}\cdots\sum_{j_N=1}^{m(i_N)} \qoi(p_1^{(j_1)},\ldots,p_N^{(j_N)})
    q_{1}^{(j_1)}\cdots q_{N}^{(j_N)}
\end{gather*}
and $q_{n}^{(j_n)}= \int_{\Gamma_n} \ell_{j_n}(p_n) \rho_n(p_n) \, dp_n$ for all $n=1,\ldots,N$.

Applying~\eqref{eq:sparse_quadr} we sparse grid approximations to the expectation and variance of the QoI:
\begin{equation*}
    \mathbb{E}[\qoi] \approx Q_I[\qoi], \qquad \mathbb{V}[\qoi] \approx Q_I[\qoi^2] - \left(Q_I[\qoi]\right)^2.
\end{equation*}



\section{Numerical results}
\label{sec:numerical_results}

The numerical results presented in this section concern a patient-specific test case. In particular, we are considering a patient affected by Alzheimer disease, with a calibration phase which from the age of 61 years to 68 years. Indeed, the computational domain is reconstructed from medical images (see Section \ref{sec:graph_reconstruction}), and the distribution of the reaction parameters in the different regions of the brain is recovered by solving an inverse UQ problem (see Section \ref{sec:estimation_parameter_distributions}). As the last step of our numerical study, we perform a forward UQ analysis  (see Section \ref{sec:future_prediction}), with the aim of determining the evolution of the Amyloid-$\beta$ protein for the next twenty years, taking as initial condition the data provided from the corresponding PET image. We mention that all the medical images we use in this section are provided by the OASIS-3 database \cite{OASIS3}. In particular, the specific brain geometry considered in this work is displayed in Figure \ref{fig:BraintoGraph:Surf}. The numerical solver used for space discretization is based on a built-in MATLAB code, whereas for the SC method we employ the Sparse Grids Matlab kit \cite{piazzolaSparseGridsMatlab2022}. The time discretization step is set equal to $dt=2\times 10^{-1}$ in all the simulations.

\subsection{Graph reconstruction and PET projection}
\label{sec:graph_reconstruction}

\begin{figure}[t]
     \centering
     \begin{subfigure}[b]{0.24\textwidth}
         \centering
         \includegraphics[width=\textwidth]{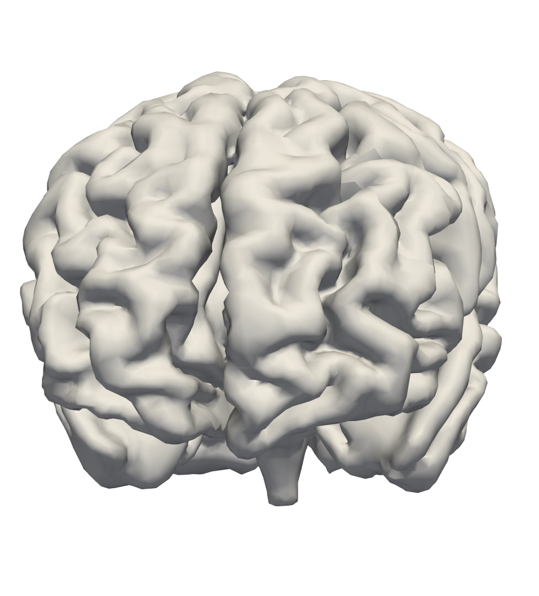}
         \caption{Brain surface}
         \label{fig:BraintoGraph:Surf}
     \end{subfigure}
     \hfill
     \begin{subfigure}[b]{0.24\textwidth}
         \centering
         \includegraphics[width=\textwidth]{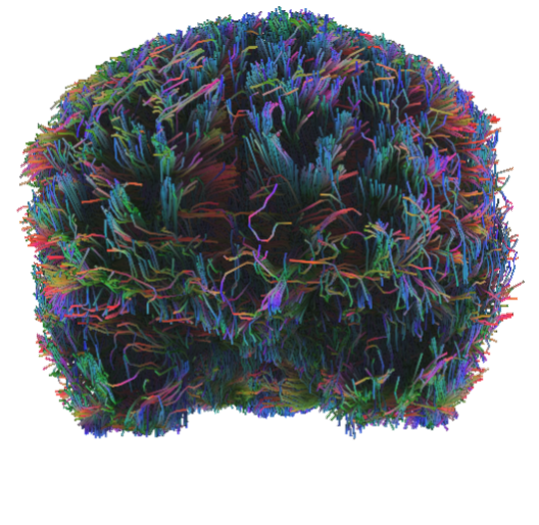}
         \caption{Tractography}
         \label{fig:BraintoGraph:Tracts}
    \end{subfigure}
    \hfill
    \begin{subfigure}[b]{0.24\textwidth}
         \centering
         \includegraphics[width=\textwidth]{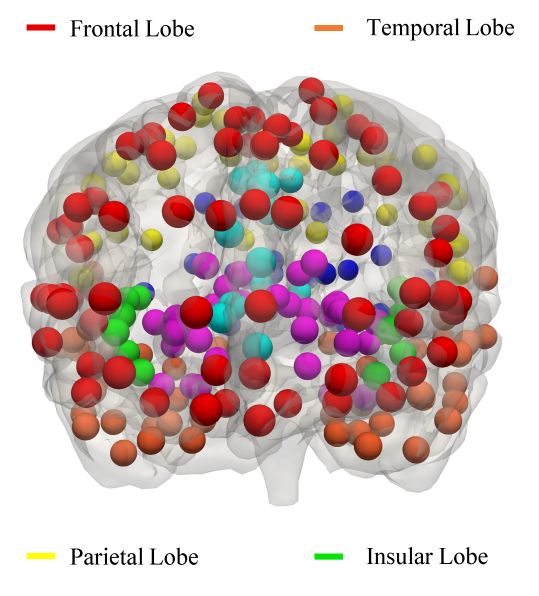}
         \caption{Nodes of the graph}
         \label{fig:BraintoGraph:Nodes}
    \end{subfigure}
    \hfill
    \begin{subfigure}[b]{0.24\textwidth}
         \centering
         \includegraphics[width=\textwidth]{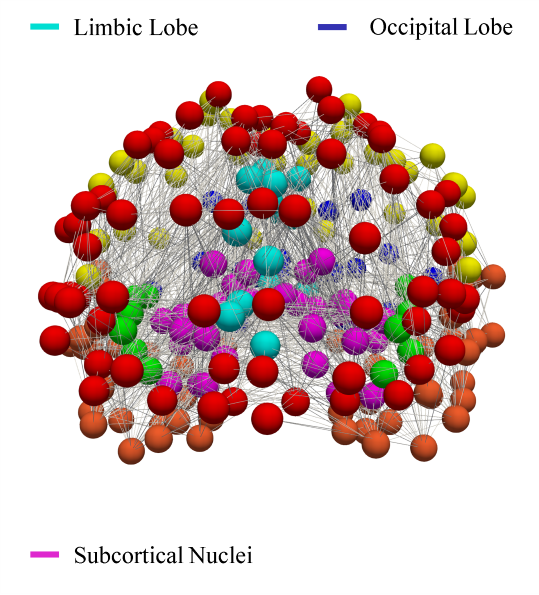}
         \caption{Graph}
         \label{fig:BraintoGraph:Graph}
     \end{subfigure}
        \caption{Graph discretization of the brain (each colour in Figures (c) and (d) denotes a different brain's region).}
        \label{fig:BraintoGraph}
\end{figure}

To generate the graph modelling the brain network we employ DWI, a Magnetic Resonance Imaging (MRI), which determines the Brownian motion of the water molecules; indeed the axonal nature of the neuronal connections creates directional paths for the motion. By using DSI studio \cite{DSIStudioTractographySoftware}, we can determine the brain's tractography, highlighting the axons' principal direction at every point of the brain \cite{yehShapeAnalysisHuman2020}, as we can observe in Figure \ref{fig:BraintoGraph:Tracts}.
\par
\begin{figure}[t!]
    \centering
    \includegraphics[width=\textwidth]{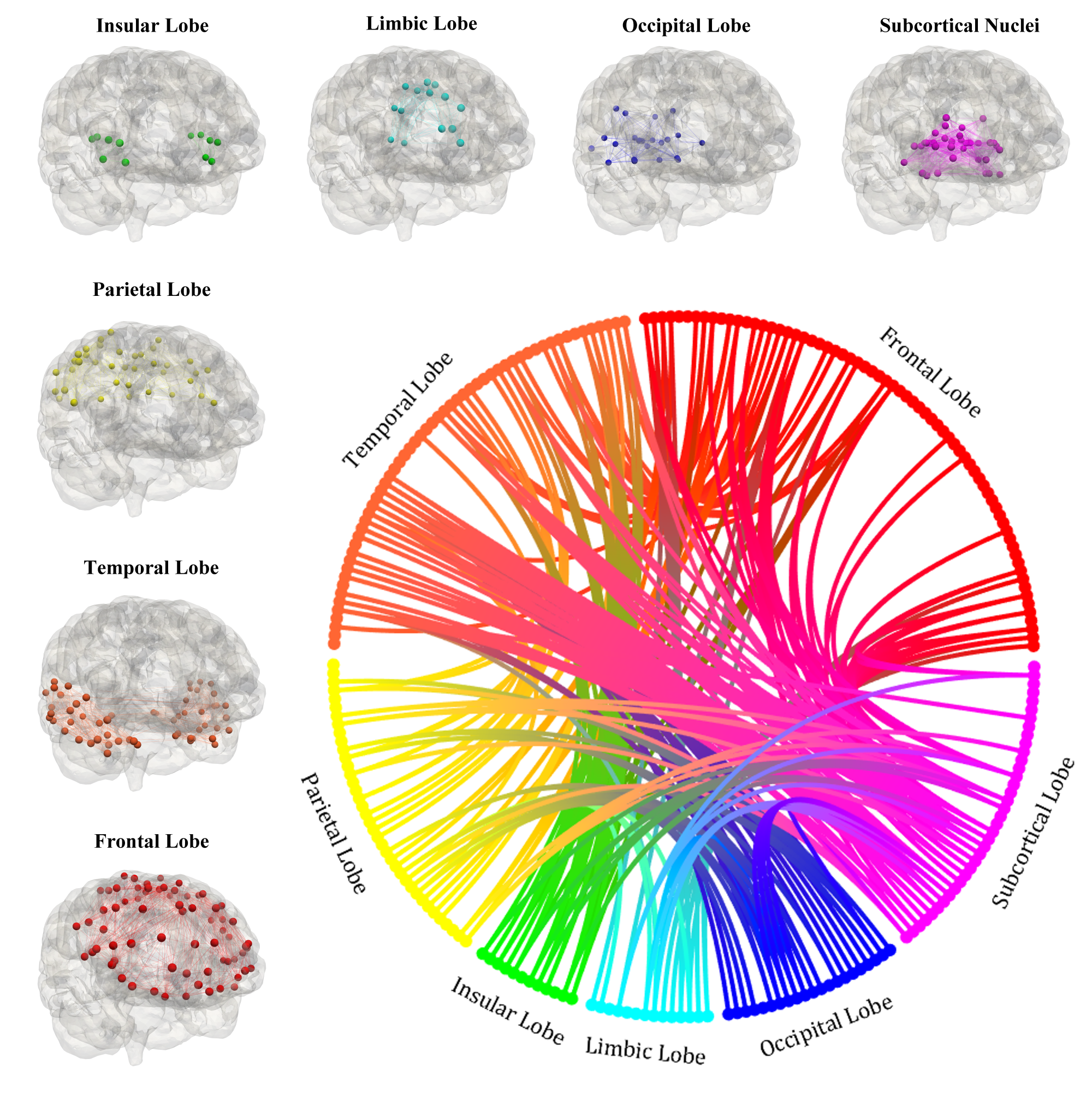}
    \caption{Local graphs of the seven regions of the brain and brain connectogram between different regions (excluding the connections with weight lower than $5\%$ of the principal one)}
    \label{fig:connectogram}
\end{figure}

The tractography is used to derive a weighted graph $\mathcal{G}$, where the nodes in $V$ are associated with a parcellation of the brain regions and the connecting edges $E$ are weighted taking into account both the number of connections and the length of the paths. In this work, we use parcellation based on the Brainnetome Atlas \cite{fanHumanBrainnetomeAtlas2016} and we group the nodes in $7$ regions: frontal lobe, temporal lobe, parietal lobe, insular lobe, limbic lobe, occipital lobe, and subcortical nuclei.
In Figure~\ref{fig:BraintoGraph:Nodes} we report the $246$ nodes of the graph coloured according to the corresponding regions, and the edges of the graph are visible in Figure \ref{fig:BraintoGraph:Graph}. 
We report the local graphs of each region in Figure \ref{fig:connectogram}. Moreover, we display the connectogram representing the principal connections between different regions. For readability, the connectogram shows only the principal connections. More specifically, we show the ones associated with weights that are larger than a threshold, which we fix to be equal to the $5\%$ of the maximum weight detected. We underline that the structure of the regions and the graph discretization have a large impact on the parameter calibration and then on the final output of the simulation, as it will be detailed in Section \ref{sec:estimation_parameter_distributions}. 

\par
\begin{figure}[t!]
\begin{subfigure}[b]{0.5\textwidth}
    \resizebox{\textwidth}{!}{\input{PETScatter.tikz}}
    \caption{Scatter plot of projected values from the PET in the nodes of the graph. The point color represents the region, the size is proportional to the volume of the node.}
    \label{fig:PETscatter}
\end{subfigure}
\hfill
\begin{subfigure}[b]{0.45\textwidth}
    \centering
    \includegraphics[width=0.85\textwidth]{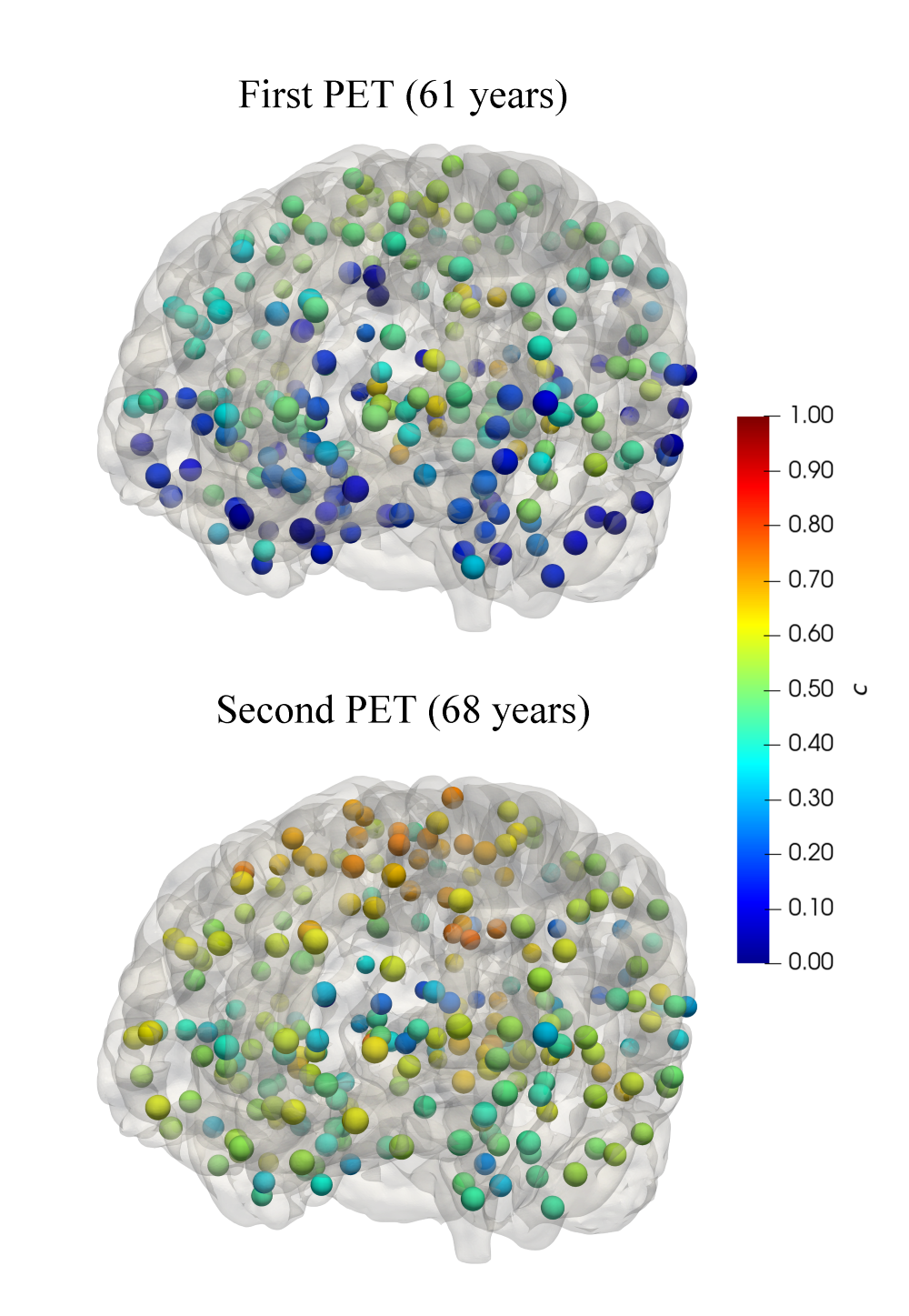}
    \caption{Projection of PET-PiB images on the reconstructed graph (first PET above and second one under). The color map refers to the level of proteins concentration.}
    \label{fig:PETprojected}
\end{subfigure}
\caption{Data extracted from PET medical images. The first PET is acquired at 61 years, and the second one is acquired at 68 years.}
\label{fig:PET}
\end{figure}

In order to estimate the reaction parameters as well as to set up the initial conditions for the (patient-specific) FK problem, we estimate the concentration $c$ of Amyloid-$\beta$ protein at two different time instances. The medical data used are Positron Emission Tomography images with Pittsburgh compound B (PET-PiB) \cite{vanoostveenImagingTechniquesAlzheimer2021}, which uses a radioligand able to identify the presence of Amyloid-$\beta$ plaques (more information on the acquisition techniques of the specific images used in this work can be found in \cite{OASIS3}). 
After having rescaled the values detected by the imaging so that they are between 0 and 1, we average the concentration on the parcellation reconstructed by DSI. Due to the fact that any region of the parcellation is associated with a single node in the brain connectome graph, we are obtaining the values of the concentration $c_j$ at any node $j$. The output of the process is the projection of the medical images (at 61 and 68 years of the patient) on the reconstructed graph: see Figure \ref{fig:PETprojected}. 

\par
By comparing the two images in Figure \ref{fig:PETprojected}, we notice that the concentration of Amyloid-$\beta$ has increased in the majority of the graph nodes. A more precise comparison is provided in Figure \ref{fig:PETscatter}, where we plot the concentration of the second PET (at 68 years) versus that of the first (at 61 years). Each bullet is coloured according to the region it belongs to (the same color code as in Figure \ref{fig:connectogram} is employed) and its size is proportional to the volume of the represented brain section in the parcellation. The scatter plot confirms the increase of the concentration in most of the nodes. However, it seems that a few nodes exhibit a decrease in the concentration of misfolded proteins. This is non-physical, we also observe that such a non-physical behaviour is mainly located in the subcortical lobe. We attribute this inconsistency to the atrophy of the region, which can be observed in the MRI images around the ventricles, and which is typical of AD \cite{vanoostveenImagingTechniquesAlzheimer2021}. Indeed, a reduction of the tissue can affect the projection of the medical data on the graph, due to the fact that the same pixel of PET image catches a larger part of the domain without neuronal tissue.

\subsection{Estimation of parameter distributions}
\label{sec:estimation_parameter_distributions}

In this section, we compute the distributions of the reaction parameters in the 7 regions of the graph by means of the inverse UQ process illustrated in Section~\ref{sec:inverse_UQ}. To better calibrate the distributions, we need to take into account the  existence of some outliers in the PET measurements, as explained in Section \ref{sec:graph_reconstruction}. When calibrating the reaction parameters, we will not take into account those nodes at which the concentration in the second PET is lower than in the first one of more than $10\%$. 
\par
According to~\eqref{eq:alpha}, the reaction coefficient $\boldsymbol{\alpha}$ is modelled as a piecewise Gaussian random field, namely, $\boldsymbol{\alpha}$ is a Gaussian random variable whenever restricted to any of the seven brain regions. In particular, we optimize the average concentration of the protein in each region of the graph as follows:
\begin{equation}
    \qoi_j = \langle c_h \rangle_j (T) = \sum_{k_j=1}^{M_j} c^{N_t}_{k_j} \dfrac{v_{k_j}}{v_\mathrm{tot}} \qquad j=1,...,7.
    \label{eq:Qj}
\end{equation}
In equation~\eqref{eq:Qj} we are computing a weighted average on the $j$-th region of the graph (containing $M_j$ nodes), for $j=1,\ldots,7$. The concentration at each node $c^{N_t}_{k_j}$ is weighted by the volume of the $k_j$-th node $v_{k_j}$ and divided by the total volume $v_\mathrm{tot}$. The corresponding index for each brain region is reported in Table \ref{tab:MCMCDistributions}.

\begin{table}[t!]
    \centering
    \begin{tabular}{l|c|C|C|C|C}
    \hline
    \multirow{2}{*}{\textbf{Brain Lobe (Index)}} & 
    \multirow{2}{*}{\textbf{Parameter}}
    & \multicolumn{2}{c|}{\textbf{Prior Distribution}} 
    & \multicolumn{2}{c}{\textbf{Estimated Distribution}}
    \\ \cline{3-6}
      & & $a_\ell$ & $b_\ell$ & $\mu_\ell$ & $\sigma^2_\ell$
    \\ \hline
    \textbf{Frontal Lobe} (1) & $p_1$
    & $-0.0700$ & $0.4300$
    & $ 0.1801$ & $0.0077$
    \\ \hline
    \textbf{Temporal Lobe} (2) & $p_2$
    & $-0.1100$ & $0.3900$
    & $0.1421$ & $0.0079$
    \\ \hline    
    \textbf{Parietal Lobe} (3) & $p_3$
    & $-0.1900$ & $0.3100$
    & $0.0627$ & $0.0060$
    \\ \hline    
    \textbf{Insular Lobe} (4) & $p_4$
    & $-0.1500$ & $0.3500$ 
    & $0.1005$ & $0.0070$ 
    \\ \hline
    \textbf{Limbic Lobe} (5)  & $p_5$
    & $-0.1200$ & $0.3800$ 
    & $0.1351$ & $0.0075$ 
    \\ \hline
    \textbf{Occipital Lobe} (6)  & $p_6$
    & $-0.1900$ & $0.3100$ 
    & $0.0545$ & $0.0086$ 
    \\ \hline
    \textbf{Subcortical Nuclei} (7)  & $p_7$
    & $-0.1500$ & $0.3500$ 
    & $0.1147$ & $0.0093$ 
    \\ \hline
    \end{tabular}
    \caption{Bounds of the prior uniform distributions of the parameters and mean and variance of the estimated distributions.}
    \label{tab:MCMCDistributions}
\end{table}
\par
The Metropolis-Hastings algorithm (see Algorithm~\ref{alg:MCMC}) is initialized by using a step with standard deviation $\widehat{\sigma}=10^{-2}$. Moreover, we adopt a likelihood with standard deviation $\sigma=0.1$. Finally, we choose the prior distribution for the 7 parameters as $p_\ell \sim U(a_\ell,b_\ell)$, where $a_\ell$ and $b_\ell$ are calibrated according to the available medical data (see Table \ref{tab:MCMCDistributions}).This mean value is also used to initialize the Markov chain, namely $\boldsymbol{p}^{(0)}=\frac{1}{2}(\boldsymbol{a}+\boldsymbol{b})$, where $\boldsymbol{a}=(a_\ell)_\ell$ and $\boldsymbol{b}=(b_\ell)_\ell$. To neglect the impact of the initial choice, we eliminate from the mean and variance computation the first $\widetilde{M}=10\,000$ steps out of the $M=100\,000$ we perform in this simulation.
\par
\begin{figure}[t!]
    \centering
    \includegraphics[width=\textwidth]{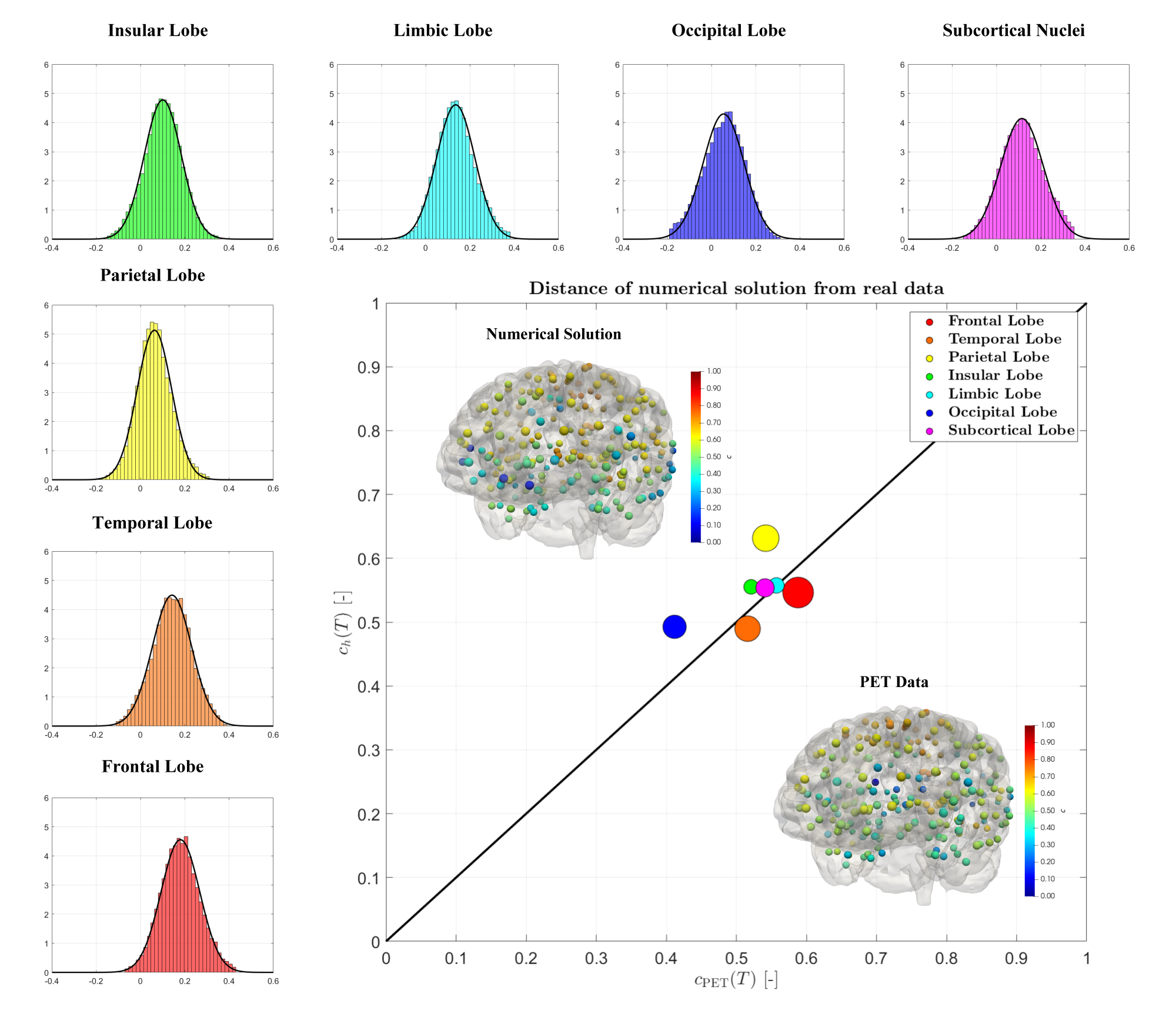}
    \caption{Results of the MCMC algorithm. Histogram and Gaussian distribution associated with each lobe of the brain and comparison between the medical data and the numerical results ($T=20$).}
    \label{fig:MCMCResult}
\end{figure}

The output values for $\mu_\ell$ and $\sigma^2_\ell$ are reported in Table \ref{tab:MCMCDistributions}. We can observe that the mean values feature high variability inside the different regions, and this result justifies the choice of using a non-homogeneous spatial description of the corresponding parameter. Moreover, in Figure \ref{fig:MCMCResult} we display the 7 different estimated distributions with the histograms associated with the Markov chain realizations.  We can observe coherence between the result of the MCMC algorithm and the resulting Gaussian probability distributions, that we use in the following section. Finally, in Figure \ref{fig:MCMCResult} we report the comparison between the $\qoi$ computed starting from the numerical solution, obtained using the mean values of the parameters and the corresponding $\qoi$ derived from the medical images, neglecting the outliers as explained in Section \ref{sec:graph_reconstruction}. In the results reported in Figure \ref{fig:MCMCResult} the size of the bullets is proportional to the volume of the lobes. The accuracy of the estimating algorithm for the parameters can be measured by the distance from the bisecting line, in which the numerical predicted value would be exactly equal to the medical data from PET images. Finally, we report the numerical solution on the graph and the medical data at the final time. In both cases, we observe a higher concentration of misfolded proteins inside the parietal and frontal lobes, coherently to what we expect from the medical literature \cite{hampelAmyloidvPathwayAlzheimer2021}. This shows that our numerical method is able to reproduce the disease progression.

\subsection{Prediction of disease development}
\label{sec:future_prediction}

In this section, we simulate future scenarios of the development of the Amyloid-$\beta$ concentration inside the brain, by following the two approaches presented in Section \ref{sec:forward_UQ}, namely MC and SC. 

In this simulation, we fix the final time at $T=20\,\mathrm{years}$ and we discretize the time interval $[0,T]$ with a timestep $dt=0.02\,\mathrm{years}$ (we denote the set of all considered time instances as $\mathcal I_T$). As QoI we choose the spatial average of the concentration inside the seven lobes of the brain graph, evaluated at all time instances $t\in\mathcal I_T$. In particular, at every $t\in\mathcal I_T$, we have a vectorial QoI $\qoi=(\qoi_1,\ldots,\qoi_7)$, where the $j$-th component is $\qoi_j = \langle c_h \rangle_j(t)$ for $j=1,...,7$. 
We aim at approximating the expected value and the variance of each of the QoI.

\par
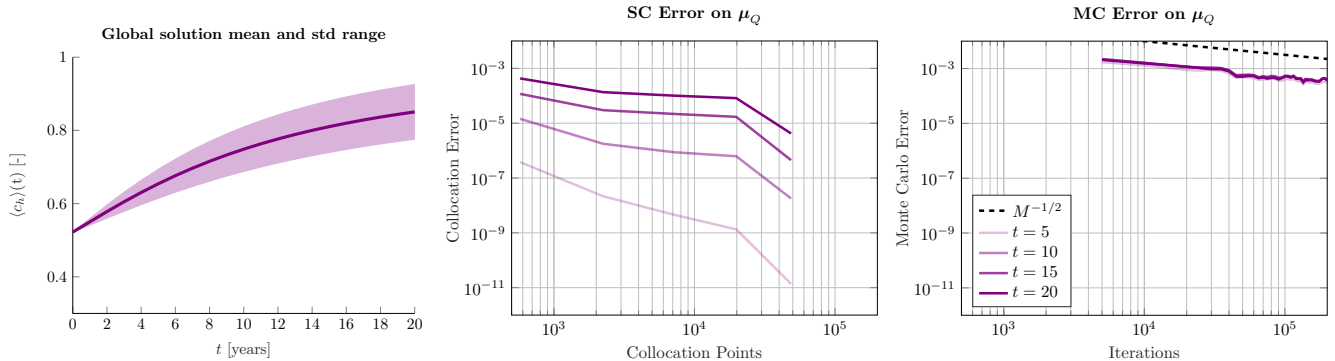
\begin{figure}[t!]
\begin{subfigure}[b]{0.32\textwidth}
    \resizebox{\textwidth}{!}{\input{Solution_Global.tikz}}
    \caption{Result of the simulation in terms of average over the whole brain.}
    \label{fig:Solution_Global}
\end{subfigure}
\hfill
\begin{subfigure}[b]{0.33\textwidth}
    \resizebox{\textwidth}{!}{\input{Coll_Err_Global.tikz}}
    \caption{Stochastic collocation errors at different times snapshots.}
    \label{fig:SCE_Global}
\end{subfigure}
\hfill
\begin{subfigure}[b]{0.33\textwidth}
     \resizebox{\textwidth}{!}{\input{MC_Err_Global.tikz}}
    \caption{MC errors versus $M$ at different times snapshots, (convergence rate $-1/2$).}
    \label{fig:MCE_Global}
\end{subfigure}
\caption{Numerical results of the SC and MC convergence tests for the average solution over the whole graph.}
\label{fig:GlobalResult}
\end{figure}

The reference expectation and variance are computed using the SC method on a Smolyak sparse grid of level 9 with Leja points (in total, there are $224\,143$ collocation nodes). The reference expectation and variance are denoted by $\boldsymbol{\mu}_Q$ and $\boldsymbol{\sigma}_Q^2$, respectively, and they are both vectors of length $7$: their different components correspond to the regions of the brain graph.
In Figure \ref{fig:Solution_Global}, we report the computed results in terms of the volume-weighted mean value over the different lobes $\langle\boldsymbol{\mu}_Q\rangle$ and a range of standard deviation $[\langle\boldsymbol\mu_Q\rangle-\langle\boldsymbol\sigma_Q\rangle;\langle\boldsymbol\mu_Q\rangle+\langle\boldsymbol\sigma_Q\rangle]$. We observe a global increase of prionic concentration in time, coherently to what expected by the disease progression \cite{hampelAmyloidvPathwayAlzheimer2021}. Moreover, the stochastic description yields an increasing uncertainty over the years, which allows for obtaining a range of future possible scenarios after many years.
\par
We perform a convergence analysis, with the aim of assessing the approximation properties of Smolyak sparse grids with increasing level $w=3,4,5,6,7,8$. 
The errors are reported in Table \ref{tab:SCErrors}. We can notice that both for $\langle\boldsymbol{\mu}_Q\rangle$ and $\langle\boldsymbol{\sigma}^2_Q\rangle$, we have a monotone decrease of the error as the number of collocation points increases. 
We can also notice, that the accuracy of the approximation deteriorates as time advances. Indeed, it is reasonable that it is necessary to use a larger number of points, to approximate phenomena with larger variance, such as the concentration after 15-20 years. 
The convergence results for the SC method are also displayed in Figures \ref{fig:SCE_Global} and \ref{fig:MCE_Global}, where the errors at time instances $t=5,10,15,20$ are reported versus the number of collocation points. For the sake of comparison, we also report the convergence of the MC error versus the number of MC samples (note that the MC error decays with the predicted rate $M^{-1/2}$, $M$ being the number of samples). 
We notice that, with a comparable number of evaluations of the differential problem, the SC method provides a more accurate solution. This is particularly pronounced for smaller time instances. 
\par
\begin{table}[t!]
    \centering
    \begin{tabular}{c|c|c|c|c|c|c}
    \hline
    \textbf{Grid level} & $\mathbf{3}$ & $\mathbf{4}$ & $\mathbf{5}$ & $\mathbf{6}$ & $\mathbf{7}$ & $\mathbf{8}$ \\ \hline
    \textbf{Collocation Points} & $375$ & $2\,241$ & $7\,183$ & $19\,825$ & $48\,639$ & $108\,545$  \\ \hline
    \multicolumn{7}{c}{\textbf{Collocation errors on approximation of }$\langle\boldsymbol{\mu}_Q\rangle$} \\ \hline
    
    Time: $t=5$ & $3.81\times10^{-7}$ & $2.15\times10^{-8}$ & $4.54\times10^{-9}$ & $1.34\times10^{-9}$ &	$1.37\times10^{-11}$ & $3.82\times 10^{-12}$ \\ \hline

    Time: $t=10$ & $1.44\times10^{-5}$ & $1.77\times10^{-6}$ & $8.69\times10^{-7}$ & $6.23\times10^{-7}$ &	$1.79\times10^{-8}$ & $1.61\times 10^{-8}$ \\ \hline

    Time: $t=15$ & $1.17\times10^{-4}$ & $2.98\times10^{-5}$ & $2.17\times10^{-5}$ & $1.70\times10^{-5}$ &	$4.48\times10^{-7}$ & $4.55\times 10^{-7}$ \\ \hline

    Time: $t=20$ & $4.29\times10^{-4}$ & $1.35\times10^{-4}$ & $1.01\times10^{-4}$ & $8.21\times10^{-5}$ &	$4.22\times10^{-6}$ & $2.17\times 10^{-6}$ \\ \hline
    
    \multicolumn{7}{c}{\textbf{Collocation errors on approximation of  }$\langle\boldsymbol{\sigma}^2_Q\rangle$} \\ \hline

    Time: $t=5$ & $5.69\times10^{-5}$ & $4.20\times10^{-6}$ & $1.26\times10^{-6}$ & $5.34\times10^{-7}$ &	$6.42\times10^{-9}$ & $3.37\times 10^{-9}$ \\ \hline
 
    Time: $t=10$ & $2.24\times10^{-3}$ & $4.32\times10^{-4}$ & $2.36\times10^{-4}$ & $1.50\times10^{-4}$ &	$4.43\times10^{-6}$ & $2.03\times 10^{-6}$ \\ \hline

    Time: $t=15$ & $8.06\times10^{-3}$ & $2.21\times10^{-3}$ & $1.40\times10^{-3}$ & $1.03\times10^{-3}$ &	$9.59\times10^{-5}$ & $2.78\times 10^{-5}$ \\ \hline

    Time: $t=20$ & $1.62\times10^{-2}$ & $5.56\times10^{-3}$ & $4.68\times10^{-3}$ & $3.94\times10^{-3}$ &	$4.99\times10^{-4}$ & $1.56\times 10^{-4}$ \\ \hline
    \end{tabular}
    \caption{Bounds of the prior uniform distributions of the parameters and expectation and variance of the estimated distributions.}
    \label{tab:SCErrors}
\end{table}

\begin{table}[t!]
    \centering
    \begin{tabular}{l|c|c|c|c|c|c|c|c}
    \hline
    \multirow{2}{*}{\textbf{Brain Lobe (Index)}} 
    & \multicolumn{2}{c|}{$t=5$} 
    & \multicolumn{2}{c|}{$t=10$}
    & \multicolumn{2}{c|}{$t=15$}
    & \multicolumn{2}{c}{$t=20$}
    \\ \cline{2-9}
      & $\mu_Q$ & $\sigma_Q^2$ & $\mu_Q$ & $\sigma^2_Q$ & $\mu_Q$ & $\sigma^2_Q$ & $\mu_Q$ & $\sigma^2_Q$
    \\ \hline
    \textbf{Frontal Lobe} (1)
    & $0.7527$ & $0.0049$
    & $0.8486$ & $0.0078$
    & $0.9000$ & $0.0077$
    & $0.9289$ & $0.0069$
    \\ \hline
    \textbf{Temporal Lobe} (2)
    & $0.6515$ & $0.0063$
    & $0.7598$ & $0.0117$
    & $0.8277$ & $0.0128$
    & $0.8699$ & $0.0122$
    \\ \hline    
    \textbf{Parietal Lobe} (3)
    & $0.6018$ & $0.0060$
    & $0.6729$ & $0.0156$
    & $0.7303$ & $0.0226$
    & $0.7738$ & $0.0265$
    \\ \hline    
    \textbf{Insular Lobe} (4)
    & $0.6417$ & $0.0061$
    & $0.7405$ & $0.0123$
    & $0.8096$ & $0.0142$ 
    & $0.8558$ & $0.1038$ 
    \\ \hline
    \textbf{Limbic Lobe} (5)
    & $0.6914$ & $0.0051$
    & $0.7892$ & $0.0095$
    & $0.8514$ & $0.0102$ 
    & $0.8905$ & $0.0093$ 
    \\ \hline
    \textbf{Occipital Lobe} (6)
    & $0.5121$ & $0.0073$
    & $0.6043$ & $0.0194$
    & $0.6782$ & $0.0275$ 
    & $0.7336$ & $0.0314$ 
    \\ \hline
    \textbf{Subcortical Nuclei} (7)
    & $0.5822$ & $0.0067$
    & $0.7066$ & $0.0147$
    & $0.7891$ & $0.0176$ 
    & $0.8413$ & $0.0175$ 
    \\ \hline
    \end{tabular}
    \caption{Expectation and variance of the solution distributions inside the brain regions.}
    \label{tab:FinalDistributions}
\end{table}

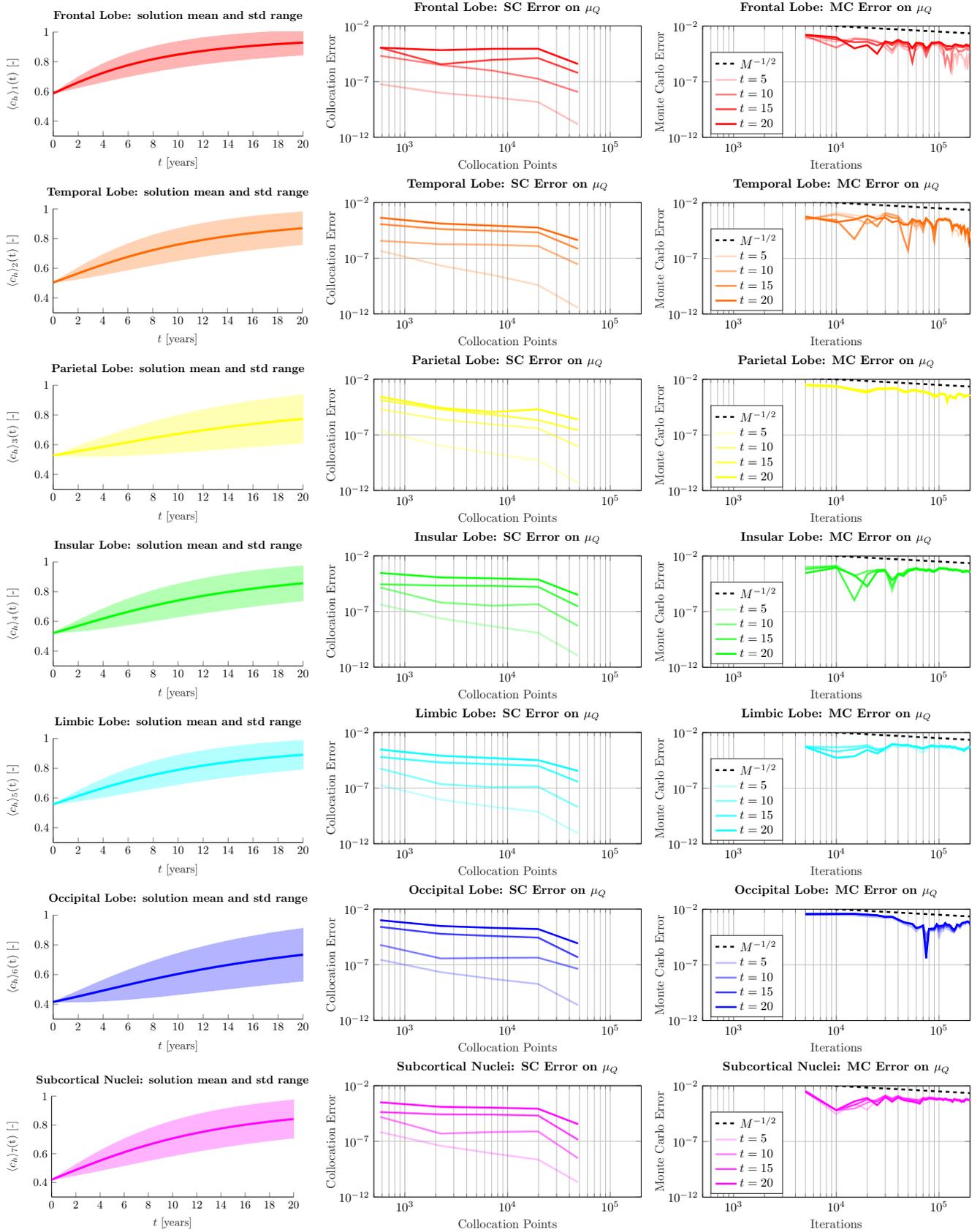
\begin{figure}[!]
\begin{subfigure}[b]{0.32\textwidth}
    \resizebox{\textwidth}{!}{\input{Solution_1_Frontal.tikz}}
\end{subfigure}
\hfill
\begin{subfigure}[b]{0.33\textwidth}
    \resizebox{\textwidth}{!}{\input{Coll_Err_1_Frontal.tikz}}
\end{subfigure}
\hfill
\begin{subfigure}[b]{0.33\textwidth}
     \resizebox{\textwidth}{!}{\input{MC_Err_1_Frontal.tikz}}
\end{subfigure}
\begin{subfigure}[b]{0.32\textwidth}
    \resizebox{\textwidth}{!}{\input{Solution_2_Temporal.tikz}}
\end{subfigure}
\hfill
\begin{subfigure}[b]{0.33\textwidth}
    \resizebox{\textwidth}{!}{\input{Coll_Err_2_Temporal.tikz}}
\end{subfigure}
\hfill
\begin{subfigure}[b]{0.33\textwidth}
     \resizebox{\textwidth}{!}{\input{MC_Err_2_Temporal.tikz}}
\end{subfigure}
\begin{subfigure}[b]{0.32\textwidth}
    \resizebox{\textwidth}{!}{\input{Solution_3_Parietal.tikz}}
\end{subfigure}
\hfill
\begin{subfigure}[b]{0.33\textwidth}
    \resizebox{\textwidth}{!}{\input{Coll_Err_3_Parietal.tikz}}
\end{subfigure}
\hfill
\begin{subfigure}[b]{0.33\textwidth}
     \resizebox{\textwidth}{!}{\input{MC_Err_3_Parietal.tikz}}
\end{subfigure}
\begin{subfigure}[b]{0.32\textwidth}
    \resizebox{\textwidth}{!}{\input{Solution_4_Insular.tikz}}
\end{subfigure}
\hfill
\begin{subfigure}[b]{0.33\textwidth}
    \resizebox{\textwidth}{!}{\input{Coll_Err_4_Insular.tikz}}
\end{subfigure}
\hfill
\begin{subfigure}[b]{0.33\textwidth}
     \resizebox{\textwidth}{!}{\input{MC_Err_4_Insular.tikz}}
\end{subfigure}
\begin{subfigure}[b]{0.32\textwidth}
    \resizebox{\textwidth}{!}{\input{Solution_5_Limbic.tikz}}
\end{subfigure}
\hfill
\begin{subfigure}[b]{0.33\textwidth}
    \resizebox{\textwidth}{!}{\input{Coll_Err_5_Limbic.tikz}}
\end{subfigure}
\hfill
\begin{subfigure}[b]{0.33\textwidth}
     \resizebox{\textwidth}{!}{\input{MC_Err_5_Limbic.tikz}}
\end{subfigure}
\begin{subfigure}[b]{0.32\textwidth}
    \resizebox{\textwidth}{!}{\input{Solution_6_Occipital.tikz}}
\end{subfigure}
\hfill
\begin{subfigure}[b]{0.33\textwidth}
    \resizebox{\textwidth}{!}{\input{Coll_Err_6_Occipital.tikz}}
\end{subfigure}
\hfill
\begin{subfigure}[b]{0.33\textwidth}
     \resizebox{\textwidth}{!}{\input{MC_Err_6_Occipital.tikz}}
\end{subfigure}
\begin{subfigure}[b]{0.32\textwidth}
    \resizebox{\textwidth}{!}{\input{Solution_7_Subcortical.tikz}}
\end{subfigure}
\hfill
\begin{subfigure}[b]{0.33\textwidth}
    \resizebox{\textwidth}{!}{\input{Coll_Err_7_Subcortical.tikz}}
\end{subfigure}
\hfill
\begin{subfigure}[b]{0.33\textwidth}
     \resizebox{\textwidth}{!}{\input{MC_Err_7_Subcortical.tikz}}
\end{subfigure}

\caption{Numerical results of the simulation in terms of the average on every single lobe (on the left), collocation errors at different times (in the middle), MC errors versus $M$ at different times snapshots $t=5,10,15,20$ (on the right).}
\label{fig:LobesResults}
\end{figure}
In Table \ref{tab:FinalDistributions} we report the expectation and variance of the averaged concentration inside every single region at time instances $t=5,10,15,20$. As we can notice, the expectation is increasing in time within each lobe, coherently with what we expect from the progression of the pathology. However, we can also notice that the variance is increasing only for the initial times, and when the expectation becomes larger than $0.8$ it starts decreasing due to the impossibility of the solution becoming larger than $1$, for the structure of FK equation \cite{salsa:EDP}. For example, concerning the frontal lobe, we can observe a reduction in the variance, associated with the high values of protein concentration. 
In Figure \ref{fig:LobesResults} 
each row corresponds to a different lobe of the brain. In the first column, we display the evolution of the average protein concentration with the relative standard deviation ranges; in the second column we report the error convergence for SC versus the number of collocation points; in the third column the report the MC error versus the number of samples.
We notice that also considering the regional average concentrations the SC method allows for reaching a more precise solution than MC, with a comparable number of evaluations of the differential problem. However, in this case, the SC error does not decrease monotonically (this is likely due to the nonlinearity of the PDE). Finally, we notice that, compared with the global solution in Figure \ref{fig:GlobalResult}, we have a larger variability of the solution in the different regions, particularly after 20 years.

\section{Conclusions}
\label{sec:conclusions}
In this work, we have proposed an approach to predict the progression of the accumulation of misfolded Amyloid$-\beta$ proteins by means of a Fisher-Kolmogorov model in Alzheimer’s disease. In order to model the reaction coefficients, which cannot be measured and whose values depend on many different physical processes, we adopted a stochastic model. The evaluation of the distributions has been made by means of an inverse uncertainty quantification algorithm and then applied to the prediction of the disease development in the next 20 years.
\par
To keep the computational cost affordable, a reduced-order model, based on a graph discretization of the brain connectome, is used. The derivation of the graph from medical images and the calibration of the parameter distributions by means of clinical data yielded a complete patient-specific model for the simulation of the progression of Alzheimer disease. In the numerical results, we presented a comparison in applying the Monte Carlo algorithm with respect to a stochastic collocation, for estimating proteins' concentration. The tests confirm that the application of collocation on sparse grids yields a more precise approximation of the final solution, in particular for smaller time frames (5-10 years), with the same number of evaluations of the PDE solution. 
\par
Some future developments of this work could concern the extension of the uncertainty quantification analysis to the three-dimensional models, for example, applying the study to discontinuous Galerkin methods on real brain geometries. Another possible development is using a white noise description for the weights of the graphs that are computed starting from medical images and could be affected by measuring errors. Finally, it can be interesting to perform a complete analysis on a large number of different patients, to highlight the differences in the computed values among them.
\section*{Acknowledgments}

The present research has been supported by MUR, grant Dipartimento di Eccellenza 2023-2027. The brain MRI images were provided by OASIS-3: Longitudinal Multimodal Neuroimaging: Principal Investigators: T. Benzinger, D. Marcus, J. Morris; NIH P30 AG066444, P50 AG00561, P30 NS09857781, P01 AG026276, P01 AG003991, R01 AG043434, UL1 TR000448, R01 EB009352.

\section*{Declaration of competing interests}
The authors declare that they have no known competing financial interests or personal relationships that could have appeared to influence the work reported in this article.

\bibliographystyle{ieeetr}
\bibliography{bibliography.bib}
\end{document}

%% file: PETScatter.tikz
\definecolor{frontal}{rgb}{1.00000,0.00000,0.00000}
\definecolor{temporal}{rgb}{1.00000,0.40000,0.20000}
\definecolor{parietal}{rgb}{1.00000,1.00000,0.00000}
\definecolor{insular}{rgb}{0.00000,1.00000,0.00000}
\definecolor{limbic}{rgb}{0.00000,1.00000,1.00000}
\definecolor{occipital}{rgb}{0.00000,0.00000,1.00000}
\definecolor{subcortical}{rgb}{1.00000,0.00000,1.00000}

\begin{tikzpicture}

\begin{axis}[%
width=4in,
height=4in,
at={(1.296in,1.089in)},
scale only axis,
xmin=0,
xmax=1,
xlabel={$c_\mathbf{PET1}$ [-]},
ymin=0,
ymax=1,
ylabel={$c_\mathbf{PET2}$ [-]},
axis background/.style={fill=white},
title style={font=\bfseries},
title={\textbf{Comparison of concentration in the graph nodes}},
xmajorgrids,
ymajorgrids,
]

\addplot[scatter, only marks, color=frontal, draw=black,
    visualization depends on={\thisrow{S} \as\perpointmarksize},
    scatter/@pre marker code/.style={/tikz/mark size=1.25*\perpointmarksize},
    scatter/@post marker code/.style={},
    forget plot] 
table[row sep=crcr]{%
x	y	S\\
0.472299466251897	0.730303560369322	2.39769758500825\\
0.442346196347881	0.675401262595871	2.55841316517235\\
0.418808122077406	0.669209527183732	2.80814481204823\\
0.342729055510406	0.639999469983393	2.43880989595229\\
0.460504609412906	0.607408651513455	2.54366750826483\\
0.366879614777177	0.591130355768607	2.66649862608477\\
0.548681683111425	0.640713500024626	2.48395621792868\\
0.463258981121392	0.697833008011364	2.36830146328667\\
0.475536039546213	0.744541599537225	2.39101011185538\\
0.563904658312597	0.715569263779	2.36830146328667\\
0.349319223428439	0.682010479096512	2.25168876990623\\
0.28009116152355	0.58350519392066	2.39324740219028\\
0.48067177671955	0.577376556571341	2.64691480891253\\
0.353305900794984	0.401226066509909	2.73334862093789\\
0.446518874139175	0.589124132821709	2.99007381857694\\
0.402480990109856	0.610516940413848	2.84606514421059\\
0.44303166655981	0.543070135047707	2.70596522085348\\
0.4009909430422	0.517198296858932	2.4260533533712\\
0.494492342527277	0.627310615956904	2.81503084112585\\
0.422096483261698	0.617415032036929	2.93958541173289\\
0.450791516132341	0.576909838961513	2.87740257143968\\
0.374422676749734	0.607288230218748	2.80120952718123\\
0.423185625857148	0.612777530672004	2.86577087235255\\
0.369171990079117	0.613780490287014	2.70750757417353\\
0.465368740009066	0.698983608904159	2.50368683632152\\
0.435790252993286	0.74194684843379	2.35200733767342\\
0.102892509888586	0.608070556965084	2.77723546954096\\
0.157290135495449	0.597118414527567	2.89522035813877\\
0.398261468943861	0.542698680146026	1.87048258949699\\
0.407679049082498	0.498339796905616	2.09376837697699\\
0.182854511802624	0.532261417408984	2.03463202989393\\
0.442104976809347	0.57301053702082	1.965819758943\\
0.116154085700735	0.314975125637265	1.94119043870726\\
0.0350610930927074	0.420375068238014	2.10040548129298\\
0.139314966492926	0.468499804985285	2.04900883581507\\
0.0626626792731517	0.542423800362877	2.03463202989393\\
0.181807188444743	0.506155861605643	2.19604996094467\\
0.168336649888798	0.521735422484574	2.30097781692202\\
0.396959272591118	0.426619820785785	1.95771280578869\\
0.164956167699813	0.405069141327073	1.80648108905035\\
0.249644573231119	0.601499468143638	2.14828016421654\\
0.220223403947863	0.563071394853713	2.38651946937749\\
0.210915720076316	0.458385557056823	2.19024734183694\\
0.148836085137238	0.539833683288267	2.1326783621607\\
0.178019053661435	0.544952844202078	2.80676257042579\\
0.14398148571011	0.545468246211741	2.87611784463149\\
0.017402543195496	0.535279557367501	2.38200622334003\\
0.0396993425210197	0.519663685481475	2.47188683505831\\
0.112706431485862	0.440443266070776	2.64856257814055\\
0.111224114641209	0.385757837382321	2.47995264875345\\
0.221682142101693	0.406717773287916	2.1513624208385\\
0.111138252445016	0.525220625373066	2.03098628797521\\
0.414511341630275	0.559003143045738	2.44511648598002\\
0.33595934966002	0.572948049716479	2.09042387919119\\
0.452918886785778	0.625856592158193	2.55841316517235\\
0.412428133786563	0.670660399792969	2.5853942127635\\
0.508096571412885	0.634981059046364	2.18439823408294\\
0.523637426181102	0.67836029773622	2.12315121172136\\
0.528475975857667	0.732395165229998	1.77449545213136\\
0.474801546003494	0.701024727546208	1.91992610946376\\
0.426954940223422	0.53704228994769	2.08030154837989\\
0.463048193353997	0.333810288015433	1.97382712896306\\
0.452833452100096	0.597477532326468	2.59244928474016\\
0.3835427927573	0.558945427612399	2.50951553477354\\
0.5420328381982	0.720187139436084	1.83189495538108\\
0.567997128608513	0.713291812742811	1.87514124329872\\
0.589357804781308	0.725316410053005	2.15747971205641\\
0.582210278804165	0.738150059547439	2.2490048269168\\};

\addplot[scatter, only marks, color=temporal, draw=black,
    visualization depends on={\thisrow{S} \as\perpointmarksize},
    scatter/@pre marker code/.style={/tikz/mark size=1.25*\perpointmarksize},
    scatter/@post marker code/.style={},
    forget plot] 
table[row sep=crcr]{%
x	y	S\\
0.324087454223393	0.431438332554036	2.60118679739897\\
0.385485839204031	0.332730791402454	2.54180468276311\\
0.503563299753163	0.550831933929344	1.75217015812413\\
0.472670018364091	0.567434747267965	1.89792964900184\\
0.511584663258195	0.565659994618203	2.18145785490508\\
0.421960506218634	0.399244490508401	2.31343588906969\\
0.163745093067783	0.463708373040506	2.05608744093053\\
0.358236856947784	0.427520260040965	2.14207951728752\\
0.511549268115833	0.434468798622336	2.1326783621607\\
0.475622407200074	0.307511074331889	2.31836387204215\\
0.558788811906595	0.550924774481356	2.20751549015728\\
0.52644917323065	0.577353061626095	1.55691210380244\\
0.0489366153567274	0.487326092539162	2.1326783621607\\
0.0966942404808378	0.458817733990148	2.24089133938669\\
0.116150887996778	0.540609450312455	2.63361196466191\\
0.297224877910132	0.505229839022247	2.78859569929087\\
0.0360175228726284	0.292086638642545	2.19893226788139\\
0.127716942284238	0.337090125415831	2.44511648598002\\
0.430559444398693	0.580393825432119	2.4174340310849\\
0.542479756443054	0.555199369965538	2.77294013047057\\
0.0859555068263444	0.499406784557319	1.84664991957027\\
0.00597244993475918	0.606813411244114	1.21835010308322\\
0.101922599616591	0.482725323663801	1.90681983263038\\
0.123055872802598	0.427423580786026	1.62525881390111\\
0.147847254776818	0.431167503515896	2.29340539060087\\
0.150622643749814	0.430652031941105	2.14207951728752\\
0.0268912125477579	0.556026018762185	2.05608744093053\\
0.0332403426268055	0.560377057440376	2.07689582351518\\
0.0756364828480951	0.346796311538301	2.17553737429575\\
0.116370217401745	0.443032199521266	1.87514124329872\\
0.180177031599249	0.527427077585102	2.30097781692202\\
0.106735059701772	0.450673274887843	1.94119043870726\\
0.0831575855817185	0.536217968630862	2.29340539060087\\
0.0846883761052405	0.611444958005234	2.27543645401692\\
0.245414602132896	0.448728465955701	2.62688424063143\\
0.131486749561419	0.458013284113906	2.64194871203641\\
0.483217919285437	0.670060090574938	2.65675980995807\\
0.579043033628186	0.602172209159109	2.54366750826483\\
0.423297474513382	0.657092450719401	2.83269703023825\\
0.538259533665816	0.623766747138702	2.8053782836504\\
0.121194250848842	0.309579152944781	1.48754005209535\\
0.0812326138244304	0.508561040881721	1.37262930039173\\
0.254450179634667	0.442278900439048	1.55691106873611\\
0.193777747336567	0.494769646515277	1.49675103352506\\
0.527884878905854	0.594527017641339	1.42831330225759\\
0.505313923867982	0.581801128827274	1.29534624398477\\
0.186281982608299	0.257096292581744	1.61808237512951\\
0.149395574003773	0.293147790074838	1.34832443735924\\
0.151385258755881	0.47273368263655	1.2662176467853\\
0.144179531851034	0.388722336259983	1.45882753298841\\
0.658757138058448	0.557329526368828	1.4685871206256\\
0.697458636744915	0.53594135596105	1.50579503469908\\
0.562241781508511	0.602849420560047	1.9112184095427\\
0.632792494965097	0.603875005479728	1.86578886453092\\
0.406399209356288	0.630874535863326	1.87514124329872\\
0.601421009792236	0.602685861018877	1.75217015812413\\
};

\addplot[scatter, only marks, color=parietal, draw=black,
    visualization depends on={\thisrow{S} \as\perpointmarksize},
    scatter/@pre marker code/.style={/tikz/mark size=1.25*\perpointmarksize},
    scatter/@post marker code/.style={},
    forget plot] 
table[row sep=crcr]{%
x	y	S\\
0.48999852118212	0.532885796343428	2.11672643485833\\
0.541221907047785	0.442503272965307	2.31343588906969\\
0.419059657162782	0.358192041145447	2.2381648337733\\
0.395058217133649	0.424820543683477	2.38200333311172\\
0.532670652532737	0.578547481935223	2.14828016421654\\
0.55878739040528	0.572324468079202	1.93699335221\\
0.461230532183266	0.564321968609798	2.1326783621607\\
0.540167038615	0.574119158720498	2.07346887953255\\
0.416258135716043	0.538036924017897	2.20180694452157\\
0.555062527251013	0.574795103980386	2.22993201272165\\
0.162857853716054	0.25426476712803	2.75701481179778\\
0.222558490406658	0.357336436866026	2.8566269396577\\
0.442246706708082	0.481346751185151	2.56387734234994\\
0.517126182029413	0.428607043935469	2.81365697450406\\
0.41653057823696	0.595734044325946	2.61669686616058\\
0.516251560347038	0.562206774512594	2.72884043134022\\
0.374199847078375	0.514808753183234	2.65185101550092\\
0.530569783633684	0.505859072808853	2.95985530128233\\
0.270435992039859	0.423114913644089	3.04573735712776\\
0.501019541665991	0.539520535967014	3.05113844456259\\
0.459191273700743	0.536957777190076	2.75847501665888\\
0.404674063966772	0.503285169520349	2.97275315359247\\
0.473902064525677	0.68133509693651	2.41090990273885\\
0.550876423190454	0.636945397353635	2.28831403994502\\
0.57238994270787	0.667705643126429	2.26497332558446\\
0.612634615473234	0.696456796324542	2.41308956869179\\
0.518425863182943	0.64472004664576	2.50951553477354\\
0.554012647690263	0.645928933417067	2.7044202256705\\
0.657071983525227	0.746389761132811	2.46169352771708\\
0.634100030054983	0.70339621307215	2.57828318642427\\
0.415875630221344	0.560080338924117	2.4174340310849\\
0.421915024790948	0.555291195957035	2.44929303697071\\
0.505794989508976	0.444562472265897	2.31096162423687\\
0.520148049729525	0.402122046123185	2.39547559383704\\
0.492645050528698	0.523740963083115	2.30097781692202\\
0.49335588003722	0.594523178965953	2.26233644171172\\
0.502196979460083	0.605270295978546	2.24630779411012\\
0.559474204180837	0.709300658337632	2.05608744093053\\
};
\addplot[scatter, only marks, color=insular, draw=black,
    visualization depends on={\thisrow{S} \as\perpointmarksize},
    scatter/@pre marker code/.style={/tikz/mark size=1.25*\perpointmarksize},
    scatter/@post marker code/.style={},
    forget plot] 
table[row sep=crcr]{%
x	y	S\\
0.225443035880174	0.473398936761164	1.90681983263038\\
0.1959263633542	0.504110855372328	1.84664991957027\\
0.289710925599721	0.544406150808847	1.74068018665178\\
0.44434913584029	0.526787674555156	1.57288210430251\\
0.381940311620883	0.611254689315737	1.75217015812413\\
0.240545548081647	0.515085348258803	1.68593109974047\\
0.360581801043419	0.412363272172655	1.75783210386821\\
0.25269340481471	0.448566915397748	1.72300647584076\\
0.377679479894005	0.553264347529386	1.8889133308129\\
0.405418402620592	0.511027234501092	1.80126765565868\\
0.460637649766466	0.577847617758266	2.00491835104254\\
0.485702871029934	0.545137616808099	1.75217015812413\\
};

\addplot[scatter, only marks, color=limbic, draw=black,
    visualization depends on={\thisrow{S} \as\perpointmarksize},
    scatter/@pre marker code/.style={/tikz/mark size=1.25*\perpointmarksize},
    scatter/@post marker code/.style={},
    forget plot] 
table[row sep=crcr]{%
x	y	S\\
0.587456072616102	0.768074444111761	1.77449545213136\\
0.539695182660642	0.752155931684154	1.57288210430251\\
0.00749676985634193	0.496684759253245	1.02450514442835\\
0.198082660979359	0.334449343424405	1.55691210380244\\
0.160036481886188	0.308832960519476	1.97382712896306\\
0.104867967811761	0.490029268573319	1.71699178894595\\
0.497102331019016	0.666032018705414	1.90239060393323\\
0.476995250373515	0.621110612686013	1.89343759033074\\
0.0654423858409746	0.642272424804242	1.65305035594409\\
0.144627713467414	0.648739100257685	1.40682796115230\\
0.512821754287123	0.717705729074967	2.05255955140931\\
0.524476927172799	0.713638369493466	1.71091400033305\\
0.311999165586399	0.471081608630221	2.36830146328667\\
0.453683421256996	0.426966231957797	2.04544412917595\\
};

\addplot[scatter, only marks, color=occipital, draw=black,
    visualization depends on={\thisrow{S} \as\perpointmarksize},
    scatter/@pre marker code/.style={/tikz/mark size=1.25*\perpointmarksize},
    scatter/@post marker code/.style={},
    forget plot] 
table[row sep=crcr]{%
x	y	S\\
0.357160019522206	0.721420615950644	2.23543182424522\\
0.173601239700207	0.585793115094443	2.36830146328667\\
0.647872988328341	0.681332020723659	2.47390940268454\\
0.694816112228754	0.699893009249123	2.47995264875345\\
0.103266253397258	0.29297767734449	2.2490048269168\\
0.260525504415174	0.408835898574318	2.24360447034614\\
0.637010242943846	0.613084534719928	2.58892895839153\\
0.661963087482388	0.605049314342313	2.58716452446685\\
0.614888866335126	0.497509544145087	2.51530390001822\\
0.634198581621659	0.625172537633757	2.45137199943362\\
0.0963627791960577	0.267104010816484	2.43458177402548\\
0.0902441219668858	0.185133802730468	2.48594949313945\\
0.136567160950028	0.384559233954541	2.45137199943362\\
0.151170481810518	0.286522820897	2.60118679739897\\
0.0634579559007619	0.1943249583551	2.37289613453954\\
0.0662222044936975	0.185287419677695	2.56206060688068\\
0.0797160445291545	0.349529975565696	2.81502909002577\\
0.128486774579046	0.35844278191794	2.70750757417353\\
0.223244877393349	0.216036278132348	2.34492012712561\\
0.168424192925227	0.210820988590181	2.37746426960028\\
0.378129800224137	0.384799380378955	2.45137199943362\\
};

\addplot[scatter, only marks, color=subcortical, draw=black,
    visualization depends on={\thisrow{S} \as\perpointmarksize},
    scatter/@pre marker code/.style={/tikz/mark size=1.25*\perpointmarksize},
    scatter/@post marker code/.style={},
    forget plot] 
table[row sep=crcr]{%
x	y	S\\
0.411195961071332	0.426913490085256	2.64856257814055\\
0.164747728154196	0.465925995021822	1.58837996122671\\
0.145704589411526	0.448445083288704	1.83189495538108\\
0.18715234342867	0.483072794782989	1.02450550768818\\
0.418964292453354	0.387410231649518	1.29534444675614\\
0.071682190964823	0.447684828972702	2.57649449174788\\
0.0943330747922438	0.445775069252078	2.43669858620224\\
0.496995616555302	0.43035330235918	2.48395621792868\\
0.412661329719337	0.35233268197362	2.50951553477354\\
0.371093091993409	0.584949386623421	2.33776806711878\\
0.230549021849982	0.575008343494003	1.85149033582191\\
0.499221918586711	0.660487893674786	1.93699388970628\\
0.470860384121188	0.580422837627646	2.07346887953255\\
0.286830548293591	0.441371978116083	1.76343914890078\\
0.189833570622024	0.436505474004615	1.99728451366719\\
0.41956238629056	0.665048747620787	1.76343914890078\\
0.455606771281823	0.66615068229478	1.66643698377555\\
0.595379039938651	0.313609321080038	2.81913692983967\\
0.420760051557108	0.0745630119329758	2.94439309105273\\
0.430932261597805	0.610582902074056	1.89343759033074\\
0.459791570130881	0.622767291198559	2.05608744093053\\
0.640098385283447	0.587348850453387	1.58837996122671\\
0.561958347329526	0.471696338595902	1.4685871206256\\
0.465067086802554	0.708191555706299	1.43870144968745\\
0.488010659958531	0.753395059806601	1.49675103352506\\
0.439986578992476	0.514433053607443	1.43870144968745\\
0.330275853074453	0.545511178695407	1.83189495538108\\
0.508433701984883	0.423126578469898	1.36064044418607\\
0.399070512262624	0.249762986429015	1.52341145624943\\
0.345442631879026	0.3732948695038	1.53199306227594\\
0.530317103376352	0.356124799026637	1.54873927365813\\
0.501103552953136	0.347044775703733	1.61081009433214\\
0.326398994153838	0.387124722371142	0.994777988102889\\
0.434097893858314	0.207175609995582	1.12315562019777\\
0.372270628090608	0.215558153843012	1.0518503294417\\
0.460749593078043	0.677835175190611	2.06657241422903\\
0.645701489194939	0.428294703840555	1.62525881390111\\
};

\addplot [color=black, line width=2.0pt, forget plot]
  table[row sep=crcr]{%
0	0\\
1	1\\
};

\end{axis}

\end{tikzpicture}%

%% file: Solution_Global.tikz
\begin{tikzpicture}

\begin{axis}[%
width=3in,
height=2.25in,
at={(0.758in,0.481in)},
scale only axis,
xmin=0,
xmax=20,
xlabel style={font=\color{white!15!black}},
xlabel={$t$ [years]},
ymin=0.3,
ymax=1,
ylabel style={font=\color{white!15!black}},
ylabel={$\langle c_h\rangle$(t) [-]},
axis background/.style={fill=white},
title style={font=\bfseries},
title={\textbf{Global solution mean and std range}},
axis x line*=bottom,
axis y line*=left
]

\addplot[area legend, draw=violet, fill=violet, draw opacity=0.3, fill opacity=0.3, forget plot]
table[row sep=crcr] {%
x	y\\
0	0.522121380481211\\
1	0.559497055108961\\
2	0.596286427436778\\
3	0.63113883564001\\
4	0.663704427798274\\
5	0.693824097329287\\
6	0.721477941194614\\
7	0.746738024003376\\
8	0.769731792636113\\
9	0.790616249856428\\
10	0.809560694183721\\
11	0.826735637335418\\
12	0.842305989098635\\
13	0.856427127666995\\
14	0.86924288786386\\
15	0.880884797021817\\
16	0.891472097656768\\
17	0.901112243516917\\
18	0.909901659285443\\
19	0.917926626564114\\
20	0.925264208590669\\
20	0.775843432702249\\
19	0.769599757857031\\
18	0.762907463672216\\
17	0.755727396631719\\
16	0.748017109459632\\
15	0.739730868693062\\
14	0.730819830860309\\
13	0.721232474860798\\
12	0.710915411288605\\
11	0.699814729931346\\
10	0.687878091217854\\
9	0.675057806420368\\
8	0.661315163941866\\
7	0.646626205233716\\
6	0.630988967386227\\
5	0.614431797413392\\
4	0.597021622357756\\
3	0.578870091984783\\
2	0.560134821845068\\
1	0.541013903271928\\
0	0.522121223487736\\
}--cycle;
\addplot [color=violet, line width=2.0pt, forget plot]
  table[row sep=crcr]{%
0	0.522121301984474\\
1	0.550255479190444\\
2	0.578210624640923\\
3	0.605004463812397\\
4	0.630363025078015\\
5	0.65412794737134\\
6	0.67623345429042\\
7	0.696682114618546\\
8	0.71552347828899\\
9	0.732837028138398\\
10	0.748719392700787\\
11	0.763275183633382\\
12	0.77661070019362\\
13	0.788829801263897\\
14	0.800031359362085\\
15	0.81030783285744\\
16	0.8197446035582\\
17	0.828419820074318\\
18	0.836404561478829\\
19	0.843763192210572\\
20	0.850553820646459\\
};
\end{axis}

\end{tikzpicture}%

%% file: Coll_Err_Global.tikz
\begin{tikzpicture}

\begin{axis}[%
width=3in,
height=2.25in,
at={(0.758in,0.481in)},
scale only axis,
xmode=log,
xmin=500,
xmax=200000,
xminorticks=true,
xlabel style={font=\color{white!15!black}},
xlabel={Collocation Points},
ymode=log,
ymin=1e-12,
ymax=0.01,
yminorticks=true,
ylabel style={font=\color{white!15!black}},
ylabel={Collocation Error},
axis background/.style={fill=white},
title style={font=\bfseries},
title={\textbf{SC Error on $\boldsymbol{\mu}_Q$}},
xmajorgrids,
xminorgrids,
ymajorgrids,
yminorgrids,
legend style={at={(0.03,0.03)}, anchor=south west, legend cell align=left, align=left, draw=white!15!black}
]
\addplot [color=violet, line width=1.5pt, draw opacity=0.25]
  table[row sep=crcr]{%
575	3.81082369244492e-07\\
2241	2.15984668611688e-08\\
7183	4.5418133484768e-09\\
19825	1.34211709291864e-09\\
48639	1.37065785303823e-11\\
};

\addplot [color=violet, line width=1.5pt,draw opacity=0.5]
  table[row sep=crcr]{%
575	1.44332165880409e-05\\
2241	1.77526223423021e-06\\
7183	8.69570931215733e-07\\
19825	6.23388277635347e-07\\
48639	1.79335051493992e-08\\
};

\addplot [color=violet, line width=1.5pt, draw opacity=0.75]
  table[row sep=crcr]{%
575	0.000117253678902503\\
2241	2.98342572140319e-05\\
7183	2.17746310374133e-05\\
19825	1.70677073338447e-05\\
48639	4.48140528529935e-07\\
};

\addplot [color=violet, line width=1.5pt]
  table[row sep=crcr]{%
575	0.0004299391725742\\
2241	0.000135612877179825\\
7183	0.000101088412185069\\
19825	8.21253916054837e-05\\
48639	4.2204694895398e-06\\
};

\end{axis}

\end{tikzpicture}%

%% file: MC_Err_Global.tikz
\begin{tikzpicture}

\begin{axis}[%
width=3in,
height=2.25in,
at={(0.758in,0.481in)},
scale only axis,
xmode=log,
xmin=500,
xmax=200000,
xminorticks=true,
xlabel style={font=\color{white!15!black}},
xlabel={Iterations},
ymode=log,
ymin=1e-12,
ymax=0.01,
yminorticks=true,
ylabel style={font=\color{white!15!black}},
ylabel={Monte Carlo Error},
axis background/.style={fill=white},
title style={font=\bfseries},
title={\textbf{MC Error on $\boldsymbol{\mu}_Q$}},
xmajorgrids,
xminorgrids,
ymajorgrids,
yminorgrids,
legend style={at={(0.03,0.03)}, anchor=south west, legend cell align=left, align=left, draw=white!15!black}
]

\addplot [color=black, dashed, line width=1.5pt]
  table[row sep=crcr]{%
5000	0.014142135623731\\
10000	0.01\\
15000	0.00816496580927726\\
20000	0.00707106781186548\\
25000	0.00632455532033676\\
30000	0.00577350269189626\\
35000	0.00534522483824849\\
40000	0.005\\
45000	0.00471404520791032\\
50000	0.00447213595499958\\
55000	0.00426401432711221\\
60000	0.00408248290463863\\
65000	0.00392232270276368\\
70000	0.00377964473009227\\
75000	0.00365148371670111\\
80000	0.00353553390593274\\
85000	0.00342997170285018\\
90000	0.00333333333333333\\
95000	0.00324442842261525\\
100000	0.00316227766016838\\
105000	0.00308606699924184\\
110000	0.00301511344577764\\
115000	0.00294883912309794\\
120000	0.00288675134594813\\
125000	0.00282842712474619\\
130000	0.00277350098112615\\
135000	0.00272165526975909\\
140000	0.00267261241912424\\
145000	0.00262612865719445\\
150000	0.00258198889747161\\
155000	0.00254000254000381\\
160000	0.0025\\
165000	0.00246182981958665\\
170000	0.00242535625036333\\
175000	0.00239045721866879\\
180000	0.00235702260395516\\
185000	0.00232495277487639\\
190000	0.00229415733870562\\
195000	0.00226455406828919\\
200000	0.00223606797749979\\
};
\addlegendentry{$M^{-1/2}$}

\addplot [color=violet, line width=1.5pt, draw opacity=0.25]
  table[row sep=crcr]{%
5000	0.00170416464301774\\
10000	0.00131115589835482\\
15000	0.00105919137541153\\
20000	0.000884898392427518\\
25000	0.0008473082801836\\
30000	0.000887792499300074\\
35000	0.000836386245149029\\
40000	0.000674622310565631\\
45000	0.000455522058206584\\
50000	0.000481563247663104\\
55000	0.000481410222367919\\
60000	0.000450749812167512\\
65000	0.000384704847353823\\
70000	0.000409000473201025\\
75000	0.000399018468548663\\
80000	0.000438800425006919\\
85000	0.000421411674033662\\
90000	0.000364125409696267\\
95000	0.000357881458285189\\
100000	0.000412204198555686\\
105000	0.000413646547626804\\
110000	0.0003915562795715\\
115000	0.000411587011017759\\
120000	0.000365322435167913\\
125000	0.000339744167372243\\
130000	0.000342253037178801\\
135000	0.000269405570584303\\
140000	0.000331434522149603\\
145000	0.000328455674006143\\
150000	0.000305497448127175\\
155000	0.000317038323423402\\
160000	0.000300218487392022\\
165000	0.000278331819373053\\
170000	0.00028250054157534\\
175000	0.000274480170664519\\
180000	0.000297955594877473\\
185000	0.000326822288134998\\
190000	0.000330589405018386\\
195000	0.000329858533179092\\
200000	0.000303093763984292\\
};
\addlegendentry{$t=5$}

\addplot [color=violet, line width=1.5pt,draw opacity=0.5]
  table[row sep=crcr]{%
5000	0.0021957396036024\\
10000	0.00161996516427688\\
15000	0.00130433067110375\\
20000	0.00110344248209136\\
25000	0.0010778651391051\\
30000	0.00109543290678017\\
35000	0.00105025351178698\\
40000	0.00087062663513451\\
45000	0.000578426395368082\\
50000	0.000598271144009984\\
55000	0.000599934528570134\\
60000	0.000575956652617321\\
65000	0.000482658017877235\\
70000	0.000522754816550031\\
75000	0.000494834329440983\\
80000	0.00054332407935083\\
85000	0.000528074853084939\\
90000	0.000463070769013728\\
95000	0.00046885188672408\\
100000	0.000537429892775756\\
105000	0.000537192442376069\\
110000	0.000502857146836657\\
115000	0.000537252476256019\\
120000	0.000460384942699337\\
125000	0.000437540999007482\\
130000	0.000438754447968597\\
135000	0.000338848582254403\\
140000	0.000418939752834495\\
145000	0.000413498604843605\\
150000	0.000393159974770522\\
155000	0.000399760717616802\\
160000	0.000381167111346598\\
165000	0.000349799359653285\\
170000	0.00035723624342728\\
175000	0.000344267640753464\\
180000	0.000384886320427978\\
185000	0.000426763011635257\\
190000	0.000432107839446411\\
195000	0.000425084401214046\\
200000	0.00038328891275141\\
};
\addlegendentry{$t=10$}

\addplot [color=violet, line width=1.5pt, draw opacity=0.75]
  table[row sep=crcr]{%
5000	0.00220504082943621\\
10000	0.00162186369316794\\
15000	0.0013382611467659\\
20000	0.00116181424369021\\
25000	0.00109369684732445\\
30000	0.00105527055526562\\
35000	0.00103356539743561\\
40000	0.000863529217588007\\
45000	0.000561272603282781\\
50000	0.000567693189412043\\
55000	0.000586078908687107\\
60000	0.000571398543399451\\
65000	0.000475898447758375\\
70000	0.000513942415138774\\
75000	0.000476882222781345\\
80000	0.000517078942133538\\
85000	0.000506445156518731\\
90000	0.00044801940009283\\
95000	0.000472372919626252\\
100000	0.000541267812783563\\
105000	0.000540177870796418\\
110000	0.000500633377581244\\
115000	0.000545054399778817\\
120000	0.000454121088043683\\
125000	0.000440182357043628\\
130000	0.00043903477696138\\
135000	0.000332769233333677\\
140000	0.000411072615450687\\
145000	0.000406068493090115\\
150000	0.00039831105292517\\
155000	0.000398247625056018\\
160000	0.000379590640613178\\
165000	0.000352534014944948\\
170000	0.000357809823506421\\
175000	0.000344837725045851\\
180000	0.000391172798926164\\
185000	0.00043974678025794\\
190000	0.000445860607837232\\
195000	0.00043220578996828\\
200000	0.000381002032742933\\
};
\addlegendentry{$t=15$}

\addplot [color=violet, line width=1.5pt]
  table[row sep=crcr]{%
5000	0.00207353048378346\\
10000	0.00154475658345637\\
15000	0.00132320153590919\\
20000	0.00117814225584178\\
25000	0.00103841825936941\\
30000	0.000950271801206094\\
35000	0.000951574292638436\\
40000	0.000796345775909804\\
45000	0.000510168150508595\\
50000	0.000505779620913128\\
55000	0.000538008528212352\\
60000	0.000532863497623295\\
65000	0.000446598594309261\\
70000	0.000473997106916521\\
75000	0.000429235582109876\\
80000	0.000459581633729309\\
85000	0.000453838106752495\\
90000	0.000405274424648187\\
95000	0.000447813333742504\\
100000	0.000512790169548251\\
105000	0.000512583954059318\\
110000	0.000471019635882308\\
115000	0.000521520925640493\\
120000	0.000425632623551868\\
125000	0.000420377627449213\\
130000	0.000415956465861796\\
135000	0.000309755623056538\\
140000	0.00037936540378327\\
145000	0.00037533062444026\\
150000	0.000381453835082779\\
155000	0.000375844481444398\\
160000	0.000357291638259201\\
165000	0.000340340625146088\\
170000	0.000342062000996306\\
175000	0.000332370391065845\\
180000	0.000376368025179246\\
185000	0.000427757260479098\\
190000	0.000434584944710438\\
195000	0.000415956830440194\\
200000	0.000359006097710423\\
};
\addlegendentry{$t=20$}

\end{axis}

\end{tikzpicture}%

%% file: Solution_1_Frontal.tikz
\begin{tikzpicture}

\begin{axis}[%
width=3in,
height=1.25in,
at={(0.758in,0.481in)},
scale only axis,
xmin=0,
xmax=20,
xlabel style={font=\color{white!15!black}},
xlabel={$t$ [years]},
ymin=0.3,
ymax=1,
ylabel style={font=\color{white!15!black}},
ylabel={$\langle c_h\rangle_1$(t) [-]},
axis background/.style={fill=white},
title style={font=\bfseries},
title={\textbf{Frontal Lobe: solution mean and std range}},
axis x line*=bottom,
axis y line*=left
]

\addplot[area legend, draw=red, fill=red, draw opacity=0.3, fill opacity=0.3, forget plot]
table[row sep=crcr] {%
x	y\\
0	0.586513985491474\\
1	0.644049922680315\\
2	0.697581344686663\\
3	0.745240321263468\\
4	0.786909671957426\\
5	0.822864956938747\\
6	0.853615834659563\\
7	0.879775245688323\\
8	0.901970406619053\\
9	0.920791121593933\\
10	0.936764961056394\\
11	0.950349337211904\\
12	0.96193305934584\\
13	0.971842530670057\\
14	0.980349706902487\\
15	0.987680239861237\\
16	0.994021027862794\\
17	0.999526853840008\\
18	1.00432603995908\\
19	1.00852517077413\\
20	1.01221299002068\\
20	0.845519905752314\\
19	0.840010943224544\\
18	0.833995978927199\\
17	0.827421965861605\\
16	0.820230947020207\\
15	0.812360189728485\\
14	0.803742674523338\\
13	0.794308114377583\\
12	0.783984741182226\\
11	0.772702165713807\\
10	0.760395681407154\\
9	0.747012408150073\\
8	0.732519591074672\\
7	0.716915053422264\\
6	0.700239047574915\\
5	0.682585297644993\\
4	0.664106743005473\\
3	0.645008858235342\\
2	0.625522480154689\\
1	0.60585304659583\\
0	0.586511969204125\\
}--cycle;
\addplot [color=red, line width=2.0pt, forget plot]
  table[row sep=crcr]{%
0	0.586512977347799\\
1	0.624951484638073\\
2	0.661551912420676\\
3	0.695124589749405\\
4	0.725508207481449\\
5	0.75272512729187\\
6	0.776927441117239\\
7	0.798345149555294\\
8	0.817244998846863\\
9	0.833901764872003\\
10	0.848580321231774\\
11	0.861525751462856\\
12	0.872958900264033\\
13	0.88307532252382\\
14	0.892046190712913\\
15	0.900020214794861\\
16	0.907125987441501\\
17	0.913474409850806\\
18	0.919161009443138\\
19	0.924268056999336\\
20	0.928866447886499\\
};
\end{axis}

\end{tikzpicture}%

%% file: Coll_Err_1_Frontal.tikz
\begin{tikzpicture}

\begin{axis}[%
width=3in,
height=1.25in,
at={(0.758in,0.481in)},
scale only axis,
xmode=log,
xmin=500,
xmax=200000,
xminorticks=true,
xlabel style={font=\color{white!15!black}},
xlabel={Collocation Points},
ymode=log,
ymin=1e-12,
ymax=0.01,
yminorticks=true,
ylabel style={font=\color{white!15!black}},
ylabel={Collocation Error},
axis background/.style={fill=white},
title style={font=\bfseries},
title={\textbf{Frontal Lobe: SC Error on $\mu_Q$}},
xmajorgrids,
xminorgrids,
ymajorgrids,
yminorgrids,
legend style={at={(0.03,0.03)}, anchor=south west, legend cell align=left, align=left, draw=white!15!black}
]
\addplot [color=red, line width=1.5pt, draw opacity=0.25]
  table[row sep=crcr]{
575	5.79957325732549e-08\\
2241	9.65999468278901e-09\\
7183	3.88911682046578e-09\\
19825	1.4492327099818e-09\\
48639	1.46752193855782e-11\\
};

\addplot [color=red, line width=1.5pt,draw opacity=0.5]
  table[row sep=crcr]{%
575	2.09030007278479e-05\\
2241	3.13002145962407e-06\\
7183	9.59258929501098e-07\\
19825	1.7844132515515e-07\\
48639	1.15107680483264e-08\\
};

\addplot [color=red, line width=1.5pt, draw opacity=0.75]
  table[row sep=crcr]{%
575	0.000111253822922689\\
2241	3.62316521790194e-06\\
7183	9.05170164583194e-06\\
19825	1.2926625714221e-05\\
48639	6.02368489037968e-07\\
};

\addplot [color=red, line width=1.5pt]
  table[row sep=crcr]{%
575	0.000112139723457469\\
2241	6.72950288927772e-05\\
7183	8.8470836924166e-05\\
19825	8.99297624474176e-05\\
48639	3.81934181287391e-06\\
};

\end{axis}

\end{tikzpicture}%

%% file: MC_Err_1_Frontal.tikz
\begin{tikzpicture}

\begin{axis}[%
width=3in,
height=1.25in,
at={(0.758in,0.481in)},
scale only axis,
xmode=log,
xmin=500,
xmax=200000,
xminorticks=true,
xlabel style={font=\color{white!15!black}},
xlabel={Iterations},
ymode=log,
ymin=1e-12,
ymax=0.01,
yminorticks=true,
ylabel style={font=\color{white!15!black}},
ylabel={Monte Carlo Error},
axis background/.style={fill=white},
title style={font=\bfseries},
title={\textbf{Frontal Lobe: MC Error on $\mu_Q$}},
xmajorgrids,
xminorgrids,
ymajorgrids,
yminorgrids,
legend style={at={(0.03,0.03)}, anchor=south west, legend cell align=left, align=left, draw=white!15!black}
]

\addplot [color=black, dashed, line width=1.5pt]
  table[row sep=crcr]{%
5000	0.014142135623731\\
10000	0.01\\
15000	0.00816496580927726\\
20000	0.00707106781186548\\
25000	0.00632455532033676\\
30000	0.00577350269189626\\
35000	0.00534522483824849\\
40000	0.005\\
45000	0.00471404520791032\\
50000	0.00447213595499958\\
55000	0.00426401432711221\\
60000	0.00408248290463863\\
65000	0.00392232270276368\\
70000	0.00377964473009227\\
75000	0.00365148371670111\\
80000	0.00353553390593274\\
85000	0.00342997170285018\\
90000	0.00333333333333333\\
95000	0.00324442842261525\\
100000	0.00316227766016838\\
105000	0.00308606699924184\\
110000	0.00301511344577764\\
115000	0.00294883912309794\\
120000	0.00288675134594813\\
125000	0.00282842712474619\\
130000	0.00277350098112615\\
135000	0.00272165526975909\\
140000	0.00267261241912424\\
145000	0.00262612865719445\\
150000	0.00258198889747161\\
155000	0.00254000254000381\\
160000	0.0025\\
165000	0.00246182981958665\\
170000	0.00242535625036333\\
175000	0.00239045721866879\\
180000	0.00235702260395516\\
185000	0.00232495277487639\\
190000	0.00229415733870562\\
195000	0.00226455406828919\\
200000	0.00223606797749979\\
};
\addlegendentry{$M^{-1/2}$}

\addplot [color=red, line width=1.5pt, draw opacity=0.25]
  table[row sep=crcr]{%
5000	0.000863111312337756\\
10000	0.000399367489417947\\
15000	0.000852465401713917\\
20000	0.000697077507984219\\
25000	0.000463606070568177\\
30000	0.000116041072327084\\
35000	3.4074221085337e-05\\
40000	0.000239684300138945\\
45000	0.000102923543431576\\
50000	6.56129488240786e-05\\
55000	0.000141818210159477\\
60000	7.1296273940682e-05\\
65000	0.000160814555929974\\
70000	9.2042865282671e-05\\
75000	1.77491000513058e-05\\
80000	8.13393561884675e-05\\
85000	0.000241824607786872\\
90000	0.00018190090260481\\
95000	0.000201614310594779\\
100000	0.000154771109112536\\
105000	4.6786669084506e-05\\
110000	2.59327383677972e-05\\
115000	8.84809614787631e-05\\
120000	9.221222370301e-05\\
125000	6.35680916230773e-05\\
130000	6.4345397237092e-05\\
135000	1.71044491428574e-06\\
140000	3.14672160116849e-05\\
145000	3.15024410945961e-05\\
150000	5.84154962793407e-06\\
155000	4.1718351054623e-05\\
160000	2.53654256171217e-05\\
165000	3.04301095882841e-05\\
170000	1.30687764899975e-05\\
175000	1.17491181336656e-05\\
180000	2.8124282236588e-06\\
185000	1.41971258165845e-05\\
190000	1.27775501015632e-06\\
195000	1.83535359254533e-05\\
200000	3.94773358935385e-06\\
};
\addlegendentry{$t=5$}

\addplot [color=red, line width=1.5pt,draw opacity=0.5]
  table[row sep=crcr]{%
5000	0.00127680466519761\\
10000	0.000117616450502707\\
15000	0.000682507197933219\\
20000	0.0006111863850011\\
25000	0.000364039270992204\\
30000	0.000258393233768713\\
35000	6.40763628857187e-05\\
40000	0.000375431993286698\\
45000	0.000161598797385223\\
50000	0.000131946680235888\\
55000	9.88457499267363e-05\\
60000	0.000152933112398821\\
65000	0.000218880742822937\\
70000	0.000146783703810714\\
75000	2.15223489332871e-05\\
80000	9.36905719231736e-05\\
85000	0.000290542235070621\\
90000	0.00025129198359322\\
95000	0.000282716947045246\\
100000	0.000261091454386658\\
105000	0.000152936115693182\\
110000	6.68083279522569e-05\\
115000	3.77624897627487e-06\\
120000	2.02771632543261e-05\\
125000	7.93059250032231e-06\\
130000	1.2603677077462e-06\\
135000	5.12891333101825e-05\\
140000	0.000106499664242496\\
145000	0.000121157050344771\\
150000	8.09074085510745e-05\\
155000	4.97437456193376e-05\\
160000	6.11446576629748e-05\\
165000	5.27220113153513e-05\\
170000	5.8620099733838e-05\\
175000	8.12307871415949e-05\\
180000	8.02495904359292e-05\\
185000	0.000103887739388729\\
190000	7.27656273441529e-05\\
195000	8.97067488696085e-05\\
200000	6.65202118797975e-05\\
};
\addlegendentry{$t=10$}

\addplot [color=red, line width=1.5pt, draw opacity=0.75]
  table[row sep=crcr]{%
5000	0.00150594388785803\\
10000	0.000642869660883806\\
15000	0.000359018047167364\\
20000	0.000389768711751242\\
25000	0.000182126307262039\\
30000	0.000363934436246216\\
35000	0.000192229100500901\\
40000	0.000491498939115497\\
45000	0.000256313390217182\\
50000	0.000236937972732874\\
55000	5.80703796347308e-06\\
60000	0.000241900629019781\\
65000	0.000267319981814492\\
70000	0.000201035317915011\\
75000	5.51890352853681e-05\\
80000	0.000124693167201658\\
85000	0.000311333651200141\\
90000	0.000294195128970914\\
95000	0.00032392252323048\\
100000	0.000333210062078241\\
105000	0.000249506182362927\\
110000	0.000160426916260406\\
115000	0.000109480763264073\\
120000	7.28259551460828e-05\\
125000	9.62103424878367e-05\\
130000	8.37219684577223e-05\\
135000	0.000108467096886452\\
140000	0.000169927099654948\\
145000	0.000193181850876342\\
150000	0.000162259462203725\\
155000	0.000138240806129968\\
160000	0.000140873459449504\\
165000	0.000130599889561744\\
170000	0.000122425477830178\\
175000	0.000136446278540999\\
180000	0.000146625228742503\\
185000	0.00017931387633312\\
190000	0.000141648293199047\\
195000	0.00015325977761983\\
200000	0.000123732278771514\\
};
\addlegendentry{$t=15$}

\addplot [color=red, line width=1.5pt]
  table[row sep=crcr]{%
5000	0.00165224934030243\\
10000	0.000997585107216188\\
15000	0.000100774262061082\\
20000	0.000199792814176834\\
25000	3.25591922928093e-05\\
30000	0.000430256189469774\\
35000	0.000303579798555147\\
40000	0.000584895309320771\\
45000	0.000353702044446639\\
50000	0.000337672549255617\\
55000	0.000108759013302188\\
60000	0.000321216541228865\\
65000	0.000311727397644202\\
70000	0.000249058329670289\\
75000	0.000100042817178863\\
80000	0.000165639159696914\\
85000	0.000330926012638172\\
90000	0.000327794188304993\\
95000	0.000349270059122594\\
100000	0.00037894819461621\\
105000	0.000317660087534398\\
110000	0.000228775296428764\\
115000	0.00018819480542772\\
120000	0.000145126649807747\\
125000	0.000164583815492264\\
130000	0.000148198406578579\\
135000	0.000151583192793742\\
140000	0.000210056369485991\\
145000	0.000235991168309567\\
150000	0.000217480987724914\\
155000	0.000197555977979228\\
160000	0.000192489092069032\\
165000	0.000181838322602507\\
170000	0.00016322742664286\\
175000	0.000169255714740555\\
180000	0.00018829527055378\\
185000	0.000226399326626263\\
190000	0.000188355713489767\\
195000	0.000195010660348935\\
200000	0.000161191641686173\\
};
\addlegendentry{$t=20$}

\end{axis}

\end{tikzpicture}%

%% file: Solution_2_Temporal.tikz
\definecolor{temporal}{rgb}{1.00000,0.41800,0.01950}%

\begin{tikzpicture}

\begin{axis}[%
width=3in,
height=1.25in,
at={(0.758in,0.481in)},
scale only axis,
xmin=0,
xmax=20,
xlabel style={font=\color{white!15!black}},
xlabel={$t$ [years]},
ymin=0.3,
ymax=1,
ylabel style={font=\color{white!15!black}},
ylabel={$\langle c_h\rangle_2$(t) [-]},
axis background/.style={fill=white},
title style={font=\bfseries},
title={\textbf{Temporal Lobe: solution mean and std range}},
axis x line*=bottom,
axis y line*=left
]

\addplot[area legend, draw=temporal, fill=temporal, draw opacity=0.3, fill opacity=0.3, forget plot]
table[row sep=crcr] {%
x	y\\
0	0.504012628524858\\
1	0.554098826961502\\
2	0.603314925610334\\
3	0.649541287499901\\
4	0.692097310351395\\
5	0.730667431296277\\
6	0.765221047158035\\
7	0.795922616287855\\
8	0.823052334767692\\
9	0.846945951340113\\
10	0.867953910179292\\
11	0.886416211145661\\
12	0.90264858597321\\
13	0.916936168102318\\
14	0.92953179308562\\
15	0.940656961747542\\
16	0.950504191695819\\
17	0.959239971108497\\
18	0.967007853013702\\
19	0.973931434894957\\
20	0.980117095650774\\
20	0.75965017229938\\
19	0.751873000062148\\
18	0.743526131503573\\
17	0.734567271237294\\
16	0.724953853751461\\
15	0.71464446339508\\
14	0.703600884012179\\
13	0.691790978195086\\
12	0.679192633029924\\
11	0.665799031088409\\
10	0.651625487733149\\
9	0.63671799584907\\
8	0.62116337069574\\
7	0.605100404621821\\
6	0.588730651748892\\
5	0.572326414659152\\
4	0.556232601996585\\
3	0.540859465857374\\
2	0.526666682721763\\
1	0.514147326153364\\
0	0.504011693781299\\
}--cycle;
\addplot [color=temporal, line width=2.0pt, forget plot]
  table[row sep=crcr]{%
0	0.504012161153078\\
1	0.534123076557433\\
2	0.564990804166049\\
3	0.595200376678638\\
4	0.62416495617399\\
5	0.651496922977715\\
6	0.676975849453464\\
7	0.700511510454838\\
8	0.722107852731716\\
9	0.741831973594591\\
10	0.759789698956221\\
11	0.776107621117035\\
12	0.790920609501567\\
13	0.804363573148702\\
14	0.8165663385489\\
15	0.827650712571311\\
16	0.83772902272364\\
17	0.846903621172895\\
18	0.855266992258638\\
19	0.862902217478553\\
20	0.869883633975077\\
};
\end{axis}

\end{tikzpicture}%

%% file: Coll_Err_2_Temporal.tikz
\definecolor{temporal}{rgb}{1.00000,0.41800,0.01950}%

\begin{tikzpicture}

\begin{axis}[%
width=3in,
height=1.25in,
at={(0.758in,0.481in)},
scale only axis,
xmode=log,
xmin=500,
xmax=200000,
xminorticks=true,
xlabel style={font=\color{white!15!black}},
xlabel={Collocation Points},
ymode=log,
ymin=1e-12,
ymax=0.01,
yminorticks=true,
ylabel style={font=\color{white!15!black}},
ylabel={Collocation Error},
axis background/.style={fill=white},
title style={font=\bfseries},
title={\textbf{Temporal Lobe: SC Error on $\mu_Q$}},
xmajorgrids,
xminorgrids,
ymajorgrids,
yminorgrids,
legend style={at={(0.03,0.03)}, anchor=south west, legend cell align=left, align=left, draw=white!15!black}
]
\addplot [color=temporal, line width=1.5pt, draw opacity=0.25]
  table[row sep=crcr]{
575	4.72008946775785e-07\\
2241	2.17376569384952e-08\\
7183	3.05832399226492e-09\\
19825	3.86722993709909e-10\\
48639	3.52604656489092e-12\\
};

\addplot [color=temporal, line width=1.5pt,draw opacity=0.5]
  table[row sep=crcr]{%
575	3.65821234141441e-06\\
2241	1.80665114880325e-06\\
7183	1.61951260332772e-06\\
19825	1.22936997172694e-06\\
48639	2.78456326039765e-08\\
};

\addplot [color=temporal, line width=1.5pt, draw opacity=0.75]
  table[row sep=crcr]{%
575	0.000115275311005639\\
2241	4.14159360321884e-05\\
7183	2.91006020125222e-05\\
19825	2.10369898712901e-05\\
48639	6.8769296229086e-07\\
};

\addplot [color=temporal, line width=1.5pt]
  table[row sep=crcr]{%
575	0.000432870457098754\\
2241	0.000130556559096943\\
7183	8.41893441843328e-05\\
19825	5.67027976273496e-05\\
48639	4.12912973838012e-06\\
};

\end{axis}

\end{tikzpicture}%

%% file: MC_Err_2_Temporal.tikz
\definecolor{Temporal}{rgb}{1.00000,0.41800,0.01950}%

\begin{tikzpicture}

\begin{axis}[%
width=3in,
height=1.25in,
at={(0.758in,0.481in)},
scale only axis,
xmode=log,
xmin=500,
xmax=200000,
xminorticks=true,
xlabel style={font=\color{white!15!black}},
xlabel={Iterations},
ymode=log,
ymin=1e-12,
ymax=0.01,
yminorticks=true,
ylabel style={font=\color{white!15!black}},
ylabel={Monte Carlo Error},
axis background/.style={fill=white},
title style={font=\bfseries},
title={\textbf{Temporal Lobe: MC Error on $\mu_Q$}},
xmajorgrids,
xminorgrids,
ymajorgrids,
yminorgrids,
legend style={at={(0.03,0.03)}, anchor=south west, legend cell align=left, align=left, draw=white!15!black}
]

\addplot [color=black, dashed, line width=1.5pt]
  table[row sep=crcr]{%
5000	0.014142135623731\\
10000	0.01\\
15000	0.00816496580927726\\
20000	0.00707106781186548\\
25000	0.00632455532033676\\
30000	0.00577350269189626\\
35000	0.00534522483824849\\
40000	0.005\\
45000	0.00471404520791032\\
50000	0.00447213595499958\\
55000	0.00426401432711221\\
60000	0.00408248290463863\\
65000	0.00392232270276368\\
70000	0.00377964473009227\\
75000	0.00365148371670111\\
80000	0.00353553390593274\\
85000	0.00342997170285018\\
90000	0.00333333333333333\\
95000	0.00324442842261525\\
100000	0.00316227766016838\\
105000	0.00308606699924184\\
110000	0.00301511344577764\\
115000	0.00294883912309794\\
120000	0.00288675134594813\\
125000	0.00282842712474619\\
130000	0.00277350098112615\\
135000	0.00272165526975909\\
140000	0.00267261241912424\\
145000	0.00262612865719445\\
150000	0.00258198889747161\\
155000	0.00254000254000381\\
160000	0.0025\\
165000	0.00246182981958665\\
170000	0.00242535625036333\\
175000	0.00239045721866879\\
180000	0.00235702260395516\\
185000	0.00232495277487639\\
190000	0.00229415733870562\\
195000	0.00226455406828919\\
200000	0.00223606797749979\\
};
\addlegendentry{$M^{-1/2}$}

\addplot [color=Temporal, line width=1.5pt, draw opacity=0.25]
  table[row sep=crcr]{%
5000	0.000284990385991925\\
10000	0.00119676689092665\\
15000	0.000735455189812716\\
20000	0.000553398514096902\\
25000	0.000551591321704584\\
30000	0.00121829824649971\\
35000	0.00092169153107541\\
40000	0.000602348194693606\\
45000	0.000180208006986245\\
50000	0.00011936784600705\\
55000	0.000157603790351536\\
60000	3.45619359207477e-05\\
65000	0.000108153639247515\\
70000	2.82728211174886e-05\\
75000	0.000302432561526139\\
80000	0.000162464253301271\\
85000	8.63145280636465e-05\\
90000	2.65880888616888e-06\\
95000	0.000132699774254718\\
100000	0.000250683645316488\\
105000	0.000243569721646015\\
110000	0.000243882995237697\\
115000	0.000290161445903693\\
120000	0.000229860408495389\\
125000	0.000228532327426848\\
130000	0.000150235776862347\\
135000	1.42490349479526e-05\\
140000	5.00196206719751e-05\\
145000	3.22269922148463e-05\\
150000	2.81547406033077e-05\\
155000	1.11353288523226e-05\\
160000	2.18058765683124e-05\\
165000	8.58511783198251e-05\\
170000	9.49239236392998e-05\\
175000	0.000105152296550114\\
180000	4.93278660188892e-05\\
185000	3.74991807253841e-05\\
190000	2.41618559912573e-05\\
195000	5.20263669751168e-05\\
200000	4.77961933460219e-05\\
};
\addlegendentry{$t=5$}

\addplot [color=Temporal, line width=1.5pt, draw opacity=0.5]
  table[row sep=crcr]{%
5000	0.000322455173489018\\
10000	0.000860205988223281\\
15000	0.000427376119199688\\
20000	0.000130416607466774\\
25000	0.0002944432670214\\
30000	0.00114246295799816\\
35000	0.000810900947675018\\
40000	0.000599744746841694\\
45000	0.000143169121475522\\
50000	0.000108347445767257\\
55000	0.000210099743627255\\
60000	7.05844397734207e-05\\
65000	0.00013625819196288\\
70000	2.67834202038158e-05\\
75000	0.000356085800900516\\
80000	0.000168158745240875\\
85000	7.30436579410655e-05\\
90000	2.27306324976045e-05\\
95000	0.000221999156746266\\
100000	0.000354368449080844\\
105000	0.000339013607170621\\
110000	0.000331064816068807\\
115000	0.000400027176669444\\
120000	0.000325100772907115\\
125000	0.000321404532615808\\
130000	0.000216081895817632\\
135000	5.58664183145701e-05\\
140000	0.000106673825956855\\
145000	9.18763055178238e-05\\
150000	9.90071648260635e-05\\
155000	4.541860220687e-05\\
160000	7.7070096982761e-05\\
165000	4.55010731871285e-05\\
170000	6.71181285510311e-05\\
175000	8.24365783101391e-05\\
180000	1.49575018341954e-05\\
185000	1.67871223617963e-05\\
190000	2.99779125017067e-05\\
195000	7.31847499205001e-06\\
200000	1.56662556091615e-05\\
};
\addlegendentry{$t=10$}

\addplot [color=Temporal, line width=1.5pt, draw opacity=0.75]
  table[row sep=crcr]{%
5000	0.000444802041759943\\
10000	0.000279233441958365\\
15000	5.39586213509462e-06\\
20000	0.000370198038677949\\
25000	7.41884623209881e-05\\
30000	0.000769271200720182\\
35000	0.000465765796142726\\
40000	0.000408379950589533\\
45000	8.16548873298023e-06\\
50000	5.05257207307508e-07\\
55000	0.000284924880627154\\
60000	0.000154415372243132\\
65000	0.000199285026983314\\
70000	7.26795773819432e-05\\
75000	0.000386876385458397\\
80000	0.000195361574729557\\
85000	9.42789280379188e-05\\
90000	4.21337564922286e-06\\
95000	0.000224420494161722\\
100000	0.000355450489325591\\
105000	0.000328368006616039\\
110000	0.000316509996388868\\
115000	0.000387278126720336\\
120000	0.000315768670681137\\
125000	0.000314562729176894\\
130000	0.000209361551769038\\
135000	5.68044309582064e-05\\
140000	0.000108356495282281\\
145000	0.000100488771486096\\
150000	0.000119726091963993\\
155000	6.84442863245448e-05\\
160000	9.13890182610061e-05\\
165000	2.32755759728618e-05\\
170000	5.41980764811386e-05\\
175000	6.96268674174606e-05\\
180000	8.18391752848818e-06\\
185000	3.42642589468348e-05\\
190000	4.64167177176527e-05\\
195000	6.11817812088698e-06\\
200000	1.25995783195931e-05\\
};
\addlegendentry{$t=15$}

\addplot [color=Temporal, line width=1.5pt]
  table[row sep=crcr]{%
5000	0.00057936838829244\\
10000	0.00017975272193922\\
15000	0.000308025377928711\\
20000	0.000714115688315034\\
25000	0.000344716503641755\\
30000	0.000447017996177481\\
35000	0.000176585271719761\\
40000	0.00023480582050633\\
45000	0.000106592515716004\\
50000	0.000102157825290673\\
55000	0.000341729101545602\\
60000	0.000223815341597592\\
65000	0.00024852842635975\\
70000	0.00011534196540426\\
75000	0.000401797410479932\\
80000	0.000224351433320059\\
85000	0.000124611480409392\\
90000	4.31796156813213e-05\\
95000	0.000196272644844255\\
100000	0.000323433162496288\\
105000	0.000284640806869264\\
110000	0.000271886838656194\\
115000	0.000336732530265108\\
120000	0.000269871061281836\\
125000	0.000272451380164245\\
130000	0.000177203983286335\\
135000	3.82186271145353e-05\\
140000	8.41823799468422e-05\\
145000	8.18444700176987e-05\\
150000	0.000110607830846352\\
155000	6.64300623000803e-05\\
160000	8.14483128238012e-05\\
165000	2.13762320894416e-05\\
170000	5.74431337523627e-05\\
175000	7.03227934631633e-05\\
180000	1.86772097271765e-05\\
185000	2.56252886846579e-05\\
190000	3.73297376755405e-05\\
195000	3.79962688878423e-06\\
200000	2.82910745253088e-05\\
};
\addlegendentry{$t=20$}

\end{axis}

\end{tikzpicture}%

%% file: Solution_3_Parietal.tikz
\definecolor{parietal}{rgb}{1.00000,1.00000,0.00000}%

\begin{tikzpicture}

\begin{axis}[%
width=3in,
height=1.25in,
at={(0.758in,0.481in)},
scale only axis,
xmin=0,
xmax=20,
xlabel style={font=\color{white!15!black}},
xlabel={$t$ [years]},
ymin=0.3,
ymax=1,
ylabel style={font=\color{white!15!black}},
ylabel={$\langle c_h\rangle_3$(t) [-]},
axis background/.style={fill=white},
title style={font=\bfseries},
title={\textbf{Parietal Lobe: solution mean and std range}},
axis x line*=bottom,
axis y line*=left
]

\addplot[area legend, draw=parietal, fill=parietal, draw opacity=0.3, fill opacity=0.3, forget plot]
table[row sep=crcr] {%
x	y\\
0	0.527519231955912\\
1	0.559007920761696\\
2	0.590576819111542\\
3	0.621276514247848\\
4	0.650852683972685\\
5	0.679110481950878\\
6	0.705915840418175\\
7	0.73119052866899\\
8	0.754903786078364\\
9	0.777062774409654\\
10	0.797703317058344\\
11	0.816881712922019\\
12	0.834667932151014\\
13	0.851140210646567\\
14	0.866380910346991\\
15	0.880473452071929\\
16	0.893500118214914\\
17	0.905540538607453\\
18	0.916670699414535\\
19	0.926962343749937\\
20	0.936482659663485\\
20	0.611213084622938\\
19	0.605262820808806\\
18	0.599152932864107\\
17	0.592899739901817\\
16	0.586525315124313\\
15	0.580058661592794\\
14	0.573537050866688\\
13	0.56700751459984\\
12	0.560528457877498\\
11	0.554171330784364\\
10	0.548022248362719\\
9	0.542183386553624\\
8	0.536773903963925\\
7	0.53193005503474\\
6	0.527804092534656\\
5	0.524561551121278\\
4	0.522376626890966\\
3	0.521425692053946\\
2	0.521879524278823\\
1	0.523895438158931\\
0	0.52751883294967\\
}--cycle;
\addplot [color=parietal, line width=2.0pt, forget plot]
  table[row sep=crcr]{%
0	0.527519032452791\\
1	0.541451679460313\\
2	0.556228171695182\\
3	0.571351103150897\\
4	0.586614655431826\\
5	0.601836016536078\\
6	0.616859966476416\\
7	0.631560291851865\\
8	0.645838845021145\\
9	0.659623080481639\\
10	0.672862782710531\\
11	0.685526521853191\\
12	0.697598195014256\\
13	0.709073862623203\\
14	0.71995898060684\\
15	0.730266056832361\\
16	0.740012716669613\\
17	0.749220139254635\\
18	0.757911816139321\\
19	0.766112582279371\\
20	0.773847872143212\\
};
\end{axis}

\end{tikzpicture}%

%% file: Coll_Err_3_Parietal.tikz
\definecolor{Parietal}{rgb}{1.00000,1.00000,0.00000}%

\begin{tikzpicture}

\begin{axis}[%
width=3in,
height=1.25in,
at={(0.758in,0.481in)},
scale only axis,
xmode=log,
xmin=500,
xmax=200000,
xminorticks=true,
xlabel style={font=\color{white!15!black}},
xlabel={Collocation Points},
ymode=log,
ymin=1e-12,
ymax=0.01,
yminorticks=true,
ylabel style={font=\color{white!15!black}},
ylabel={Collocation Error},
axis background/.style={fill=white},
title style={font=\bfseries},
title={\textbf{Parietal Lobe: SC Error on $\mu_Q$}},
xmajorgrids,
xminorgrids,
ymajorgrids,
yminorgrids,
legend style={at={(0.03,0.03)}, anchor=south west, legend cell align=left, align=left, draw=white!15!black}
]
\addplot [color=Parietal, line width=1.5pt, draw opacity=0.25]
  table[row sep=crcr]{
575	2.16324242531285e-07\\
2241	1.07959433155511e-08\\
7183	1.98876348271483e-09\\
19825	5.32740282542094e-10\\
48639	5.01086053071657e-12\\
};

\addplot [color=Parietal, line width=1.5pt,draw opacity=0.5]
  table[row sep=crcr]{%
575	2.06287807983047e-05\\
2241	2.45433148072221e-06\\
7183	8.67758806471719e-07\\
19825	3.85557173827205e-07\\
48639	9.25846657844491e-09\\
};

\addplot [color=Parietal, line width=1.5pt, draw opacity=0.75]
  table[row sep=crcr]{%
575	0.00013479836590453\\
2241	1.96203217770351e-05\\
7183	6.40781418583821e-06\\
19825	2.11762432694638e-06\\
48639	2.63180617449726e-07\\
};

\addplot [color=Parietal, line width=1.5pt]
  table[row sep=crcr]{%
575	0.000268807652038765\\
2241	2.63702897362863e-05\\
7183	1.1195703989082e-05\\
19825	1.97661730027468e-05\\
48639	2.36760577061493e-06\\
};

\end{axis}

\end{tikzpicture}%

%% file: MC_Err_3_Parietal.tikz
\definecolor{Parietal}{rgb}{1.00000,1.00000,0.00000}%

\begin{tikzpicture}

\begin{axis}[%
width=3in,
height=1.25in,
at={(0.758in,0.481in)},
scale only axis,
xmode=log,
xmin=500,
xmax=200000,
xminorticks=true,
xlabel style={font=\color{white!15!black}},
xlabel={Iterations},
ymode=log,
ymin=1e-12,
ymax=0.01,
yminorticks=true,
ylabel style={font=\color{white!15!black}},
ylabel={Monte Carlo Error},
axis background/.style={fill=white},
title style={font=\bfseries},
title={\textbf{Parietal Lobe: MC Error on $\mu_Q$}},
xmajorgrids,
xminorgrids,
ymajorgrids,
yminorgrids,
legend style={at={(0.03,0.03)}, anchor=south west, legend cell align=left, align=left, draw=white!15!black}
]

\addplot [color=black, dashed, line width=1.5pt]
  table[row sep=crcr]{%
5000	0.014142135623731\\
10000	0.01\\
15000	0.00816496580927726\\
20000	0.00707106781186548\\
25000	0.00632455532033676\\
30000	0.00577350269189626\\
35000	0.00534522483824849\\
40000	0.005\\
45000	0.00471404520791032\\
50000	0.00447213595499958\\
55000	0.00426401432711221\\
60000	0.00408248290463863\\
65000	0.00392232270276368\\
70000	0.00377964473009227\\
75000	0.00365148371670111\\
80000	0.00353553390593274\\
85000	0.00342997170285018\\
90000	0.00333333333333333\\
95000	0.00324442842261525\\
100000	0.00316227766016838\\
105000	0.00308606699924184\\
110000	0.00301511344577764\\
115000	0.00294883912309794\\
120000	0.00288675134594813\\
125000	0.00282842712474619\\
130000	0.00277350098112615\\
135000	0.00272165526975909\\
140000	0.00267261241912424\\
145000	0.00262612865719445\\
150000	0.00258198889747161\\
155000	0.00254000254000381\\
160000	0.0025\\
165000	0.00246182981958665\\
170000	0.00242535625036333\\
175000	0.00239045721866879\\
180000	0.00235702260395516\\
185000	0.00232495277487639\\
190000	0.00229415733870562\\
195000	0.00226455406828919\\
200000	0.00223606797749979\\
};
\addlegendentry{$M^{-1/2}$}

\addplot [color=Parietal, line width=1.5pt, draw opacity=0.25]
  table[row sep=crcr]{%
5000	0.00197683763568382\\
10000	0.00159294739338744\\
15000	0.000987847911837956\\
20000	0.000769482171780761\\
25000	0.00117863670329447\\
30000	0.00103134513715257\\
35000	0.00125598761413006\\
40000	0.00104876555914364\\
45000	0.000682228661781801\\
50000	0.000657605777567886\\
55000	0.000848146417151194\\
60000	0.000742587747596767\\
65000	0.000623392952278264\\
70000	0.000577152834929333\\
75000	0.000539789650325819\\
80000	0.000499783538674895\\
85000	0.000391603572800403\\
90000	0.00030860186351944\\
95000	0.00027675516768835\\
100000	0.000378141358942504\\
105000	0.000379102656652307\\
110000	0.000379577303736919\\
115000	0.000330709707816205\\
120000	0.000132274275450579\\
125000	0.000224684334840988\\
130000	0.000282711690326336\\
135000	0.000177205525901859\\
140000	0.000198534453261392\\
145000	0.000204505135663256\\
150000	0.000303962565050046\\
155000	0.000414358922829573\\
160000	0.000374062845181703\\
165000	0.000348468647678013\\
170000	0.000303605444631133\\
175000	0.000206639776684516\\
180000	0.000250587873624639\\
185000	0.000231732878609909\\
190000	0.000239501759105828\\
195000	0.00024457403606106\\
200000	0.000218243938219074\\
};
\addlegendentry{$t=5$}

\addplot [color=Parietal, line width=1.5pt, draw opacity=0.5]
  table[row sep=crcr]{%
5000	0.00285511330953371\\
10000	0.00233024469773655\\
15000	0.0013080828535725\\
20000	0.00096836911130741\\
25000	0.00157258636592976\\
30000	0.00145558464655191\\
35000	0.00171960067514237\\
40000	0.00143596885924735\\
45000	0.000897851046549795\\
50000	0.00087282225498129\\
55000	0.00114299463141549\\
60000	0.00102906424256351\\
65000	0.000863129528196261\\
70000	0.000813006875957463\\
75000	0.000724317856828432\\
80000	0.000653962948052949\\
85000	0.000510658095681417\\
90000	0.000384306393011625\\
95000	0.000368815512584538\\
100000	0.000508811559979174\\
105000	0.000515625177083176\\
110000	0.000519039251193205\\
115000	0.000462541899763271\\
120000	0.00016233953250931\\
125000	0.00030656150024497\\
130000	0.000383284545729175\\
135000	0.000241964854207708\\
140000	0.000263637720997874\\
145000	0.000258904736798677\\
150000	0.000411756344714332\\
155000	0.000565751119447684\\
160000	0.000505787521344219\\
165000	0.000471284833285642\\
170000	0.00041323056209472\\
175000	0.000286677586366097\\
180000	0.000357127897586155\\
185000	0.000337100746152805\\
190000	0.0003506298838382\\
195000	0.00035775293337466\\
200000	0.000321023354782493\\
};
\addlegendentry{$t=10$}

\addplot [color=Parietal, line width=1.5pt, draw opacity=0.75]
  table[row sep=crcr]{%
5000	0.00318597104852139\\
10000	0.00251150207062651\\
15000	0.00125503076946598\\
20000	0.000848651549373305\\
25000	0.00152861527596888\\
30000	0.00148208135220371\\
35000	0.00172368474138919\\
40000	0.00144300057719907\\
45000	0.000852107473412889\\
50000	0.000837597520466207\\
55000	0.00114823753663085\\
60000	0.00106153847801713\\
65000	0.000884613991978589\\
70000	0.000847944502347282\\
75000	0.000719858011120677\\
80000	0.00062440332239995\\
85000	0.000478360117696957\\
90000	0.000327707760381315\\
95000	0.000347223718283776\\
100000	0.00049749854291827\\
105000	0.000512669049371032\\
110000	0.000523193909798712\\
115000	0.000471681389116861\\
120000	0.000130834604182398\\
125000	0.000304251299519959\\
130000	0.000381711941807283\\
135000	0.000242931024638152\\
140000	0.00025487470195508\\
145000	0.000230585396413502\\
150000	0.000408656638937538\\
155000	0.000572135248049024\\
160000	0.000505364765884219\\
165000	0.000469733726759208\\
170000	0.000411414397505495\\
175000	0.000289393549800699\\
180000	0.000373083219216039\\
185000	0.000359833357944715\\
190000	0.000376964647421675\\
195000	0.000382957632490663\\
200000	0.000343138089495634\\
};
\addlegendentry{$t=15$}

\addplot [color=Parietal, line width=1.5pt]
  table[row sep=crcr]{%
5000	0.00327369608487251\\
10000	0.00241951540578565\\
15000	0.001064700092637\\
20000	0.000621019292336064\\
25000	0.00131797262018464\\
30000	0.00133170773558524\\
35000	0.00154171520384231\\
40000	0.00130234345589105\\
45000	0.00071341433660033\\
50000	0.000711889009585747\\
55000	0.00104770978884186\\
60000	0.000995554392631179\\
65000	0.00082336221664716\\
70000	0.000803041034108364\\
75000	0.000651782267051692\\
80000	0.000537941297595751\\
85000	0.000398873256170032\\
90000	0.000234138867113373\\
95000	0.000289256078046263\\
100000	0.000439855114623674\\
105000	0.000463352047236893\\
110000	0.000481211954764435\\
115000	0.000434800880276953\\
120000	8.54207939199204e-05\\
125000	0.000274253564416057\\
130000	0.00034612476197733\\
135000	0.000222495614073093\\
140000	0.000223095260802935\\
145000	0.000179253467204364\\
150000	0.000367969952596045\\
155000	0.000527212861511394\\
160000	0.000460638062150128\\
165000	0.000425012413698126\\
170000	0.000368887017360289\\
175000	0.00026331572835509\\
180000	0.000352281239613854\\
185000	0.0003476072206815\\
190000	0.000366080695257165\\
195000	0.000367987410362032\\
200000	0.000326940240192535\\
};
\addlegendentry{$t=20$}

\end{axis}

\end{tikzpicture}%

%% file: Solution_4_Insular.tikz
\begin{tikzpicture}

\begin{axis}[%
width=3in,
height=1.25in,
at={(0.758in,0.481in)},
scale only axis,
xmin=0,
xmax=20,
xlabel style={font=\color{white!15!black}},
xlabel={$t$ [years]},
ymin=0.3,
ymax=1,
ylabel style={font=\color{white!15!black}},
ylabel={$\langle c_h\rangle_4$(t) [-]},
axis background/.style={fill=white},
title style={font=\bfseries},
title={\textbf{Insular Lobe: solution mean and std range}},
axis x line*=bottom,
axis y line*=left
]

\addplot[area legend, draw=green, fill=green, draw opacity=0.3, fill opacity=0.3, forget plot]
table[row sep=crcr] {%
x	y\\
0	0.521029767610828\\
1	0.564383065530223\\
2	0.60690958232822\\
3	0.647156865634265\\
4	0.684748928598983\\
5	0.719483077935084\\
6	0.751300707790412\\
7	0.780250849734169\\
8	0.8064559657087\\
9	0.830083946911348\\
10	0.851326993430388\\
11	0.870386575301329\\
12	0.887463235861785\\
13	0.90275002944261\\
14	0.916428582605428\\
15	0.928666993075178\\
16	0.939618980487487\\
17	0.949423863878266\\
18	0.958207063603483\\
19	0.966080916280474\\
20	0.973145657381214\\
20	0.738515499840145\\
19	0.729943616601942\\
18	0.720871725814386\\
17	0.711282513266021\\
16	0.701163045270972\\
15	0.690506617108616\\
14	0.679315088280923\\
13	0.667601792779855\\
12	0.65539510731394\\
11	0.642742736720477\\
10	0.629716719844726\\
9	0.616419050529628\\
8	0.602987618785681\\
7	0.589601872934124\\
6	0.576487157782923\\
5	0.563916112896271\\
4	0.552204958116502\\
3	0.541702343065263\\
2	0.532769434464504\\
1	0.525752888758701\\
0	0.521029299467661\\
}--cycle;
\addplot [color=green, line width=2.0pt, forget plot]
  table[row sep=crcr]{%
0	0.521029533539244\\
1	0.545067977144462\\
2	0.569839508396362\\
3	0.594429604349764\\
4	0.618476943357743\\
5	0.641699595415678\\
6	0.663893932786668\\
7	0.684926361334146\\
8	0.70472179224719\\
9	0.723251498720488\\
10	0.740521856637557\\
11	0.756564656010903\\
12	0.771429171587862\\
13	0.785175911111232\\
14	0.797871835443175\\
15	0.809586805091897\\
16	0.82039101287923\\
17	0.830353188572143\\
18	0.839539394708935\\
19	0.848012266441208\\
20	0.855830578610679\\
};
\end{axis}

\end{tikzpicture}%

%% file: Coll_Err_4_Insular.tikz
\begin{tikzpicture}

\begin{axis}[%
width=3in,
height=1.25in,
at={(0.758in,0.481in)},
scale only axis,
xmode=log,
xmin=500,
xmax=200000,
xminorticks=true,
xlabel style={font=\color{white!15!black}},
xlabel={Collocation Points},
ymode=log,
ymin=1e-12,
ymax=0.01,
yminorticks=true,
ylabel style={font=\color{white!15!black}},
ylabel={Collocation Error},
axis background/.style={fill=white},
title style={font=\bfseries},
title={\textbf{Insular Lobe: SC Error on $\mu_Q$}},
xmajorgrids,
xminorgrids,
ymajorgrids,
yminorgrids,
legend style={at={(0.03,0.03)}, anchor=south west, legend cell align=left, align=left, draw=white!15!black}
]
\addplot [color=green, line width=1.5pt, draw opacity=0.25]
  table[row sep=crcr]{%
575	4.37906563262168e-07\\
2241	2.3895797307427e-08\\
7183	4.59665425232004e-09\\
19825	1.21387263366001e-09\\
48639	1.06120917245272e-11\\
};

\addplot [color=green, line width=1.5pt, draw opacity=0.5]
  table[row sep=crcr]{%
575	1.42203020757222e-05\\
2241	6.36952959717358e-07\\
7183	3.30182760250965e-07\\
19825	4.47566213656434e-07\\
48639	5.07416736933672e-09\\
};

\addplot [color=green, line width=1.5pt, draw opacity=0.75]
  table[row sep=crcr]{%
575	2.76407632911926e-05\\
2241	2.13711114743882e-05\\
7183	1.98482527689664e-05\\
19825	1.62991510420071e-05\\
48639	2.83501773872361e-07\\
};

\addplot [color=green, line width=1.5pt]
  table[row sep=crcr]{%
575	0.000290295622786698\\
2241	0.000115545578135959\\
7183	9.43850456316411e-05\\
19825	7.56996020439722e-05\\
48639	3.01674165951734e-06\\
};

\end{axis}

\end{tikzpicture}%

%% file: MC_Err_4_Insular.tikz
\definecolor{Insular}{rgb}{0.00000,1.00000,0.00000}%

\begin{tikzpicture}

\begin{axis}[%
width=3in,
height=1.25in,
at={(0.758in,0.481in)},
scale only axis,
xmode=log,
xmin=500,
xmax=200000,
xminorticks=true,
xlabel style={font=\color{white!15!black}},
xlabel={Iterations},
ymode=log,
ymin=1e-12,
ymax=0.01,
yminorticks=true,
ylabel style={font=\color{white!15!black}},
ylabel={Monte Carlo Error},
axis background/.style={fill=white},
title style={font=\bfseries},
title={\textbf{Insular Lobe: MC Error on $\mu_Q$}},
xmajorgrids,
xminorgrids,
ymajorgrids,
yminorgrids,
legend style={at={(0.03,0.03)}, anchor=south west, legend cell align=left, align=left, draw=white!15!black}
]

\addplot [color=black, dashed, line width=1.5pt]
  table[row sep=crcr]{%
5000	0.014142135623731\\
10000	0.01\\
15000	0.00816496580927726\\
20000	0.00707106781186548\\
25000	0.00632455532033676\\
30000	0.00577350269189626\\
35000	0.00534522483824849\\
40000	0.005\\
45000	0.00471404520791032\\
50000	0.00447213595499958\\
55000	0.00426401432711221\\
60000	0.00408248290463863\\
65000	0.00392232270276368\\
70000	0.00377964473009227\\
75000	0.00365148371670111\\
80000	0.00353553390593274\\
85000	0.00342997170285018\\
90000	0.00333333333333333\\
95000	0.00324442842261525\\
100000	0.00316227766016838\\
105000	0.00308606699924184\\
110000	0.00301511344577764\\
115000	0.00294883912309794\\
120000	0.00288675134594813\\
125000	0.00282842712474619\\
130000	0.00277350098112615\\
135000	0.00272165526975909\\
140000	0.00267261241912424\\
145000	0.00262612865719445\\
150000	0.00258198889747161\\
155000	0.00254000254000381\\
160000	0.0025\\
165000	0.00246182981958665\\
170000	0.00242535625036333\\
175000	0.00239045721866879\\
180000	0.00235702260395516\\
185000	0.00232495277487639\\
190000	0.00229415733870562\\
195000	0.00226455406828919\\
200000	0.00223606797749979\\
};
\addlegendentry{$M^{-1/2}$}

\addplot [color=Insular, line width=1.5pt, draw opacity=0.25]
  table[row sep=crcr]{%
5000	0.000917711461747211\\
10000	0.00106103547685949\\
15000	0.000171554644723715\\
20000	0.000434856482451605\\
25000	0.000452147700402113\\
30000	0.000329845660371807\\
35000	6.54952775267042e-05\\
40000	0.000156468818976771\\
45000	0.000275101486743103\\
50000	0.000590818221172302\\
55000	0.000638376412690892\\
60000	0.000710724295256694\\
65000	0.000572715676899821\\
70000	0.000541227055264477\\
75000	0.000498653391012345\\
80000	0.000585071970178989\\
85000	0.000467540612872443\\
90000	0.000376490704863067\\
95000	0.000493785125949235\\
100000	0.000690312422667602\\
105000	0.000667510106420178\\
110000	0.000693814189818752\\
115000	0.000652779827159583\\
120000	0.000677993212638025\\
125000	0.000547244490944749\\
130000	0.000552495418408835\\
135000	0.000476268589814258\\
140000	0.0005356612211476\\
145000	0.000528991420690005\\
150000	0.000504100757224487\\
155000	0.000492153734768956\\
160000	0.000449972284384963\\
165000	0.000437161248485916\\
170000	0.000357259659775472\\
175000	0.000392824012661186\\
180000	0.000378118569235762\\
185000	0.000421052269776424\\
190000	0.000405756624397097\\
195000	0.000371021834192129\\
200000	0.000330822761465854\\
};
\addlegendentry{$t=5$}

\addplot [color=Insular, line width=1.5pt, draw opacity=0.5]
  table[row sep=crcr]{%
5000	0.000959308959233486\\
10000	0.00134085690653275\\
15000	0.000154539118720344\\
20000	0.000425373705528817\\
25000	0.00058192012890323\\
30000	0.000545212356787747\\
35000	7.40965619606874e-06\\
40000	0.0002593035476213\\
45000	0.000301436812952395\\
50000	0.000656587142479467\\
55000	0.000724240327514565\\
60000	0.000856469906252425\\
65000	0.000677961654927693\\
70000	0.00065687513283744\\
75000	0.000569865077927737\\
80000	0.000685734531526404\\
85000	0.000542339118844281\\
90000	0.000465638678006812\\
95000	0.00060946921823321\\
100000	0.000850281351437568\\
105000	0.000814012388907771\\
110000	0.00084434357543882\\
115000	0.000796588875812913\\
120000	0.000803214272659323\\
125000	0.000649634712137287\\
130000	0.000673964562799916\\
135000	0.000566843227001748\\
140000	0.000646944835810226\\
145000	0.000626772505686408\\
150000	0.000603931567480701\\
155000	0.000568251115422209\\
160000	0.000529930133997145\\
165000	0.000511372776265784\\
170000	0.000414728581634069\\
175000	0.000459708457875092\\
180000	0.000444133152337485\\
185000	0.00050285509380252\\
190000	0.000482300178530419\\
195000	0.000439571236154864\\
200000	0.000381030339388645\\
};
\addlegendentry{$t=10$}

\addplot [color=Insular, line width=1.5pt, draw opacity=0.75]
  table[row sep=crcr]{%
5000	0.000643303302732877\\
10000	0.00118172566683429\\
15000	1.06024266098022e-06\\
20000	0.000235293092085771\\
25000	0.000522513588880237\\
30000	0.000560516263550663\\
35000	1.79831735114954e-05\\
40000	0.000268002915096112\\
45000	0.000229376033997291\\
50000	0.000544448177789053\\
55000	0.000632235190864782\\
60000	0.000795279771487936\\
65000	0.000624532611515413\\
70000	0.000619628882180196\\
75000	0.000507843220364909\\
80000	0.000626283494950142\\
85000	0.000478883950040713\\
90000	0.000437264263895452\\
95000	0.000583029168301575\\
100000	0.000818217279403391\\
105000	0.000779420985633843\\
110000	0.000804942705940543\\
115000	0.000767797723060039\\
120000	0.000753815046392979\\
125000	0.00061182195576877\\
130000	0.000654883319154007\\
135000	0.000537050674677899\\
140000	0.000620488282523139\\
145000	0.000595104971685718\\
150000	0.000582119403901899\\
155000	0.000529270958369317\\
160000	0.000502525095265942\\
165000	0.000481023629339092\\
170000	0.00038575732755501\\
175000	0.000432589308340948\\
180000	0.00041633121404365\\
185000	0.000481694153241086\\
190000	0.000459540645702984\\
195000	0.000421271907835696\\
200000	0.000356475946796233\\
};
\addlegendentry{$t=15$}

\addplot [color=Insular, line width=1.5pt]
  table[row sep=crcr]{%
5000	0.000291576148612651\\
10000	0.00088380936885363\\
15000	0.000188679768751444\\
20000	1.84448761900479e-05\\
25000	0.000397013412202561\\
30000	0.000477833403829147\\
35000	6.70008520359251e-05\\
40000	0.000238832633844\\
45000	0.000157429298918126\\
50000	0.000419425613825095\\
55000	0.000521267913677021\\
60000	0.000693410320085116\\
65000	0.00054760439544409\\
70000	0.000554188069026809\\
75000	0.00043114431530711\\
80000	0.000541553288880125\\
85000	0.00039275423179735\\
90000	0.000376915735170806\\
95000	0.000520374387989118\\
100000	0.00074268605106864\\
105000	0.000707528438573939\\
110000	0.000725862211482172\\
115000	0.000703075533587818\\
120000	0.000676816467959277\\
125000	0.00055176956373311\\
130000	0.000607980279074544\\
135000	0.000487326163458663\\
140000	0.000566462034269064\\
145000	0.000541410059678066\\
150000	0.000538721671863828\\
155000	0.000477466858165627\\
160000	0.000458417702019989\\
165000	0.000434248184271522\\
170000	0.000343696141626849\\
175000	0.000390074275961527\\
180000	0.000372166431750641\\
185000	0.000439373762986838\\
190000	0.000418479414089405\\
195000	0.000387363121031617\\
200000	0.000321997669245851\\
};
\addlegendentry{$t=20$}

\end{axis}

\end{tikzpicture}%

%% file: Solution_5_Limbic.tikz
\definecolor{Limbic}{rgb}{0.00000,1.00000,1.00000}%

\begin{tikzpicture}

\begin{axis}[%
width=3in,
height=1.25in,
at={(0.758in,0.481in)},
scale only axis,
xmin=0,
xmax=20,
xlabel style={font=\color{white!15!black}},
xlabel={$t$ [years]},
ymin=0.3,
ymax=1,
ylabel style={font=\color{white!15!black}},
ylabel={$\langle c_h\rangle_5$(t) [-]},
axis background/.style={fill=white},
title style={font=\bfseries},
title={\textbf{Limbic Lobe: solution mean and std range}},
axis x line*=bottom,
axis y line*=left
]

\addplot[area legend, draw=Limbic, fill=Limbic, draw opacity=0.3, fill opacity=0.3, forget plot]
table[row sep=crcr] {%
x	y\\
0	0.557175803098968\\
1	0.603637507285033\\
2	0.648438296901211\\
3	0.690101575562346\\
4	0.728304023540234\\
5	0.762912209750693\\
6	0.793958526105267\\
7	0.821600194285251\\
8	0.846075280552393\\
9	0.867664625465957\\
10	0.886663165863862\\
11	0.903360572922022\\
12	0.91802951972924\\
13	0.930919528411623\\
14	0.94225457620137\\
15	0.952233057231026\\
16	0.961029106661624\\
17	0.96879462515676\\
18	0.97566158436363\\
19	0.981744361089282\\
20	0.987141957941169\\
20	0.793793452184275\\
19	0.786263064604132\\
18	0.778209379643165\\
17	0.769599096312105\\
16	0.760399988887026\\
15	0.750582371521329\\
14	0.740121095605896\\
13	0.728998214577993\\
12	0.717206457972887\\
11	0.704753639677719\\
10	0.691668061765494\\
9	0.678004830621909\\
8	0.663852731603833\\
7	0.649340869430277\\
6	0.634643670817144\\
5	0.619982181491148\\
4	0.60561923417836\\
3	0.591846688920469\\
2	0.57896523829125\\
1	0.567260935174474\\
0	0.557174770071199\\
}--cycle;
\addplot [color=Limbic, line width=2.0pt, forget plot]
  table[row sep=crcr]{%
0	0.557175286585084\\
1	0.585449221229754\\
2	0.613701767596231\\
3	0.640974132241408\\
4	0.666961628859297\\
5	0.69144719562092\\
6	0.714301098461205\\
7	0.735470531857764\\
8	0.754964006078113\\
9	0.772834728043933\\
10	0.789165613814678\\
11	0.804057106299871\\
12	0.817617988851063\\
13	0.829958871494808\\
14	0.841187835903633\\
15	0.851407714376177\\
16	0.860714547774325\\
17	0.869196860734432\\
18	0.876935482003397\\
19	0.884003712846707\\
20	0.890467705062722\\
};
\end{axis}

\end{tikzpicture}%

%% file: Coll_Err_5_Limbic.tikz
\definecolor{Limbic}{rgb}{0.00000,1.00000,1.00000}%

\begin{tikzpicture}

\begin{axis}[%
width=3in,
height=1.25in,
at={(0.758in,0.481in)},
scale only axis,
xmode=log,
xmin=500,
xmax=200000,
xminorticks=true,
xlabel style={font=\color{white!15!black}},
xlabel={Collocation Points},
ymode=log,
ymin=1e-12,
ymax=0.01,
yminorticks=true,
ylabel style={font=\color{white!15!black}},
ylabel={Collocation Error},
axis background/.style={fill=white},
title style={font=\bfseries},
title={\textbf{Limbic Lobe: SC Error on $\mu_Q$}},
xmajorgrids,
xminorgrids,
ymajorgrids,
yminorgrids,
legend style={at={(0.03,0.03)}, anchor=south west, legend cell align=left, align=left, draw=white!15!black}
]
\addplot [color=Limbic, line width=1.5pt, draw opacity=0.25]
  table[row sep=crcr]{
575	1.85004772791545e-07\\
2241	9.26722307354239e-09\\
7183	2.10518390968155e-09\\
19825	7.50814707612665e-10\\
48639	8.60629047261825e-12\\
};

\addplot [color=Limbic, line width=1.5pt,draw opacity=0.5]
  table[row sep=crcr]{%
575	5.6511488597048e-06\\
2241	2.2790762903153e-07\\
7183	1.20392017468587e-07\\
19825	1.32462222137465e-07\\
48639	1.9859257966405e-09\\
};

\addplot [color=Limbic, line width=1.5pt, draw opacity=0.75]
  table[row sep=crcr]{%
575	6.22725317896004e-05\\
2241	2.00000085786319e-05\\
7183	1.39882317381696e-05\\
19825	9.9889568961363e-06\\
48639	3.68571481877273e-07\\
};

\addplot [color=Limbic, line width=1.5pt]
  table[row sep=crcr]{%
575	0.000290935247681523\\
2241	8.08562070505462e-05\\
7183	4.88714248624391e-05\\
19825	3.23885311557742e-05\\
48639	3.45301946824313e-06\\
};

\end{axis}

\end{tikzpicture}%

%% file: MC_Err_5_Limbic.tikz
\definecolor{Limbic}{rgb}{0.00000,1.00000,1.00000}%

\begin{tikzpicture}

\begin{axis}[%
width=3in,
height=1.25in,
at={(0.758in,0.481in)},
scale only axis,
xmode=log,
xmin=500,
xmax=200000,
xminorticks=true,
xlabel style={font=\color{white!15!black}},
xlabel={Iterations},
ymode=log,
ymin=1e-12,
ymax=0.01,
yminorticks=true,
ylabel style={font=\color{white!15!black}},
ylabel={Monte Carlo Error},
axis background/.style={fill=white},
title style={font=\bfseries},
title={\textbf{Limbic Lobe: MC Error on $\mu_Q$}},
xmajorgrids,
xminorgrids,
ymajorgrids,
yminorgrids,
legend style={at={(0.03,0.03)}, anchor=south west, legend cell align=left, align=left, draw=white!15!black}
]

\addplot [color=black, dashed, line width=1.5pt]
  table[row sep=crcr]{%
5000	0.014142135623731\\
10000	0.01\\
15000	0.00816496580927726\\
20000	0.00707106781186548\\
25000	0.00632455532033676\\
30000	0.00577350269189626\\
35000	0.00534522483824849\\
40000	0.005\\
45000	0.00471404520791032\\
50000	0.00447213595499958\\
55000	0.00426401432711221\\
60000	0.00408248290463863\\
65000	0.00392232270276368\\
70000	0.00377964473009227\\
75000	0.00365148371670111\\
80000	0.00353553390593274\\
85000	0.00342997170285018\\
90000	0.00333333333333333\\
95000	0.00324442842261525\\
100000	0.00316227766016838\\
105000	0.00308606699924184\\
110000	0.00301511344577764\\
115000	0.00294883912309794\\
120000	0.00288675134594813\\
125000	0.00282842712474619\\
130000	0.00277350098112615\\
135000	0.00272165526975909\\
140000	0.00267261241912424\\
145000	0.00262612865719445\\
150000	0.00258198889747161\\
155000	0.00254000254000381\\
160000	0.0025\\
165000	0.00246182981958665\\
170000	0.00242535625036333\\
175000	0.00239045721866879\\
180000	0.00235702260395516\\
185000	0.00232495277487639\\
190000	0.00229415733870562\\
195000	0.00226455406828919\\
200000	0.00223606797749979\\
};
\addlegendentry{$M^{-1/2}$}

\addplot [color=Limbic, line width=1.5pt, draw opacity=0.25]
  table[row sep=crcr]{%
5000	0.000479888931569532\\
10000	0.000522059643756475\\
15000	0.00054505621501869\\
20000	0.000734329120226264\\
25000	0.000261470362835434\\
30000	0.00019958743381529\\
35000	0.000672678777149784\\
40000	0.000445690538685264\\
45000	0.000530956821132073\\
50000	0.000597649580758089\\
55000	0.000394187062484617\\
60000	0.000335318204820163\\
65000	0.000194666837532347\\
70000	0.000169197190276067\\
75000	0.00028497434109102\\
80000	0.000394671818974836\\
85000	0.000520450943445031\\
90000	0.000413918761054826\\
95000	0.000382610769927273\\
100000	0.000360899632122701\\
105000	0.000478322657683575\\
110000	0.000385676507634648\\
115000	0.000416675321675494\\
120000	0.000364233679789923\\
125000	0.000442256182982448\\
130000	0.000424298218139735\\
135000	0.00033582048399235\\
140000	0.000384646979233522\\
145000	0.000394682414950231\\
150000	0.000282658450134864\\
155000	0.000263703299417807\\
160000	0.000291317906454031\\
165000	0.000232143246542625\\
170000	0.00026860261139499\\
175000	0.000213235703358344\\
180000	0.000254782344007105\\
185000	0.000301701833721357\\
190000	0.00035743240586266\\
195000	0.000432864037894604\\
200000	0.000421714189909271\\
};
\addlegendentry{$t=5$}

\addplot [color=Limbic, line width=1.5pt, draw opacity=0.5]
  table[row sep=crcr]{%
5000	0.000570871839033795\\
10000	0.000453121328051481\\
15000	0.00053160138996332\\
20000	0.000754901674176117\\
25000	0.000336282603673808\\
30000	0.000366923246693303\\
35000	0.000931022109083248\\
40000	0.000695324204333278\\
45000	0.000734895075873896\\
50000	0.000798473277279024\\
55000	0.000560147576448771\\
60000	0.000527517696887027\\
65000	0.000370456022800039\\
70000	0.000338212112101535\\
75000	0.00046246936156534\\
80000	0.000572359934477633\\
85000	0.000718904493303662\\
90000	0.00060155636290439\\
95000	0.00056222436234397\\
100000	0.000531960743006782\\
105000	0.000651767932979392\\
110000	0.000539197069450562\\
115000	0.000584828986284632\\
120000	0.000502429602573212\\
125000	0.000593123376973696\\
130000	0.000567746730252089\\
135000	0.000450467296313863\\
140000	0.000507583906763646\\
145000	0.000518875884548234\\
150000	0.000391761281686261\\
155000	0.000359937391795214\\
160000	0.000384817583195194\\
165000	0.00031319042159305\\
170000	0.000343837657859088\\
175000	0.000275344449274359\\
180000	0.000339075616456331\\
185000	0.000391968538748526\\
190000	0.000457986822946695\\
195000	0.00053700897684065\\
200000	0.00051565305841079\\
};
\addlegendentry{$t=10$}

\addplot [color=Limbic, line width=1.5pt, draw opacity=0.75]
  table[row sep=crcr]{%
5000	0.00055268292719933\\
10000	0.000197745528470904\\
15000	0.000309708015899052\\
20000	0.000525812415293239\\
25000	0.000256229703966432\\
30000	0.000395770007542733\\
35000	0.000953082628658565\\
40000	0.000791092844014595\\
45000	0.000758565857762046\\
50000	0.000805904607569456\\
55000	0.000590889225913612\\
60000	0.000579694900416749\\
65000	0.00044236852207373\\
70000	0.000400459236720251\\
75000	0.000500723516674049\\
80000	0.000590364238685685\\
85000	0.000727785357097337\\
90000	0.000635100042605983\\
95000	0.000602874468796631\\
100000	0.00057475709569053\\
105000	0.000667887625904579\\
110000	0.000563322974575385\\
115000	0.000614645831347947\\
120000	0.000520868049439201\\
125000	0.00060499230053461\\
130000	0.000578848375665101\\
135000	0.000454780532932765\\
140000	0.000504809579417818\\
145000	0.000514688734429116\\
150000	0.000405385476762768\\
155000	0.000367490927909671\\
160000	0.000382397585300529\\
165000	0.000317214052458276\\
170000	0.000330419126377666\\
175000	0.000261541020399633\\
180000	0.00033341966431585\\
185000	0.000383202962466764\\
190000	0.000446420386522383\\
195000	0.000506088308067895\\
200000	0.000477286912254771\\
};
\addlegendentry{$t=15$}

\addplot [color=Limbic, line width=1.5pt]
  table[row sep=crcr]{%
5000	0.000502577496571924\\
10000	5.34120180480826e-05\\
15000	7.45641879172239e-05\\
20000	0.000276875742079671\\
25000	0.000130685529744531\\
30000	0.000360135099537449\\
35000	0.000891438587380695\\
40000	0.000812268091915876\\
45000	0.000721672412957454\\
50000	0.000751078746487403\\
55000	0.000569173113182783\\
60000	0.000565533182375545\\
65000	0.000449532929330617\\
70000	0.000396262821520146\\
75000	0.000468739695406053\\
80000	0.000541990899679229\\
85000	0.000666692540515256\\
90000	0.00060137070991783\\
95000	0.000582329893082454\\
100000	0.000559991951763548\\
105000	0.000630532606198549\\
110000	0.00054060124755778\\
115000	0.000590703901279728\\
120000	0.000497032775668923\\
125000	0.000571985873935071\\
130000	0.000547408760757835\\
135000	0.000424221665786612\\
140000	0.000464394542531425\\
145000	0.000471013767411538\\
150000	0.000384018709620942\\
155000	0.000345377491925843\\
160000	0.000350448139731181\\
165000	0.000296227981016408\\
170000	0.000294853602956596\\
175000	0.000228787045634403\\
180000	0.000300443593280634\\
185000	0.000346800539208557\\
190000	0.000404802387080794\\
195000	0.000444304864779664\\
200000	0.000410509016828248\\
};
\addlegendentry{$t=20$}

\end{axis}

\end{tikzpicture}%

%% file: Solution_6_Occipital.tikz
\definecolor{Occipital}{rgb}{0.00000,0.00000,1.00000}%

\begin{tikzpicture}

\begin{axis}[%
width=3in,
height=1.25in,
at={(0.758in,0.481in)},
scale only axis,
xmin=0,
xmax=20,
xlabel style={font=\color{white!15!black}},
xlabel={$t$ [years]},
ymin=0.3,
ymax=1,
ylabel style={font=\color{white!15!black}},
ylabel={$\langle c_h\rangle_6$(t) [-]},
axis background/.style={fill=white},
title style={font=\bfseries},
title={\textbf{Occipital Lobe: solution mean and std range}},
axis x line*=bottom,
axis y line*=left
]

\addplot[area legend, draw=Occipital, fill=Occipital, draw opacity=0.3, fill opacity=0.3, forget plot]
table[row sep=crcr] {%
x	y\\
0	0.415903163806042\\
1	0.452909848742768\\
2	0.490410960158834\\
3	0.527270375670542\\
4	0.563112997124532\\
5	0.597607816683846\\
6	0.630493513583986\\
7	0.661588760493851\\
8	0.690789706655713\\
9	0.718059427344552\\
10	0.743413995211858\\
11	0.766908451896251\\
12	0.788624451365153\\
13	0.808660234237269\\
14	0.827122946494238\\
15	0.844123016040818\\
16	0.859770204227978\\
17	0.874170954999373\\
18	0.887426713094879\\
19	0.899632944483803\\
20	0.910878652301761\\
20	0.556299041115558\\
19	0.548007784867186\\
18	0.539442589306645\\
17	0.530620757465604\\
16	0.521567281977561\\
15	0.512316459619394\\
14	0.502913692369009\\
13	0.493417424627287\\
12	0.48390111703461\\
11	0.474455087105739\\
10	0.465187952755471\\
9	0.456227300325954\\
8	0.447719085593792\\
7	0.439825222506589\\
6	0.432718930413\\
5	0.426577841385001\\
4	0.421575694568172\\
3	0.417874465926188\\
2	0.415619295452675\\
1	0.414937507449581\\
0	0.415902114887066\\
}--cycle;
\addplot [color=Occipital, line width=2.0pt, forget plot]
  table[row sep=crcr]{%
0	0.415902639346554\\
1	0.433923678096175\\
2	0.453015127805755\\
3	0.472572420798365\\
4	0.492344345846352\\
5	0.512092829034424\\
6	0.531606221998493\\
7	0.55070699150022\\
8	0.569254396124753\\
9	0.587143363835253\\
10	0.604300973983664\\
11	0.620681769500995\\
12	0.636262784199881\\
13	0.651038829432278\\
14	0.665018319431623\\
15	0.678219737830106\\
16	0.69066874310277\\
17	0.702395856232489\\
18	0.713434651200762\\
19	0.723820364675495\\
20	0.73358884670866\\
};
\end{axis}

\end{tikzpicture}%

%% file: Coll_Err_6_Occipital.tikz
\begin{tikzpicture}

\begin{axis}[%
width=3in,
height=1.25in,
at={(0.758in,0.481in)},
scale only axis,
xmode=log,
xmin=500,
xmax=200000,
xminorticks=true,
xlabel style={font=\color{white!15!black}},
xlabel={Collocation Points},
ymode=log,
ymin=1e-12,
ymax=0.01,
yminorticks=true,
ylabel style={font=\color{white!15!black}},
ylabel={Collocation Error},
axis background/.style={fill=white},
title style={font=\bfseries},
title={\textbf{Occipital Lobe: SC Error on $\mu_Q$}},
xmajorgrids,
xminorgrids,
ymajorgrids,
yminorgrids,
legend style={at={(0.03,0.03)}, anchor=south west, legend cell align=left, align=left, draw=white!15!black}
]
\addplot [color=blue, line width=1.5pt, draw opacity=0.25]
  table[row sep=crcr]{%
575	2.79361727601455e-07\\
2241	2.14471812094916e-08\\
7183	5.30100522908445e-09\\
19825	1.84537851523169e-09\\
48639	2.3474647703164e-11\\
};

\addplot [color=blue, line width=1.5pt, draw opacity=0.5]
  table[row sep=crcr]{%
575	5.99989682679536e-06\\
2241	3.6593111620183e-07\\
7183	4.09083867434639e-07\\
19825	4.2773130867943e-07\\
48639	4.20571669924371e-08\\
};

\addplot [color=blue, line width=1.5pt, draw opacity=0.75]
  table[row sep=crcr]{%
575	0.000252472460013181\\
2241	5.92233755310807e-05\\
7183	3.79711754427075e-05\\
19825	2.74350390951002e-05\\
48639	4.54508696579774e-07\\
};

\addplot [color=blue, line width=1.5pt]
  table[row sep=crcr]{%
575	0.000993440051995627\\
2241	0.000301839155715709\\
7183	0.000205378789427514\\
19825	0.00016092033850976\\
48639	8.1611203606151e-06\\
};

\end{axis}

\end{tikzpicture}%

%% file: MC_Err_6_Occipital.tikz
\definecolor{Occipital}{rgb}{0.00000,0.00000,1.00000}%

\begin{tikzpicture}

\begin{axis}[%
width=3in,
height=1.25in,
at={(0.758in,0.481in)},
scale only axis,
xmode=log,
xmin=500,
xmax=200000,
xminorticks=true,
xlabel style={font=\color{white!15!black}},
xlabel={Iterations},
ymode=log,
ymin=1e-12,
ymax=0.01,
yminorticks=true,
ylabel style={font=\color{white!15!black}},
ylabel={Monte Carlo Error},
axis background/.style={fill=white},
title style={font=\bfseries},
title={\textbf{Occipital Lobe: MC Error on $\mu_Q$}},
xmajorgrids,
xminorgrids,
ymajorgrids,
yminorgrids,
legend style={at={(0.03,0.03)}, anchor=south west, legend cell align=left, align=left, draw=white!15!black}
]

\addplot [color=black, dashed, line width=1.5pt]
  table[row sep=crcr]{%
5000	0.014142135623731\\
10000	0.01\\
15000	0.00816496580927726\\
20000	0.00707106781186548\\
25000	0.00632455532033676\\
30000	0.00577350269189626\\
35000	0.00534522483824849\\
40000	0.005\\
45000	0.00471404520791032\\
50000	0.00447213595499958\\
55000	0.00426401432711221\\
60000	0.00408248290463863\\
65000	0.00392232270276368\\
70000	0.00377964473009227\\
75000	0.00365148371670111\\
80000	0.00353553390593274\\
85000	0.00342997170285018\\
90000	0.00333333333333333\\
95000	0.00324442842261525\\
100000	0.00316227766016838\\
105000	0.00308606699924184\\
110000	0.00301511344577764\\
115000	0.00294883912309794\\
120000	0.00288675134594813\\
125000	0.00282842712474619\\
130000	0.00277350098112615\\
135000	0.00272165526975909\\
140000	0.00267261241912424\\
145000	0.00262612865719445\\
150000	0.00258198889747161\\
155000	0.00254000254000381\\
160000	0.0025\\
165000	0.00246182981958665\\
170000	0.00242535625036333\\
175000	0.00239045721866879\\
180000	0.00235702260395516\\
185000	0.00232495277487639\\
190000	0.00229415733870562\\
195000	0.00226455406828919\\
200000	0.00223606797749979\\
};
\addlegendentry{$M^{-1/2}$}

\addplot [color=Occipital, line width=1.5pt, draw opacity=0.25]
  table[row sep=crcr]{%
5000	0.0030860287128344\\
10000	0.00319572117114928\\
15000	0.00277322295269794\\
20000	0.00220072451037441\\
25000	0.00200594581436005\\
30000	0.0015748148132551\\
35000	0.00159497947495555\\
40000	0.000923022677043802\\
45000	0.000479801961857857\\
50000	0.000336570578188696\\
55000	0.000243797275105456\\
60000	9.88024592801018e-05\\
65000	0.000195149875329068\\
70000	0.000270670515994229\\
75000	5.6177015058946e-05\\
80000	0.000132234322000289\\
85000	0.000106859866185306\\
90000	0.000103547146190129\\
95000	0.000226726657253245\\
100000	3.21919734449819e-05\\
105000	9.92193676485133e-05\\
110000	0.000115476716564423\\
115000	0.000275896691826741\\
120000	0.000155672953935533\\
125000	0.000156396422573493\\
130000	0.000107474619985069\\
135000	8.69365734522599e-05\\
140000	7.59399380714226e-05\\
145000	0.000126179873942434\\
150000	0.000281195169144492\\
155000	0.000238376472276869\\
160000	0.000189437839260643\\
165000	0.000304389865885955\\
170000	0.00040202211884926\\
175000	0.000423022530800934\\
180000	0.000412795364810181\\
185000	0.000496445186129565\\
190000	0.000536620607460574\\
195000	0.000491434290717423\\
200000	0.000375719973392326\\
};
\addlegendentry{$t=5$}

\addplot [color=Occipital, line width=1.5pt, draw opacity=0.5]
  table[row sep=crcr]{%
5000	0.00397021548460941\\
10000	0.00402754770928975\\
15000	0.00371201516595272\\
20000	0.00305358929137513\\
25000	0.00272597342022199\\
30000	0.00208256875826922\\
35000	0.0021129776030821\\
40000	0.00120857174734915\\
45000	0.000684952125816505\\
50000	0.000507253086280988\\
55000	0.000403702382555672\\
60000	0.00017765918505392\\
65000	0.00023508570123163\\
70000	0.000342143483503814\\
75000	2.92206467877757e-05\\
80000	0.000159944332990684\\
85000	0.000140503919929575\\
90000	0.0001423133097664\\
95000	0.00033688302262586\\
100000	9.71480556608328e-05\\
105000	0.000171917608279725\\
110000	0.000189956185851273\\
115000	0.00041852008576097\\
120000	0.000244791300388303\\
125000	0.000250296153699805\\
130000	0.000173653262914861\\
135000	0.000114047019709718\\
140000	9.63263182071592e-05\\
145000	0.000157601558427492\\
150000	0.000374205170257607\\
155000	0.000310763056849082\\
160000	0.00026001737884043\\
165000	0.000413352474540492\\
170000	0.000536607370090014\\
175000	0.000570819952901587\\
180000	0.000567880574097019\\
185000	0.000684098712358575\\
190000	0.000736021180679908\\
195000	0.00066542487014383\\
200000	0.000504974892391964\\
};
\addlegendentry{$t=10$}

\addplot [color=Occipital, line width=1.5pt, draw opacity=0.75]
  table[row sep=crcr]{%
5000	0.00388357071519577\\
10000	0.00393314982346176\\
15000	0.00390316632360414\\
20000	0.0033003834177487\\
25000	0.0028672537730317\\
30000	0.00215117933950835\\
35000	0.00217580107770396\\
40000	0.00122203211038968\\
45000	0.000737051143088899\\
50000	0.000568450185123762\\
55000	0.000485349618905516\\
60000	0.000215992525877731\\
65000	0.000235988510747673\\
70000	0.000339262496172173\\
75000	2.70275214669603e-06\\
80000	0.000158903516031845\\
85000	0.000154870068740394\\
90000	0.000160046695108377\\
95000	0.000403940532054711\\
100000	0.000183783993292505\\
105000	0.000249794547093675\\
110000	0.000258539883845795\\
115000	0.000510980776316472\\
120000	0.000320796785717539\\
125000	0.000337712533787791\\
130000	0.00024114854362346\\
135000	0.000144533135820208\\
140000	0.000115609402161505\\
145000	0.000173646101543483\\
150000	0.00040665268124292\\
155000	0.000334956155831862\\
160000	0.000299432163085075\\
165000	0.000465862405288091\\
170000	0.000578055959082517\\
175000	0.000621181356113464\\
180000	0.000623880720824055\\
185000	0.000745662304243699\\
190000	0.000796604787394959\\
195000	0.000710081233062654\\
200000	0.000536233165763703\\
};
\addlegendentry{$t=15$}

\addplot [color=Occipital, line width=1.5pt]
  table[row sep=crcr]{%
5000	0.00347616529932544\\
10000	0.00359122202316272\\
15000	0.00381535339708815\\
20000	0.00326879340856351\\
25000	0.00275906753050748\\
30000	0.00202827936626369\\
35000	0.00204983986363724\\
40000	0.00110973250463802\\
45000	0.000691787748445023\\
50000	0.000544762992287339\\
55000	0.000487952058743722\\
60000	0.000201497588941194\\
65000	0.000255029133194067\\
70000	0.000341397717062646\\
75000	3.93163995729075e-07\\
80000	0.000172413653694571\\
85000	0.000183225371964195\\
90000	0.000189061711237169\\
95000	0.000464487781605397\\
100000	0.000279432438222532\\
105000	0.000337030685989949\\
110000	0.000333491391810939\\
115000	0.000589888125952688\\
120000	0.000391743482343757\\
125000	0.000425490733007478\\
130000	0.000314152828384099\\
135000	0.000190640209561255\\
140000	0.000153609816180637\\
145000	0.00020504736380866\\
150000	0.000437428175861566\\
155000	0.000366204239208681\\
160000	0.000345610825524732\\
165000	0.000514775179553734\\
170000	0.000606022192531341\\
175000	0.000651839125206333\\
180000	0.00065436439624316\\
185000	0.000770744117741402\\
190000	0.000816331453589977\\
195000	0.000721112949096778\\
200000	0.00054627327358769\\
};
\addlegendentry{$t=20$}

\end{axis}

\end{tikzpicture}%

%% file: Solution_7_Subcortical.tikz
\definecolor{Subcortical}{rgb}{1.00000,0.00000,1.00000}%

\begin{tikzpicture}

\begin{axis}[%
width=3in,
height=1.25in,
at={(0.758in,0.481in)},
scale only axis,
xmin=0,
xmax=20,
xlabel style={font=\color{white!15!black}},
xlabel={$t$ [years]},
ymin=0.3,
ymax=1,
ylabel style={font=\color{white!15!black}},
ylabel={$\langle c_h\rangle_7$(t) [-]},
axis background/.style={fill=white},
title style={font=\bfseries},
title={\textbf{Subcortical Nuclei: solution mean and std range}},
axis x line*=bottom,
axis y line*=left
]

\addplot[area legend, draw=Subcortical, fill=Subcortical, draw opacity=0.3, fill opacity=0.3, forget plot]
table[row sep=crcr] {%
x	y\\
0	0.419650232382519\\
1	0.473044670875041\\
2	0.525307899111023\\
3	0.574769690529219\\
4	0.621082694286648\\
5	0.664045332453749\\
6	0.703557891596876\\
7	0.739604823349109\\
8	0.772249074906168\\
9	0.801624709997835\\
10	0.827922990544585\\
11	0.851374120857673\\
12	0.872228789284263\\
13	0.890742512546124\\
14	0.907164024277087\\
15	0.921727661377256\\
16	0.934649059193585\\
17	0.946123292574949\\
18	0.956324675198661\\
19	0.96540759561453\\
20	0.973507939670836\\
20	0.70906859575619\\
19	0.699871347388667\\
18	0.690039833614725\\
17	0.679536642480313\\
16	0.668327147698638\\
15	0.656381677506701\\
14	0.643678304414089\\
13	0.630206296919317\\
12	0.615970183732398\\
11	0.600994222643337\\
10	0.585326834074065\\
9	0.569044283944424\\
8	0.552252682140221\\
7	0.535087400673347\\
6	0.517709549744932\\
5	0.500300131399354\\
4	0.483053001659074\\
3	0.466166222334374\\
2	0.449828065181516\\
1	0.434194074459581\\
0	0.419649492933508\\
}--cycle;
\addplot [color=Subcortical, line width=2.0pt, forget plot]
  table[row sep=crcr]{%
0	0.419649862658013\\
1	0.453619372667311\\
2	0.48756798214627\\
3	0.520467956431796\\
4	0.552067847972861\\
5	0.582172731926552\\
6	0.610633720670904\\
7	0.637346112011228\\
8	0.662250878523195\\
9	0.68533449697113\\
10	0.706624912309325\\
11	0.726184171750505\\
12	0.74409948650833\\
13	0.760474404732721\\
14	0.775421164345588\\
15	0.789054669441978\\
16	0.801488103446112\\
17	0.812829967527631\\
18	0.823182254406693\\
19	0.832639471501598\\
20	0.841288267713513\\
};
\end{axis}

\end{tikzpicture}%

%% file: Coll_Err_7_Subcortical.tikz
\definecolor{Subcortical}{rgb}{1.00000,0.00000,1.00000}%

\begin{tikzpicture}

\begin{axis}[%
width=3in,
height=1.25in,
at={(0.758in,0.481in)},
scale only axis,
xmode=log,
xmin=500,
xmax=200000,
xminorticks=true,
xlabel style={font=\color{white!15!black}},
xlabel={Collocation Points},
ymode=log,
ymin=1e-12,
ymax=0.01,
yminorticks=true,
ylabel style={font=\color{white!15!black}},
ylabel={Collocation Error},
axis background/.style={fill=white},
title style={font=\bfseries},
title={\textbf{Subcortical Nuclei: SC Error on $\mu_Q$}},
xmajorgrids,
xminorgrids,
ymajorgrids,
yminorgrids,
legend style={at={(0.03,0.03)}, anchor=south west, legend cell align=left, align=left, draw=white!15!black}
]
\addplot [color=Subcortical, line width=1.5pt, draw opacity=0.25]
  table[row sep=crcr]{
575	6.77306627330354e-07\\
2241	3.97717678617187e-08\\
7183	8.16224443732683e-09\\
19825	2.25864221546589e-09\\
48639	2.11669073229324e-11\\
};

\addplot [color=Subcortical, line width=1.5pt,draw opacity=0.5]
  table[row sep=crcr]{%
575	1.57678422496981e-05\\
2241	4.93348405850874e-07\\
7183	6.79910302675559e-07\\
19825	7.81899904436314e-07\\
48639	2.99233249119021e-09\\
};

\addplot [color=Subcortical, line width=1.5pt, draw opacity=0.75]
  table[row sep=crcr]{%
575	4.28139550542434e-05\\
2241	2.63244302046307e-05\\
7183	2.52505228008759e-05\\
19825	2.08710511525061e-05\\
48639	1.37515627531131e-07\\
};

\addplot [color=Subcortical, line width=1.5pt]
  table[row sep=crcr]{%
575	0.000322790501943365\\
2241	0.000124934926643746\\
7183	0.000105072892407313\\
19825	8.50633299952829e-05\\
48639	3.47830249872901e-06\\
};

\end{axis}

\end{tikzpicture}%

%% file: MC_Err_7_Subcortical.tikz
\definecolor{Subcortical}{rgb}{1.00000,0.00000,1.00000}%

\begin{tikzpicture}

\begin{axis}[%
width=3in,
height=1.25in,
at={(0.758in,0.481in)},
scale only axis,
xmode=log,
xmin=500,
xmax=200000,
xminorticks=true,
xlabel style={font=\color{white!15!black}},
xlabel={Iterations},
ymode=log,
ymin=1e-12,
ymax=0.01,
yminorticks=true,
ylabel style={font=\color{white!15!black}},
ylabel={Monte Carlo Error},
axis background/.style={fill=white},
title style={font=\bfseries},
title={\textbf{Subcortical Nuclei: MC Error on $\mu_Q$}},
xmajorgrids,
xminorgrids,
ymajorgrids,
yminorgrids,
legend style={at={(0.03,0.03)}, anchor=south west, legend cell align=left, align=left, draw=white!15!black}
]

\addplot [color=black, dashed, line width=1.5pt]
  table[row sep=crcr]{%
5000	0.014142135623731\\
10000	0.01\\
15000	0.00816496580927726\\
20000	0.00707106781186548\\
25000	0.00632455532033676\\
30000	0.00577350269189626\\
35000	0.00534522483824849\\
40000	0.005\\
45000	0.00471404520791032\\
50000	0.00447213595499958\\
55000	0.00426401432711221\\
60000	0.00408248290463863\\
65000	0.00392232270276368\\
70000	0.00377964473009227\\
75000	0.00365148371670111\\
80000	0.00353553390593274\\
85000	0.00342997170285018\\
90000	0.00333333333333333\\
95000	0.00324442842261525\\
100000	0.00316227766016838\\
105000	0.00308606699924184\\
110000	0.00301511344577764\\
115000	0.00294883912309794\\
120000	0.00288675134594813\\
125000	0.00282842712474619\\
130000	0.00277350098112615\\
135000	0.00272165526975909\\
140000	0.00267261241912424\\
145000	0.00262612865719445\\
150000	0.00258198889747161\\
155000	0.00254000254000381\\
160000	0.0025\\
165000	0.00246182981958665\\
170000	0.00242535625036333\\
175000	0.00239045721866879\\
180000	0.00235702260395516\\
185000	0.00232495277487639\\
190000	0.00229415733870562\\
195000	0.00226455406828919\\
200000	0.00223606797749979\\
};
\addlegendentry{$M^{-1/2}$}

\addplot [color=Subcortical, line width=1.5pt, draw opacity=0.25]
  table[row sep=crcr]{%
5000	0.00318306042823864\\
10000	2.94866014735784e-05\\
15000	0.000272760444403744\\
20000	7.06023183869759e-05\\
25000	0.000584636592143622\\
30000	0.00122030097552876\\
35000	0.000712081340529895\\
40000	0.0010603466312128\\
45000	0.000719601547937644\\
50000	0.000665849821516663\\
55000	0.000592553643172802\\
60000	0.000570864900323932\\
65000	0.000536082111437502\\
70000	0.000761399196171005\\
75000	0.000700420634127163\\
80000	0.000809726933874318\\
85000	0.000750226341797141\\
90000	0.000730337552643805\\
95000	0.000600398771803836\\
100000	0.000625191090985845\\
105000	0.000576769734478725\\
110000	0.000502056451115956\\
115000	0.0006039197469909\\
120000	0.000506511251816582\\
125000	0.000428803597894591\\
130000	0.000478187943839482\\
135000	0.000367085928690155\\
140000	0.000559615533859756\\
145000	0.000534479539957286\\
150000	0.000448116871886105\\
155000	0.000469705876948813\\
160000	0.00046011263239004\\
165000	0.000360104315518713\\
170000	0.000421627873909167\\
175000	0.000427158189884594\\
180000	0.000527010537179469\\
185000	0.000549916419062707\\
190000	0.000502072672441568\\
195000	0.000476250527713868\\
200000	0.000473383889148703\\
};
\addlegendentry{$t=5$}

\addplot [color=Subcortical, line width=1.5pt, draw opacity=0.5]
  table[row sep=crcr]{%
5000	0.00370761925104424\\
10000	6.29549211577162e-05\\
15000	9.13425600955055e-05\\
20000	0.000306522327153625\\
25000	0.000572057867856745\\
30000	0.00138032314098616\\
35000	0.00075447142551022\\
40000	0.00123006918274917\\
45000	0.00079239652942404\\
50000	0.000726277182144645\\
55000	0.000626560198318962\\
60000	0.000582790520114116\\
65000	0.00053275252162059\\
70000	0.000859376414261378\\
75000	0.000780706849621515\\
80000	0.00092806773805307\\
85000	0.000876279541996648\\
90000	0.000834237972758605\\
95000	0.000697047125749425\\
100000	0.000746847624550048\\
105000	0.000707882624520888\\
110000	0.000606463634710137\\
115000	0.000746874685908941\\
120000	0.000626988365493997\\
125000	0.000546505389680791\\
130000	0.000591114946161316\\
135000	0.0004410934521988\\
140000	0.000671098813016146\\
145000	0.000647488241625937\\
150000	0.000550508544621629\\
155000	0.000573256776976606\\
160000	0.00056104993310997\\
165000	0.000430610633337506\\
170000	0.000505606844499711\\
175000	0.000491688117474447\\
180000	0.000634098474762749\\
185000	0.000671066415927012\\
190000	0.000618254257983668\\
195000	0.000588706065495703\\
200000	0.000579791523106745\\
};
\addlegendentry{$t=10$}

\addplot [color=Subcortical, line width=1.5pt, draw opacity=0.75]
  table[row sep=crcr]{%
5000	0.00335821120360668\\
10000	6.43227530167325e-05\\
15000	0.00018767200841317\\
20000	0.00058661334352754\\
25000	0.000374742567135274\\
30000	0.00118162332778905\\
35000	0.000564691627723516\\
40000	0.00109096134488863\\
45000	0.000662055212819616\\
50000	0.000586377055288138\\
55000	0.000507723281315446\\
60000	0.000442844375218512\\
65000	0.000383666308860415\\
70000	0.000745113779228958\\
75000	0.000681601584837217\\
80000	0.000833777852161212\\
85000	0.000801545922001757\\
90000	0.000743970022606379\\
95000	0.000642543696282369\\
100000	0.000709628476816856\\
105000	0.000689305508401655\\
110000	0.000573468766468237\\
115000	0.000720085931778277\\
120000	0.000605200835361099\\
125000	0.000541336871332267\\
130000	0.000568972655269157\\
135000	0.00041173850031957\\
140000	0.000633211489102692\\
145000	0.000617037266063714\\
150000	0.00053023333793717\\
155000	0.000546504773861907\\
160000	0.000532260970527247\\
165000	0.000404534297024869\\
170000	0.000474547657337978\\
175000	0.000441342972834794\\
180000	0.000597551915065658\\
185000	0.000644131508041084\\
190000	0.00060053650166256\\
195000	0.000577713778325319\\
200000	0.000561033665125667\\
};
\addlegendentry{$t=15$}

\addplot [color=Subcortical, line width=1.5pt]
  table[row sep=crcr]{%
5000	0.00282800116955196\\
10000	6.52366779935012e-05\\
15000	0.000398559669453093\\
20000	0.000809000167301061\\
25000	0.000182901584680915\\
30000	0.000957427391952797\\
35000	0.000370117531339558\\
40000	0.000908368545176008\\
45000	0.000516305350923262\\
50000	0.000429872101054624\\
55000	0.000383913386515289\\
60000	0.000311534557419955\\
65000	0.000248647733468527\\
70000	0.000605032785344079\\
75000	0.000562242017510966\\
80000	0.000708318837103042\\
85000	0.000695849188148357\\
90000	0.000631816301487977\\
95000	0.000573472093917095\\
100000	0.000649362758692565\\
105000	0.000641225147582572\\
110000	0.000517764617732515\\
115000	0.000652711080758452\\
120000	0.000547673189194325\\
125000	0.000500545774203327\\
130000	0.000514552689749582\\
135000	0.00036047417009635\\
140000	0.000565241752401647\\
145000	0.000554100454414959\\
150000	0.000479012050655006\\
155000	0.000486671902352708\\
160000	0.000471467128463339\\
165000	0.000357793017309085\\
170000	0.000418135014672105\\
175000	0.000375160219173383\\
180000	0.000529138710194233\\
185000	0.000581653679104866\\
190000	0.000549379663571949\\
195000	0.00053570746277207\\
200000	0.000512500067945146\\
};
\addlegendentry{$t=20$}

\end{axis}

\end{tikzpicture}%